\documentclass[ijoo,nonblindrev]{informs-ijoo}

\OneAndAHalfSpacedXI
\usepackage{amsfonts}
\usepackage{graphicx}
\usepackage{epstopdf}
\usepackage{algorithmic}
\usepackage{amsmath}
\usepackage[ampersand]{easylist}
\usepackage{mathtools}
\usepackage{accents}
\usepackage{tabu}
\usepackage{amssymb}
\usepackage{bbm}
\usepackage{dsfont}
\usepackage[numbers]{natbib}
\usepackage{multicol}
\usepackage{enumerate}
\usepackage{enumitem}
\usepackage{hyperref}
\usepackage{tikz}
\usepackage{pbox}

\newcommand{\mb}{\mathbf}
\newcommand{\bs}{\boldsymbol}
\newcommand{\mU}{\mathcal{U}}

\newcommand{\cu}{CU}
\newcommand{\cus}{CU sets}
\newcommand{\sigp}{\mb{M}}
\newcommand{\tp}{\widetilde{p}}
\newcommand{\tuq}{\widetilde{\mathbf{q}}^u}
\newcommand{\tlq}{\widetilde{\mathbf{q}}^l}
\newcommand{\uq}{\mathbf{q}^u}
\renewcommand{\lq}{\mathbf{q}^l}
\newcommand{\tR}{\widetilde{\mathbf{R}}}

\usepackage{natbib}
 \bibpunct[, ]{(}{)}{,}{a}{}{,}%
 \def\bibfont{\small}%
\TheoremsNumberedThrough     % Preferred (Theorem 1, Lemma 1, Theorem 2)
\ECRepeatTheorems

\EquationsNumberedThrough    % Default: (1), (2), ...
\MANUSCRIPTNO{}

\begin{document}
%%%%%%%%%%%%%%%%

%\setcounter{page}{1}
%\DOI{10.1287/xxxx.0000.0000}%

\RUNAUTHOR{Nohadani and Sharma}

\RUNTITLE{Optimization under Connected Uncertainty}

% \TITLE{Bundling Information Goods of Decreasing Value}
\TITLE{Optimization under Connected Uncertainty}

\ARTICLEAUTHORS{%
\AUTHOR{Omid Nohadani}
\AFF{Benefits Science Technology, Boston, MA 02110,\\ \EMAIL{onohadani@gmail.com}} %, \URL{}}
\AUTHOR{Kartikey Sharma}
\AFF{Zuse Institute Berlin, Berlin,\\
 \EMAIL{kartikey.sharma@zib.de}}
} % end of the block

\ABSTRACT{%
Robust optimization methods have shown practical advantages in a wide range of decision-making applications under uncertainty.
Recently, their efficacy has been extended to multi-period settings.
Current approaches model uncertainty either independent of the past or in an \emph{implicit} fashion by budgeting the aggregate uncertainty. 
In many applications, however, past realizations directly influence future uncertainties.
For this class of problems, we develop a modeling framework that \emph{explicitly} incorporates this dependence via \emph{connected uncertainty} sets, whose parameters at each period depend on previous uncertainty realizations.
To find optimal here-and-now solutions, we reformulate robust and distributionally robust constraints for popular set structures and demonstrate this modeling framework numerically on broadly applicable knapsack and portfolio-optimization problems.
}%

 \HISTORY{}

\maketitle
%%%%%%%%%%%%%%%%%%%%%%%%%%%%%%%%%%%%%%%%%%%%%%%%%%%%%%%%%%%%%%%%%%%%%%

\section{Introduction}
\label{sec:introduction}

Multiperiod decision making under adversarial uncertainty has instrumentally benefited a growing number of real-world applications, such as energy~\citep{jiang2012robust} and routing~\citep{agra2018robust}, amongst others. 
Current optimization methods typically model adversarial uncertainties to be independent across periods.
In many practical settings, however, the realized uncertainty in a period can affect subsequent uncertainty realizations.
Such connections have been addressed for unit commitment~\citep{lorca2015adaptive} and inventory control problems~\citep{xin2015distributionally}.
This work studies a general framework to directly model connected uncertainties, in particular, when they are based on time series.
To this end, we develop modeling frameworks for both robust and distributionally robust optimization.

The method of Robust Optimization (RO) has proven  capable of providing computationally tractable solutions for problems under uncertainty~\citep{ben2009robust,bertsimas2011theory,gabrel2014recent}. 
The paradigm of RO posits a \emph{robust optimal solution} that is feasible for all realizations of the uncertainty by
replacing the original constraints with their \emph{robust counterpart}. 
From a practical lens, the advantage of RO is that it does not require detailed information on the uncertainty and rather describes its structure by sets.
The geometry of these uncertainty sets determines the computational tractability of the formulation. 
For example, certain combinatorial RO problems achieve a tractable reformulation when the uncertain objective coefficients reside in a cardinality constrained set~\citep{bertsimas2004price, agra2016dynamic}.
The size of such sets controls the magnitude of possible uncertainties, to which the solution is immune. 
It also establishes probabilistic guarantees of constraint satisfaction.
RO techniques have found broad applications in healthcare~\citep{chu2005robust}, inventory management~\citep{aharon2009robust}, and statistics~\citep{nohadani2018}, amongst others.

When the adversarially uncertain components exhibit probabilistic characteristics, Distributionally Robust Optimization (DRO) offers an alternative approach by replacing uncertainty sets with ambiguity sets over distributions~\citep{bertsimas2010models,delage2010distributionally,goh2010distributionally}. 
These sets can be characterized by moments~\citep{wiesemann2014distributionally}, distance measures~\citep{ben2013robust,gao2016distributionally,esfahani2018data}, or hypothesis tests~\citep{bertsimas2014robust}.
DRO techniques have been applied to a broad range of applications, such as portfolio management~\citep{natarajan2010tractable}, simulation~\citep{lam2016robust}, and supply chain~\citep{fu2018profit} problems.

Most RO approaches focus on static \emph{here-and-now} solutions. 
However, the advantage of these solutions is limited in many settings that can accommodate adaptation.
To accurately capture the \emph{wait-and-see} nature, ~\citet{ben2004adjustable} introduced adjustable robust optimization (ARO), allowing decisions to adapt to observations.
This approach improves the solution quality at the expense of higher computational complexity.
In multistage approaches, the uncertain parameters are typically modeled to reside in separate sets in each stage, which enables tractability~\citep{iyengar2005robust,nilim2005robust,georghiou2016robust}. 
Uncertainties can also be modeled to depend on decisions~\citep{poss2013robust,nohadani2016optimization} or on time~\citep{arka2016}.

In many applications, however, the uncertainty at a time period depends on the realization of
the uncertain parameters in previous periods.
Autoregression is a popular method for modeling such dependence, 
and such a model can be expressed as 
\begin{equation}
\label{eq:basic_ar_mode}
\bs{\mu}_{t+1}(\mb{d}_{t}) = \mb{F}_{t} \mb{d}_{t} + \mb{c}_{t},
\end{equation}
where $\bs{\mu}_{t+1}(\mb{d}_{t})$ represents the mean of period $t+1$, $\mb{c}_{t}$ is a constant, and the matrix $\mb{F}_{t} $ controls the dependence of the mean on the past. 
For example, in queuing systems, \citet{livny1993impact} study the impact of correlations in inter-arrival or service times.
These times may be positively or negatively correlated $(\mb{F}_t \geq 0 \text{ or } \leq 0)$--that is, the observations will be in the form of waves of high and low inter-arrival or service times. 
They observe significant variations in the waiting times, depending on the type and magnitude of autocorrelation present and recommend sensitivity studies for the same. 
Similarly, \citet{balvers1997autocorrelated} study the intertemporal optimal portfolio choice problem under the assumption of autocorrelated stock returns and observe significant variations in the \emph{age effect} when auto correlation is present.
In newsvendor problems, \citet{alwan2016dynamic} provide models with autoregressive demand for commodities.
They use $\mb{F}_t$ to model correlations between period-to-period demand dynamic ordering decision based on mean square error forecasts that outperform static decisions when a sizable correlation is present. 
These problems leverage the connection amongst uncertain parameters and would benefit from the RO paradigm to address worst-case outcomes.
In principle, the presence of autocorrelations indicates a wave-like nature of observations. 
In a standard RO model, this wave-like behavior is not captured because the model does not explicitly take into account the nature of the autoregressive dependence.
When wave-like components of the uncertain parameter are correlated, their constructive interference can lead to magnified extremes.
As a result, standard RO models are limited in protecting against these intensified extremes. 
However, to the best of our knowledge, connected uncertainty models in the context of RO have not yet been studied in an application-independent fashion.

\subsection{Background on Connected Uncertainty}
\label{sec:backgroundintro}
Uncertainties that are connected across time periods have been studied from multiple application-driven views.
In the context of multistage economic dispatch and unit commitment for electric power, \cite{lorca2015adaptive,lorca2017multistage} study the problem under uncertain renewable power generation.
They incorporate a time series into the uncertainty sets and solve the problem by cut generation, making this approach suitable for large-scale unit commitment problems in a tractable fashion for realistic size datasets.
\citet{jiang2012robust} study a similar unit commitment problem also including hydropower, where the time dependence of the uncertainty is modeled via a single budget constraint.
\citet{ben2005retailer} study a multi-period and two-echelon supply chain system with uncertain demand and model the uncertainty with two types of sets, namely, a box and ellipsoidal sets.
In the box model, the periods are independent of each other, whereas in the ellipsoidal model, different periods are connected to each other via a common budget (size of the set).
The results reveal a more conservative performance for the box set.
\cite{lappas2016multi} study a production system in the context of process scheduling for mixing multiple chemicals.
A single budget constraint is used to bound the uncertainty over all time periods and all sources for a box model.
They demonstrate improved solution quality in comparison with competing approaches.

The temporal connection of uncertainty is leveraged to incorporate changes of cancerous tumor in response to the course of radiation therapy by adjusting the set size as a function of time~\citep{arka2016}.
The temporal periodicity of demand uncertainty in supply chain networks is modeled via periodic affine policies~\citep{eojin2018}.
In inventory management,~\citet{mamani2016closed} study correlated and nonidentically distributed demand over multiple periods.
Under similar conditions,~\citet{ang2012robust} investigate storage assignment. 
For production problems, \citet{xin2015distributionally} consider connected uncertainties by incorporating stochastic processes, such as martingales, to describe the ambiguity sets within a DRO framework.
These works have in common that they model connectedness in the context of a specific application to derive computational or theoretical insights.

\subsection{Positioning}
\label{sec:positioning}
The closest works to ours are the aforementioned studies by~\cite{lorca2015adaptive,lorca2017multistage} in dispatch and unit commitment problem under uncertain renewable power generation.
They leverage a time-series model to capture the high correlations in wind speeds in order to estimate the power production.
This technique reduces the direct impact of uncertainty on the unit commitment problem and demonstrates the benefit of leveraging connected uncertainties.
In general, these and the previously mentioned studies model the connectedness of uncertainties for specific applications, and the resulting uncertainty models are confined to the respective context.
Our work is on general constraints, subsumes many of the available approaches, leverages duality theory for tractable reformulations, and can be used in a broad variety of applications.
Furthermore, we provide extensions to distributional robust methods, relevant to many data-driven applications. 

In addition, previous studies model the connectedness within the set by constraining different uncertainty components, which belong to different times.
Often, they consider budgeted uncertainty sets, which connect uncertainties by imposing a single budget over uncertainties across multiple time periods, making the connections \emph{implicit}.
Besides using a common budget for the uncertainty over time, the authors focus on a specific application. 
Our approach, however, models connectedness over different time periods \emph{explicitly} by an autoregressive time series.
In other words, the connectedness is modeled over individual sets for each period, instead of within components, making the modeling process more intuitive and accessible for decision makers.

These advantages can be measured directly in the context of a specific application by comparing the objective performance and constraint feasibility to those of nonconnected approaches.
The numerical section provides detailed analysis on two applications to demonstrate these advantages.

We study the class of interdependent and adversarial uncertainties in RO problems while maintaining generality to ensure broad applicability.
Specifically, we seek to provide a step towards modeling this family of problems with uncertainty sets that capture connections to previous realizations.
We then extend this approach to the DRO perspective and provide reformulations for moment-based ambiguity sets, which depend on past realizations.
To develop intuition, we begin with an example.

\subsection{Examples}
In this section, we discuss two examples in order to (a) illustrate the difference between connected uncertainty sets and nonconnected sets and (b) examine the advantages that arise when explicitly modeling the connection between uncertainty sets.
\subsubsection*{(a) Example for uncertainty sets.}
\begin{figure}[h]
	\centering
	\begin{tikzpicture}[thick,scale=0.7, every node/.style={transform shape}]
	% Title
	\node [above] at (4.9,4.6) {(a) \textbf{Uncertainty Realizations}};    
	
	%The axes
	\draw[->] (1,-0.2) -- (9,-0.2) node[right] {$T$};
	\draw[->] (1,-0.2) -- (1,4.5) node[above] {$d$};
	
	%First set
	\draw [-] (2,-0.35) node[below] {1} -- (2,-0.05);
	
	%Second Set
	\draw [-] (5,-0.35) node[below] {2} -- (5,-0.05);
	
	%Third Set
	\draw [-] (8,-0.35) node[below] {3} -- (8,-0.05);
	
	% Random realizations
	\draw[-] [draw=gray] (2,1.83) -- (5,2.71) -- (8,3.32);
	\draw[-] [draw=gray] (2,2.81) -- (5,2.94) -- (8,3.18);  
	\draw[-] [draw=gray] (2,1.88) -- (5,2.67) -- (8,3.49);
	\draw[-] [draw=gray] (2,1.70) -- (5,1.06) -- (8,1.74);
	\draw[-] [draw=gray] (2,2.77) -- (5,2.98) -- (8,2.48);
	\draw[-] [draw=gray] (2,1.90) -- (5,2.57) -- (8,2.99);
	\draw[-] [draw=gray] (2,2.00) -- (5,1.80) -- (8,2.02);
	\draw[-] [draw=gray] (2,1.90) -- (5,2.27) -- (8,1.47);
	\draw (2,2.82) circle [radius = 0.05];
	\draw (2,1.83) circle [radius = 0.05];
	\draw (2,1.70) circle [radius = 0.05];
	\draw (2,2.00) circle [radius = 0.05];
	\draw (2,1.90) circle [radius = 0.05];
	\draw (5,2.94) circle [radius = 0.05];
	\draw (5,2.67) circle [radius = 0.05];
	\draw (5,1.06) circle [radius = 0.05];
	\draw (5,2.57) circle [radius = 0.05];
	\draw (5,1.80) circle [radius = 0.05];
	\draw (5,2.27) circle [radius = 0.05];
	\draw (8,3.32) circle [radius = 0.05];
	\draw (8,3.18) circle [radius = 0.05];
	\draw (8,3.49) circle [radius = 0.05];
	\draw (8,1.74) circle [radius = 0.05];
	\draw (8,2.48) circle [radius = 0.05];
	\draw (8,2.99) circle [radius = 0.05];
	\draw (8,2.02) circle [radius = 0.05];
	\draw (8,1.47) circle [radius = 0.05];

	\end{tikzpicture}
	\hspace*{5mm}
	\begin{tikzpicture}[thick,scale=0.7, every node/.style={transform shape}]
	% Title
	\node [above] at (4.9,4.6) {(b) \textbf{ Fixed Uncertainty Sets}};    
	
	%The axes
	\draw[->] (1,-0.2) -- (9,-0.2) node[right] {$T$};
	\draw[->] (1,-0.2) -- (1,4.5) node[above] {$d$};
	
	%First set
	\draw[-] (2,1.1) node[below] {$\mathcal{U}_1$} -- (2,2.9) ;
	\draw[-] (1.8,1.1) -- (2.2,1.1);
	\draw[-] (1.8,2.9) -- (2.2,2.9);
	\draw [fill] (2,2) circle [radius = 0.07];
	\node [left] at (2,2) {\(\mu_1\)};
	\draw [-] (2,-0.35) node[below] {1} -- (2,-0.05);
	
	%Second Set
	\draw[-] (5,1.1) node[below] {$\mathcal{U}_2$} -- (5,2.9);
	\draw[-] (4.8,1.1) -- (5.2,1.1);
	\draw[-] (4.8,2.9) -- (5.2,2.9);  
	\draw [fill] (5,2) circle [radius = 0.07];
	\node [left] at (5,2) {\(\mu_1\)};  
	\draw [-] (5,-0.35) node[below] {2} -- (5,-0.05);
	
	%Third Set
	\draw[-] (8,1.1) node[below] {$\mathcal{U}_3$} -- (8,2.9);
	\draw[-] (7.8,1.1) -- (8.2,1.1);
	\draw[-] (7.8,2.9) -- (8.2,2.9);
	\draw [fill] (8,2) circle [radius = 0.07];
	\node [right] at (8,2) {\(\mu_1\)};  
	\draw [-] (8,-0.35) node[below] {3} -- (8,-0.05);
	
	% Random realizations
	\draw[-] [draw=gray] (2,2.82) -- (5,3.04) -- (8,2.93);
	\draw[dashed] [draw=red] (2,1.83) -- (5,2.71) -- (8,3.32);
	\draw[dashed] [draw=red] (2,2.81) -- (5,2.94) -- (8,3.18);  
	\draw[dashed] [draw=red] (2,1.88) -- (5,2.67) -- (8,3.49);
	\draw[-] [draw=gray] (2,1.70) -- (5,1.06) -- (8,1.74);
	\draw[-] [draw=gray] (2,2.77) -- (5,2.98) -- (8,2.48);
	\draw[-] [draw=gray] (2,1.90) -- (5,2.57) -- (8,2.99);
	\draw[-] [draw=gray] (2,2.00) -- (5,1.80) -- (8,2.02);
	\draw[-] [draw=gray] (2,1.90) -- (5,2.27) -- (8,1.47);
	\draw[-] [draw=gray] (2,1.98) -- (5,2.59) -- (8,1.99);
	
	\end{tikzpicture}
	
	\vspace*{5mm}
	\begin{tikzpicture}[thick,scale=0.7, every node/.style={transform shape}]
	% Title
	\node [above] at (4.9,4.6) {(c) \textbf{Growing Uncertainty Sets}};    
	
	%The axes
	\draw[->] (1,-0.2) -- (9,-0.2) node[right] {$T$};
	\draw[->] (1,-0.2) -- (1,4.5) node[above] {$d$};
	
	%First set
	\draw[-] (2,1.1) node[below] {$\mathcal{U}_1$} -- (2,2.9) ;
	\draw[-] (1.8,1.1) -- (2.2,1.1);
	\draw[-] (1.8,2.9) -- (2.2,2.9);
	\draw [fill] (2,2) circle [radius = 0.07];
	\node [left] at (2,2) {\(\mu_1\)};
	\draw [-] (2,-0.35) node[below] {1} -- (2,-0.05);
	
	%Second Set
	\draw[-] (5,0.9) node[below] {$\mathcal{U}_2$} -- (5,3.1);
	\draw[-] (4.8,0.9) -- (5.2,0.9);
	\draw[-] (4.8,3.1) -- (5.2,3.1);  
	\draw [fill] (5,2) circle [radius = 0.07];
	\node [left] at (5,2) {\(\mu_1\)};  
	\draw [-] (5,-0.35) node[below] {2} -- (5,-0.05);
	
	%Third Set
	\draw[-] (8,0.5) node[below] {$\mathcal{U}_3$} -- (8,3.5);
	\draw[-] (7.8,0.5) -- (8.2,0.5);
	\draw[-] (7.8,3.5) -- (8.2,3.5);
	\draw [fill] (8,2) circle [radius = 0.07];
	\node [right] at (8,2) {\(\mu_1\)};  
	\draw [-] (8,-0.35) node[below] {3} -- (8,-0.05);
	
	% Random realizations
	\draw[-] [draw=gray] (2,2.82) -- (5,3.04) -- (8,2.93);
	\draw[dashed] [draw=red] (2,1.83) -- (5,2.71) -- (8,3.32);
	\draw[dashed] [draw=red] (2,2.81) -- (5,2.94) -- (8,3.18);  
	\draw[dashed] [draw=red] (2,1.88) -- (5,2.67) -- (8,3.49);
	\draw[-] [draw=gray] (2,1.70) -- (5,1.06) -- (8,1.74);
	\draw[-] [draw=gray] (2,2.77) -- (5,2.98) -- (8,2.48);
	\draw[-] [draw=gray] (2,1.90) -- (5,2.57) -- (8,2.99);
	\draw[-] [draw=gray] (2,2.00) -- (5,1.80) -- (8,2.02);
	\draw[-] [draw=gray] (2,1.90) -- (5,2.27) -- (8,1.47);
	\draw[-] [draw=gray] (2,1.98) -- (5,2.59) -- (8,1.99);
	%%%%%%%%%%%%%%%%%%%%%%%%%%%%%%% FIG 1d	
	\end{tikzpicture}
	\hspace*{5mm}
	\begin{tikzpicture}[thick,scale=0.7, every node/.style={transform shape}]
	\node [above] at (4.9,4.6) {(d) \textbf{Connected Uncertainty Sets}};    
	%
	%The axes
	\draw[->] (1,-0.2) -- (9,-0.2) node[right] {$T$};
	\draw[->] (1,-0.2) -- (1,4.5) node[above] {$d$};
	\draw[-] (0.9,3) node[label={[xshift=5mm, yshift=-5mm]{$\hat{d}_1$}}] {} -- (1.1,3);
	\draw[-] (0.9,3.3) node[label={[xshift=5mm, yshift=-3mm]{$\hat{d}_2$}}] {}-- (1.1,3.3);

	%First set
	\draw[-] (2,1.1) -- (2,2.9) ;
	\draw[-] (1.8,1.1) -- (2.2,1.1);
	\draw[-] (1.8,2.9) -- (2.2,2.9);
	\node [below] at (2,1) {$\mathcal{U}_1$};
	\draw [fill] (2,2) circle [radius = 0.07] ;
	\node [left] at (2,2) {\(\mu_1\)};
	\draw [-] (2,-0.35) node[below] {1} -- (2,-0.05);
	
	%Second Set
	\draw[-] (5,1.4)  -- (5,3.2);
	\draw[-] (4.8,1.4) -- (5.2,1.4);
	\draw[-] (4.8,3.2) -- (5.2,3.2);  
	\node [below] at (5,1) {$\mathcal{U}_2(\hat{d}_1)$};
	\draw [fill] (5,2.3) circle [radius = 0.07];
	\node [right] at (5,2.3) {\(\mu_2(\hat{d}_1)\)};
	\draw [-] (5,-0.35) node[below] {2} -- (5,-0.05);
	
	%Third Set
	\draw[-] (8,1.7) -- (8,3.5);
	\draw[-] (7.8,1.7) -- (8.2,1.7);
	\draw[-] (7.8,3.5) -- (8.2,3.5);
	\node [below] at (8,1) {$\mathcal{U}_3(\hat{d}_2)$};
	\draw [fill] (8,2.6) circle [radius = 0.07];
	\node [right] at (8,2.6) {\(\mu_3(\hat{d}_2)\)};
	\draw [-] (8,-0.35) node[below] {3} -- (8,-0.05);
	
	% Random realizations
	\draw[-] [draw=gray] (2,2.82) -- (5,3.04) -- (8,2.93);
	\draw[dashed] [draw=red] (2,1.83) -- (5,2.71) -- (8,3.32);
	\draw[dashed] [draw=red] (2,2.81) -- (5,2.94) -- (8,3.18);  
	\draw[dashed] [draw=red] (2,1.88) -- (5,2.67) -- (8,3.49);
	\draw[-] [draw=gray] (2,1.70) -- (5,1.06) -- (8,1.74);
	\draw[-] [draw=gray] (2,2.77) -- (5,2.98) -- (8,2.48);
	\draw[-] [draw=gray] (2,1.90) -- (5,2.57) -- (8,2.99);
	\draw[-] [draw=gray] (2,2.00) -- (5,1.80) -- (8,2.02);
	\draw[-] [draw=gray] (2,1.90) -- (5,2.27) -- (8,1.47);
	\draw[-] [draw=gray] (2,1.98) -- (5,2.59) -- (8,1.99);
	\end{tikzpicture}
	\caption{\label{fig:set_pic}Uncertainties over Time (a) can be modeled with sets that are Fixed (b),  Growing (c), or Connected (d) with 
		$\mu_t$ being updated using specific realizations $\hat{d}_t$.
	}
\end{figure}
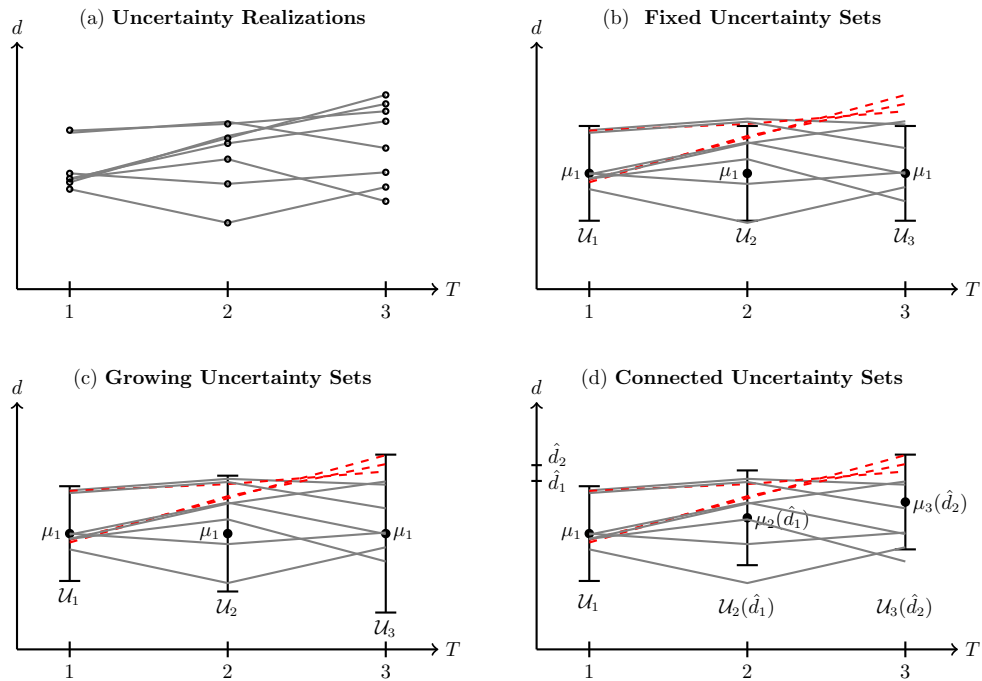

Consider a stylized knapsack problem:
\begin{equation*}
\begin{aligned}
\max_{\mb{x}_t\in \mathcal{X}_t}&\; \sum_{t=1}^{T}\mb{c}_t^{\top}\mb{x}_t\\
\text{s.t.}&\;\sum_{t=1}^{T}\mb{d}_{t}^{\top}\mb{x}_t \leq B \;\;\;\;\forall \mb{d}_t \in \mathcal{U}_t 
&&\forall t=1,\dots,T,
\end{aligned}
\label{pr:basic}
\end{equation*}
where \(\mb{x}_t \in\mathcal{X}_t \subseteq \mathbb{R}^n \,\,\forall t=1,\dots,T\) are decision
variables, \({\mb{c}_t \in \mathbb{R}^{n}}\) are
known, ${\mb{d}_t\in \mathbb{R}^{n}}$ are uncertain coefficients, and \(B\) is
the right-hand side (RHS) coefficient.
In this setting, the uncertain $\mb{d}_t$ can vary in each period while being correlated with the past. 
Figure~\ref{fig:set_pic}(a) illustrates some of the possible uncertainty realizations.
In the paradigm of RO, these uncertainties can be modeled by uncertainty sets $\mathcal{U}_t$ in different ways.
For clarity, consider sets $\mU_t$ that are parameterized by their centers $\bs{\mu}_t$ and sizes $r_t$.
In the following, we discuss three distinct models for \(T = 3\). 
%%%%%%%%%%%%%%%%%%%%%%%%%%%%%%%%%%
%
\begin{itemize}
	\item Often uncertainty sets across periods are modeled to be invariant as  
	\(\bs{\mu}_1 = \bs{\mu}_2 = \bs{\mu}_3\) and \(r_1 = r_2 = r_3\), resulting in \(\mU_1 = \mU_2 = \mU_3\) as illustrated in~Figure~\ref{fig:set_pic}(b).
	Although this fixed model offers simplicity and computational advantages, it fails to capture all possible scenarios, when uncertainties are actually auto-correlated. 
	\item To overcome this limitation, \(\mU_t\) can be modeled to grow in size, i.e., \(\bs{\mu}_1 = \bs{\mu}_2 = \bs{\mu}_3\) and \(r_1 \leq r_2 \leq r_3\), resulting in \(\mU_1 \subseteq \mU_2 \subseteq \mU_3\).
	Such growing uncertainty sets are sketched in Figure~\ref{fig:set_pic}(c).
	The expansion, however, may render the solutions over conservative. 
	\item In both previous models, the uncertainty set in each period is independent of the realizations in previous periods.
	We propose a new uncertainty set model $\mathcal{U}_{t} = \mathcal{U}_{t}(\mb{d}_{t-1})$ for \(t > 1\), where the center of the set  is correlated with the previous period uncertainty realization. 
	Such a setting occurs, for example, when
	\(\bs{\mu}_2(\mb{d}_1) = \mb{A}\mb{d}_1\) and \(\bs{\mu}_3(\mb{d}_2) = \mb{A}\mb{d}_2\) with \(\mb{A} = \alpha \mb{I}\), \(\alpha > 0\), and \(r_1 = r_2 = r_3\). 
	Such sets are illustrated in Figure~\ref{fig:set_pic}(d) for positively correlated centers.
	This model has the advantage of capturing autocorrelated uncertainties without increasing set size. 
\end{itemize}

\subsubsection*{(b) Example for feasibility of the solution.}
Consider the constraint
\begin{equation}
	\label{eq:example}
	\begin{aligned}
		& d_1 x_1 + d_2 x_2 \leq B \;\;\forall d_2 \in \mU_2(d_1),\; \forall d_1 \in \mU_1\\
		& x_1, x_2 \geq 0,
	\end{aligned}
\end{equation}
where the uncertainty sets are given by
\begin{align*}
	\mU_1 &= [-1, 1] \; \text{and}\\
	\mU_2(d_1) &= [\mu_2(d_1) - 1, \mu_2(d_1) + 1], \text{ where } \mu_2(d_1) = (\alpha - 0.5)  d_1.
\end{align*}
Notice that the parameter $\alpha$ controls the connectedness between the periods 1 and 2.
Using these sets, we can rewrite the robust counterpart of the Constraint~\eqref{eq:example} to obtain
\begin{align*}
	|x_1 + (\alpha - 0.5) x_2| + x_2  \leq B.
\end{align*}
Depending on the value of $\alpha$, the nature of the above constraint varies significantly. 

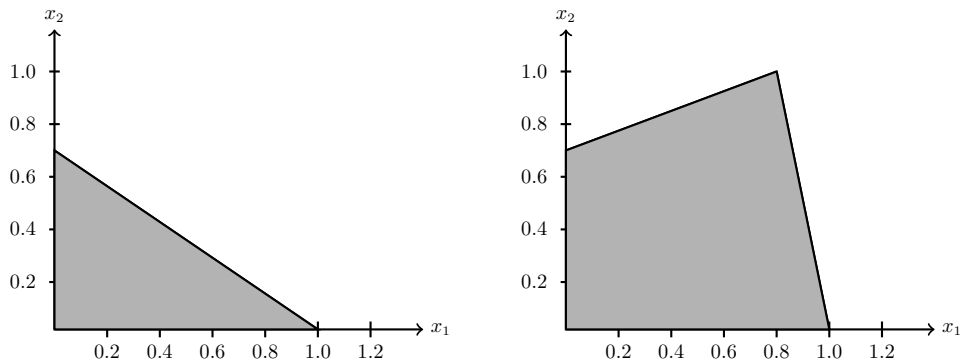
\begin{figure}[h]
	\centering
	\begin{tikzpicture}[thick,scale=0.7, every node/.style={transform shape}]
		% Title
		
		\coordinate (A) at (1,3.2) ;
		\coordinate (B) at (6,-0.2) ;
		\coordinate (O) at (1,-0.2) ;

		%The axes
		\draw[->] (1,-0.2) -- (8,-0.2) node[right] {$x_1$};
		\draw[->] (1,-0.2) -- (1,5.5) node[above] {$x_2$};

		%yticks
		\draw[-] (0.9,4.7) node[label={[xshift=-5mm, yshift=-4mm]{$1.0$}}] {}-- (1.1,4.7);
		\draw[-] (0.9,3.7) node[label={[xshift=-5mm, yshift=-4mm]{$0.8$}}] {} -- (1.1,3.7);
		\draw[-] (0.9,2.7) node[label={[xshift=-5mm, yshift=-4mm]{$0.6$}}] {} -- (1.1,2.7);
		\draw[-] (0.9,1.7) node[label={[xshift=-5mm, yshift=-4mm]{$0.4$}}] {} -- (1.1,1.7);
		\draw[-] (0.9,0.7) node[label={[xshift=-5mm, yshift=-4mm]{$0.2$}}] {} -- (1.1,0.7);
		
		%xticks
		\draw [-] (2,-0.35) node[below] {0.2} -- (2,-0.05);
		\draw [-] (3,-0.35) node[below] {0.4} -- (3,-0.05);
		\draw [-] (4,-0.35) node[below] {0.6} -- (4,-0.05);
		\draw [-] (5,-0.35) node[below] {0.8} -- (5,-0.05);
		\draw [-] (6,-0.35) node[below] {1.0} -- (6,-0.05);
		\draw [-] (7,-0.35) node[below] {1.2} -- (7,-0.05);
		
		%First set
		\draw[-] (1,3.2) -- (6,-0.2);

		\draw[black, fill=gray!60] (O) -- (A) -- (B) -- (O) -- cycle;

		%%%%%%%%%%%%%%%%%%%%%%%%%%%%%%% FIG 1d	
	\end{tikzpicture}
	\hspace*{5mm}
	\begin{tikzpicture}[thick,scale=0.7, every node/.style={transform shape}]
		
		\coordinate (A) at (1,3.2) ;
		\coordinate (B) at (6,-0.2) ;
		\coordinate (C) at (5, 4.7);
		\coordinate (O) at (1,-0.2) ;
		
		%The axes
		\draw[->] (1,-0.2) -- (8,-0.2) node[right] {$x_1$};
		\draw[->] (1,-0.2) -- (1,5.5) node[above] {$x_2$};
		
		%yticks
		\draw[-] (0.9,4.7) node[label={[xshift=-5mm, yshift=-4mm]{$1.0$}}] {}-- (1.1,4.7);
		\draw[-] (0.9,3.7) node[label={[xshift=-5mm, yshift=-4mm]{$0.8$}}] {} -- (1.1,3.7);
		\draw[-] (0.9,2.7) node[label={[xshift=-5mm, yshift=-4mm]{$0.6$}}] {} -- (1.1,2.7);
		\draw[-] (0.9,1.7) node[label={[xshift=-5mm, yshift=-4mm]{$0.4$}}] {} -- (1.1,1.7);
		\draw[-] (0.9,0.7) node[label={[xshift=-5mm, yshift=-4mm]{$0.2$}}] {} -- (1.1,0.7);

		%xticks
		\draw [-] (2,-0.35) node[below] {0.2} -- (2,-0.05);
		\draw [-] (3,-0.35) node[below] {0.4} -- (3,-0.05);
		\draw [-] (4,-0.35) node[below] {0.6} -- (4,-0.05);
		\draw [-] (5,-0.35) node[below] {0.8} -- (5,-0.05);
		\draw [-] (6,-0.35) node[below] {1.0} -- (6,-0.05);
		\draw [-] (7,-0.35) node[below] {1.2} -- (7,-0.05);
		
		%First set
		\draw[-] (1,3.2) -- (5,4.7) ;
		\draw[-] (5,4.7) -- (6,-0.2);
		
		\draw[black, fill=gray!60] (O) -- (A) -- (C) -- (B) -- (O) -- cycle;
		
	\end{tikzpicture}
	\caption{\label{fig:example}The shaded area represents the feasible region. The (left) figure represents $\alpha = 1.4$ and the (right) figure represents $\alpha = 0.4$.}
\end{figure}

The reformulated constraint includes the parameter $\alpha$ explicitly, connecting uncertainties at different time periods directly with each other in a parameterized fashion. 
Figure~\ref{fig:example} illustrates that the shape of the feasibility region can change significantly as a direct consequence of variations in the connectedness parameter $\alpha$. 
The advantage of this approach is further illustrated in the following section where we compare the method to the well studied budget constrained uncertainty sets.

\subsection{Implicit vs Explicit Connection}
As we have stated previously, our approach allows us to explicitly model the connection between uncertainties whereas existing approaches only allow for an implicit connection between the uncertainties across different periods. 
Consider the following uncertainty set which is very common in many robust optimization problems. 
$$\mU = \left\{(d_1, d_2) \mid |d_1| + |d_2| \leq B\right\}.$$
This is a budget uncertainty set which has gained significant popularity due to its tractability as well as the reduced conservatism compared to standard box uncertainty sets. 
In this set, the value of the uncertain parameters $d_2$ is limited on the basis of the realization of $d_1$, i.e., 
$$\mU_2(d_1) = \left\{d_2 \mid |d_2| \leq B - |d_1|\right\}.$$
As such, there is an implicit connection present between the uncertain parameters in different periods. 
The only way to control the dependence of $d_2$ on $d_1$ is through the parameter $B$ which simultaneously also controls the conservativeness of the set. 
This limits its flexibility. 
Instead of the above set, we can use a connected uncertainty set of the form 
$$\mU = \left\{(d_1, d_2) \mid |d_1| \leq B, |d_2 - \mu_2(d_1)| \leq B  \right\},$$
where $\mu_2(d_1)$ is a linear function of $d_1$ of the form $A d_1 + b$ which controls the dependence of $d_2$ on $d_1$. 
This creates an explicit connection between $d_1$ and $d_2$ and also allows us the flexibility of changing the connection without affecting the conservatism of the problem.

We will revisit this example in Section~\ref{sec:appl-3:-knaps} with more detailed numerical analysis.
For clarity, we now formally define the proposed uncertainty sets.
\begin{definition}[\cu~set]
	If any of the parameters describing an uncertainty set explicitly depend on the realizations of uncertainty from previous periods, then this set is a \emph{connected set}. 	
\end{definition}

Consequently, we call realizations from~\cu~sets as \emph{connected uncertainties}.
For example, the set illustrated in~\autoref{fig:set_pic}(d) is a connected uncertainty set, whereas those in (b) and (c) are not.

\subsection{Contributions}
\label{sec:contrib}
We study optimization problems over multiple time periods with both static and adaptive decisions, where the uncertainty at each period is influenced by the uncertainty realizations in the past.
Specifically, the contributions of this paper are:
\begin{itemize}
	\item We introduce the concept of connected uncertainty sets that offer generalizable insights. Specifically, we discuss uncertainties with linear and quadratic dependence on previous periods.
	\item We provide tractable reformulations with here-and-know decisions for common families of uncertainty models, namely, i) polyhedral \cus, where the RHS parameters of the set depend on previous periods, ii) ellipsoidal \cus, where the center depends on the past; or (iii) ellipsoidal sets, where the covariance matrix depends on previous periods.
	\item We also provide tractable reformulations for affinely adaptive decisions for polyhedral and ellipsoidal uncertainty sets, which are commonly used sets in practical settings.
	\item We probe the proposed concept through numerical experiments. 
\end{itemize}

Note the following two remarks.
First, such connected sets $\mathcal{U}_{t}(\mb{d}_{t-1})$ can also be expressed as a single joint uncertainty set over \(\mb{d}_1,\dots,\mb{d}_T\) to ease the reformulation. 
Second, such a joint set is different and smaller than the box set that can be created by projecting the uncertainty onto each individual $\mb{d}_t$ and taking the product.
A key advantage of the proposed \cus~is their explicit dependence on previous realizations.
This modeling power can be exploited, e.g., when the uncertainty stems from a time series.
Furthermore, if an application requires the periodic resolving of the problem, our approach directly prescribes how to update the sets over time.

\subsubsection*{Notation.}
All decisions are \emph{here-and-now} (nonadaptive) unless specified otherwise.
Bold lowercase and uppercase letters denote vectors and matrices.
For any matrix \(\mb{A}\), the \(i^{th}\) row of the matrix is denoted by  \(\mb{A}_i^{\top}\).
To streamline the exposition, we use ``uncertainty set'' for both the RO and DRO settings.
The former is over parameters and the latter over distributions and is also known as an ambiguity set.
\(T\) denotes the total time periods, \(\tau\) refers to a particular time period, \(t\) serves as an index, and $\mb{e}$ is a vector of all ones. 

%
%%%%%%%%%%%%%%%%%%%%%%%%%%%%%%%%%%%%%%%%%%%%%%%%%%%%%%%%%%%%%%%%%%%%%%%%%%%%%%%%%%%%%%%
%
\section{Connected Uncertainty with RO}
\label{sec:background}

The concept of interstage dependence in multistage problems is well studied in stochastic optimization.
By naturally embedding scenario trees, 
it can be used to model complex state-dependent processes, such as financial volatility~\citep{hoyland2001generating}.
Some models assume independence, allowing decomposition algorithms and sharing of cuts within the same stage--that is, the L-shaped method~\citep{infanger1996cut}.
In general, sharing cuts is not permitted when the scenarios are dependent, unless
the dependency follows simple time series models, in which case these cuts can be modified and shared~\citep{de2013sharing}.

In RO,  uncertainty models typically employ a single set for all periods, as discussed in Section~\ref{sec:positioning}.
The popular budgeted uncertainty set models use one budget across periods to couple the uncertainty realizations.
In this paper, we consider the dependence on past realizations for general-purpose uncertainty sets rather than application-specific sets.
We reformulate a standard linear constraint over a sequence of uncertainty sets, whose parameters (centers, covariances, etc.) are connected across time.

%%%%%%%%%%%%%%%%%%%%%%%%%%%%%%%%%%%%%%%%%%%%%%%%%%%%%%%%%%%%%%%%%%%%%%%%%%%%%%%%%%%%%%%%%%%%%%%%%%%%%%%%%%%%%

\subsection{Model}
\label{sec:ro}

In a variety of applications, such as cutting stock and packing problems, capital budgeting, and project selection etc., linear constraints are typically subject to connected uncertainties~\citep{freville2004multidimensional}. 
Because this effect can be confined to constraint coefficients, the goal of this section is to reformulate the constraint 
\begin{equation}
\label{eq:prcs}
\tag{C-RO}
\sum_{t=1}^{T}\mb{d}_t^{\top}\mb{x}_t \leq B \quad \forall \mb{d}_t \in \mathcal{U}_t(\mb{d}_{t-1}) \;\;\forall \mb{d}_1 \in \mathcal{U}_1 \;\;\forall t=2,\dots,T.
\end{equation}
Here, $\mb{d}_t$ is the vector of uncertain coefficients, $\mb{x}_t$ denotes the decision variables, and $B$ is a
constant upper bound.
In each period, $\mb{d}_t$ resides in a set
$\mathcal{U}_t(\cdot)$, which may depend on $\mb{d}_{t-1}$.
The robust counterpart of~\eqref{eq:prcs} becomes
\vspace{-2mm}
\begin{equation}
\label{eq:rc}
\tag{RC}
\max_{\mb{d}_1 \in \mU_1}\{\mb{d}_1^{\top}\mb{x}_1 + \max_{\mb{d}_2 \in
	\mU_2(\mb{d}_1)} \{\mb{d}_2^{\top}\mb{x}_2 + \dots +  \max_{\mb{d}_T \in
	\mU_T(\mb{d}_{T-1})}\mb{d}_T^{\top}\mb{x}_T\}\} \leq B.
\vspace{-2mm}
\end{equation}
Each period $t$ is affected by the worst-case realization of the uncertain parameter in period $t-1$.
In this work, we focus on \cus, whose parameters linearly depend on the previous realization.
Such uncertainties are prevalent in modeling popular applications--for example., autoregressive models and minimum mean square error predictors (linear or jointly normal)~\citep{hoggallen}.
In the RO literature, when confidence regions of the uncertainty are available, ellipsoidal sets are commonly used, whereas when only bounds on the uncertainty are known, polyhedral uncertainty sets are used. 
Therefore, we discuss \cus~of ellipsoidal and polyhedral structures to extend the benefits of RO to multi-period settings.
For ellipsoidal sets, we study the dependence of the set center and of the covariance matrix on the past, which may arise from either a Bayesian or frequentist update, or from a time-series model. 
For polyhedral sets, we focus on the RHS set coefficients that may depend on previous realizations.
This occurs when the magnitude or the location of the uncertainty is affected by past realizations.
%

%%%%%%%%%%%%%%%%%%%%%%%%%%%%%%%%%%%%%%%%%%%%%%%%%%%%%%%%%%%%%%%%%%%%%%%%%%%%%%%%%%%%%%%%%%%%%%%%%%%%%%%%%%%%%

\subsection{Center Dependence of Ellipsoidal Sets}
\label{sec:ellips_center}
When uncertainties are Gaussian, they can be modeled with ellipsoidal sets, as Gaussian distributions have ellipsoidal contours. Confidence or predictive regions can also be naturally described by ellipsoids~\citep{charnes1993multivariate}. 
For such sets, there are three parameters that can be affected by previous uncertainty
realizations: radius $r$, center $\bs{\mu}_t$, or covariance matrix $\bs{\Sigma}_t$.
The case of radius dependence leads to a nested norm structure, resulting in nonconvex problems that are beyond our scope.
Here, we discuss the setting where $\bs{\mu}_t$ depends on the previous period realizations, whereas $r$ and $\bs{\Sigma}_t$ are constant.
Such dependence arises, when, for example, the demand is autoregressive~\citep{alwan2016dynamic}.
Section~\ref{sec:matrix-dependence} discusses $\bs{\Sigma}_t$ dependence, whereas $r$ and $\bs{\mu}_t$ are invariant.

When the set centers $\bs{\mu}_t$ are connected, the uncertainty set for each period is
\vspace{-2mm}
\begin{equation}
\label{elips_u}\tag{E}
\mathcal{U}_t(\mb{d}_{t-1}) = \{\mb{d}_t \mid \mb{d}_t = \bs{\mu}_t(\mb{d}_{t-1}) + \mb{L}_t \mb{u}_t : \|\mb{u}_t\|_2 \leq r_t\},
\vspace{-2mm}
\end{equation}
where \(\mb{L}_t \mb{L}_t^{\top} = \bs{\Sigma}_t\). 
If $\bs{\mu}_t$ is autoregressive and depends on the previous period realization as
\vspace{-2mm}
\begin{equation}
\label{eq:center_update}
\bs{\mu}_{t+1}(\mb{d}_{t}) = \mb{A}_{t} \bs{\mu}_{t}(\mb{d}_{t-1}) + \mb{F}_{t} \mb{d}_{t} + \mb{c}_{t},
\vspace{-2mm}
\end{equation}
the constraint~\eqref{eq:prcs} can be reformulated as with the upcoming Theorem~\ref{prop:ellips_center}.

The nature of the connectedness motivates \(\mb{A}_t\), \(\mb{F}_t\) and \(\mb{c}_t\)--for example, by autoregressive or Bayesian models, etc.
We reformulate the uncertainty in a period assuming all previous realizations are known.
The robust counterpart at $t$ has to take future uncertainties into account. 
As such, all decisions are affected by succeeding decisions.
We define recursive variables for each period \({k \in \{1\dots,T-1\}}\) as:
\begin{eqnarray*}
	\mb{y}_{k} &=& \mb{x}_{k} + (\mb{F}_{k} + \mb{A}_{k})^{\top}\mb{y}_{k+1}, \text{ with } \mb{y}_{T} = \mb{x}_T,\\
	C_{k} &=& \mb{c}_{k}^{\top}\mb{y}_{k+1} + C_{k+1}, \text{ with } C_{T} = 0 \text{ and}\\
	R_{k} &=& r_{k}\|L_{k}^{\top}(\mb{x}_{k} + \mb{F}_{k}\mb{y}_{k+1})\|_2 + R_{k+1} \text{ with } R_T = r_T\|\mb{L}_T^{\top}\mb{x}_T\|_2.
\end{eqnarray*}
The aggregates $C_{k}$ and $R_{k}$ are functions of the variable $\mb{y}_{k}$ and represent the effect of the constants $\mb{c}_k$ in  the update~\eqref{eq:center_update} and the 2-norm protection term, respectively.
\begin{theorem}
	\label{prop:ellips_center}
	The robust counterpart of Constraint~\eqref{eq:prcs} for the ellipsoidal set~\eqref{elips_u} is
	\[\bs{\mu}_{1}^{\top}\mb{y}_{1} + C_{1} + R_{1} \leq B.\]
	%.
\end{theorem}
\proof{Proof:}
	We start by reformulating~\eqref{eq:prcs} for an arbitrary \(k\). 
	Let \(s_{k} = \sum_{t=1}^{k}\mb{d}_t^{\top}\mb{x}_t\);  
	then, the robust counterpart of~\eqref{eq:prcs} for \(k\) is 
	\begin{equation}
	\label{eq:prop_ec_lem}
	s_{k} + \bs{\mu}_{k+1}(\mb{d}_{k})^{\top}\mb{y}_{k+1} + C_{k+1} + R_{k+1} \leq B.
	\end{equation}
	We prove the statement for \(\mb{d}_T\), then assume it to be true for \(\mb{d}_{\tau+1}\), before proving it for \(\mb{d}_{\tau}\).
	\paragraph{Base case (\(k =T-1\))}
	The Constraint~\eqref{eq:prcs} can be expanded as 
	\(s_{T-1} + \mb{d}_{T}^{\top}\mb{x}_T + C_{T} \leq B.\)
	Because the constraint must hold for all \(\mb{u}_T \in \{\mb{u}_T \mid \|\mb{u}_T\|_2 \leq r_T\}\), it also holds for the robust counterpart
	\[s_{T-1} + \bs{\mu}_{T}(\mb{d}_{T-1})^{\top}\mb{x}_T + r_T\|\mb{L}_{T}^{\top}\mb{x}_{T}\|_2 +  C_{T} \leq B,\]
	as \(\mb{d}_T = \bs{\mu}_T(\mb{d}_{T-1}) + \mb{L}_T \mb{u}_T \).
	With \(\mb{y}_T = \mb{x}_T\), \(C_{T} = 0\), and \({R_T =  r_T\|\mb{L}_{T}^{\top}\mb{x}_T\|_2}\), we obtain the result
	\[s_{T-1} + \bs{\mu}_{T}(\mb{d}_{T-1})^{\top}\mb{y}_{T} + C_{T} + R_{T} \leq B.\]
	\paragraph{Inductive case (\(k=\tau\))}
	Assume that the reformulation~\eqref{eq:prop_ec_lem} holds for \(\tau+1\).
	Then the robust counterpart of~\eqref{eq:prcs} with respect to \(\mb{d}_T,\mb{d}_{T-1}, \dots, \mb{d}_{\tau+1}\) is
	\[s_{\tau} + \bs{\mu}_{\tau+1}(\mb{d}_{\tau})^{\top}\mb{y}_{\tau+1} + C_{\tau+1} + R_{\tau+1} \leq B.\]
	Substituting the mean and rearranging the terms, this can be expressed as
	\begin{align*}
	s_{\tau-1} &+ \mb{d}_{\tau}^{\top}(\mb{x}_{\tau} + \mb{F}_{\tau}^{\top}\mb{y}_{\tau+1}) +  \bs{\mu}_{\tau}(\mb{d}_{\tau-1})^{\top}\mb{A}_{\tau}^{\top}\mb{y}_{\tau+1}  \\
	&+\mb{c}_{\tau}^{\top}\mb{y}_{\tau+1} + C_{\tau+1} + R_{\tau+1} \leq B.
	\end{align*}
	\noindent Using the uncertainty set~\eqref{elips_u}, this can be rewritten as
	\begin{align*}
	s_{\tau-1} &+ \bs{\mu}_{\tau}(\mb{d}_{\tau-1})^{\top}(\mb{x}_{\tau} + \mb{F}_{\tau}^{\top}\mb{y}_{\tau+1}) \\
	& +\bs{\mu}_{\tau}(\mb{d}_{\tau-1})^{\top}\mb{A}_{\tau}^{\top}\mb{y}_{\tau+1} + \mb{u}_{\tau}^{\top}\mb{L}_{\tau}^{\top}(\mb{x}_{\tau} + \mb{F}_{\tau}^{\top}\mb{y}_{\tau+1}) 	  \\
	&+ \mb{c}_{\tau}^{\top}\mb{y}_{\tau+1} + C_{\tau+1} + R_{\tau+1} \leq B.
	\end{align*}
	Taking the robust counterpart with respect to \(\mb{u}_{\tau}\), we can write
	\begin{align*}
	s_{\tau-1} &+ \bs{\mu}_{\tau}(\mb{d}_{\tau-1})^{\top}(\mb{x}_{\tau} + \mb{F}_{\tau}^{\top}\mb{y}_{\tau+1} + \mb{A}_{\tau}^{\top}\mb{y}_{\tau+1}) \label{eq:center_worst}\\
	&+ r_{\tau}\|\mb{L}_{\tau}^{\top}(\mb{x}_{\tau} + \mb{F}_{\tau}^{\top}\mb{y}_{\tau+1})\|_2  \\
	&+ \mb{c}_{\tau}^{\top}\mb{y}_{\tau+1} + C_{\tau+1} + R_{\tau+1} \leq B.\nonumber
	\end{align*}
	With \(\mb{y}_{\tau} = \mb{x}_{\tau} + \mb{F}_{\tau}^{\top}\mb{y}_{\tau+1} + \mb{A}_{\tau}^{\top}\mb{y}_{\tau+1}\) 	and the definitions of  \(R_{\tau}\) and \(C_{\tau}\), we obtain the desired result
	\begin{align*}
	s_{\tau-1} + \bs{\mu}_{\tau}(\mb{d}_{\tau-1})^{\top}\mb{y}_{\tau}  + C_{\tau} + R_{\tau} \leq B.
	\end{align*}
	This concludes the induction, and the final reformulation is obtained by substituting \(\tau = 1\). 
\Halmos
\endproof

In summary, when modeling uncertainty with connected centers, a recursive variable  \(\mb{y}_{k+1}\) is required to protect against the accumulating effect of future uncertainties.
Current RO models, which model uncertainties at each period to be independent, 
neglect this effect. 
Thus, they do not capture all uncertainties or lead to over conservative solutions. 
In order to optimally account for uncertainty over all periods, Theorem~\ref{prop:ellips_center} prescribes how to use this \(\mb{y}_{k+1}\) to modify the 2-norm in order to protect against uncertainty at each period.
We will numerically demonstrate this result for connected centers on a knapsack application in Section~\ref{sec:appl-3:-knaps} and show that they can improve constraint satisfaction and increase objective function value for any given level of constraint satisfaction, when compared with nonconnected sets.

\textbf{Remark}: Alternatively, a joint uncertainty set can be constructed by a product of the individual uncertainty sets~\eqref{elips_u} over all time periods.
This joint set is second-order cone representable as such, the Constraint~\eqref{eq:prcs} can be tractably reformulated over this set as well. 

%%%%%%%%%%%%%%%%%%%%%%%%%%%%%%%%%%%%%%%%%%%%%%%%%%%%%%

\subsection{Matrix Dependence of Ellipsoidal Sets}
\label{sec:matrix-dependence}
The previous section focused on changes in the location of ellipsoidal uncertainty sets.
Many problems require models that allow for the shape of the set to vary--for example, when the volatility in an uncertain process depends on past realizations.
Such uncertainties are usually described by autoregressive conditional heteroskedastic models--for example, in asset prices~\citep{bollerslev1992arch} or demand and sales growth of firms~\citep{larson2015empirical}.
To this end, we model the covariance \(\bs{\Sigma}_t\) of the ellipsoidal uncertainty set to depend on previous realizations, i.e., $\bs{\Sigma}_{t}(\mb{d}_{t-1})$, while
the radius $r$ and mean \(\bs{\mu}_t\) are predetermined. 
A general form of this dependence is
\begin{equation}\label{eq:mu}
\bs{\Sigma}_{t+1}(\mb{d}_{t}) = a_{t} \cdot \bs{\Sigma}_{t}(\mb{d}_{t-1}) + f_{t} \cdot (\mb{d}_{t} - \bs{\mu}_{t})(\mb{d}_{t} - \bs{\mu}_{t})^{\top} + \mb{C}_{t} \;\;\forall t,
\end{equation}
where \(a_t \geq 0\), \(f_t \geq 0\), \(\mb{C}_t \succeq 0\), and \(\bs{\Sigma}_1 \succeq 0\) are constants. 
This form of dependence helps to capture marginal updates for both the frequentist and Bayesian paradigms.
The former assumes no prior distribution and the latter assumes one.
We omit the dependence of the mean $\bs{\mu}$ on the past in this setting to simplify the resulting reformulation.
Because the second term in~\eqref{eq:mu} is rank 1 and $\mb{C}_t$ is positive semi-definite, the update~\eqref{eq:mu}  preserves the positive semi-definiteness of $\bs{\Sigma}_t$.

The quadratic dependence in~\eqref{eq:mu} leads to nonlinear terms in the robust counterpart, obscuring an analytic reformulation.
For this, Theorem~\ref{prop:matrix_depen} provides a conservative reformulation.
Since the quadratic term can arise from a positive or a negative deviation from $\bs{\mu}_{t}$, we introduce a sign variable \(n\) with  \(\mb{n}_k = (n_{k,k},n_{k,k+1},\dots,n_{k,T})\) to differentiate between these cases.
Each \(n_{k,t}\) can  be \(1\) or \(-1\) except for $n_{T,T} = 1$, and
\(\mb{n}_{[k]}^\top = (\mb{n}_{k}^\top,\dots,\mb{n}_{T}^\top)\). 
The set \(\mathcal{N} = \mathcal{N}_{[1]}\) consists of all possible \(2^{\frac{T}{2}(T-1)}\) vectors \(\mb{n}\), where \(\mathcal{N}_{[k]}\) consists of all combinations of \(\mb{n}_{[k]}\). 
For \(k \in \{1,\dots,T-1\}\), we also define  
\begin{align*}
\mb{y}_{k}(\mb{n}_{[k]}) &= \mb{x}_{k} + \sum_{t=k+1}^{T} n_{k,t} \cdot \mb{y}_{t} (\mb{n}_{[t]})\;r_{t}\sqrt{A_{k+1,t}f_{k}},\\
&\text{ with } \mb{y}_{T} = \mb{x}_T,\\
R_{k} &= R_{k+1} + \sum\limits_{t = k+1}^{T}r_{t}\sqrt{A_{k+1,t}\;\mb{y}_{t}(\mb{n}_{t})^{\top}\mb{C}_{k} \mb{y}_{t}(\mb{n}_{t})}, \\
&\text{ with } R_{T+1} =0, \text{and}\\
A_{k,t} &= \prod_{j=k}^{t}a_j ,\;k=1,\dots,t, \\
&\text{ with }  A_{t,t} = 1.
\end{align*}
In what follows, we suppress the dependence of \(\mb{y}_{k}\) on \(\mb{n}_{[k]}\) 
and \(\bs{\Sigma}_{t+1}\) on \(\mb{d}_t\) for brevity.
To reformulate~\eqref{eq:prcs}, let \(\mb{d}_{k}\) for some \(k \in \{1,\dots,T\}\) reside in a set~\eqref{elips_u}  with radius \(r_{k}\), center \(\bs{\mu}_{k}\), and covariance matrix \(\bs{\Sigma}_{k}\) with
\(\bs{\Sigma}_{k} = \mb{L}_{k}\mb{L}_{k}^{\top}\) updated as in~\eqref{eq:mu}.
\begin{theorem}
	\label{prop:matrix_depen}	
	A conservative robust reformulation of~\eqref{eq:prcs} 
	is
	\begin{eqnarray*}
		&&\sum_{t=1}^{T}\bs{\mu}_t^{\top}\mb{x}_t + \sum_{t=1}^{T}r_{t}\sqrt{A_{1,t	}}\;\|\mb{L}_1^{\top} \mb{y}_{t}(\mb{n}_{[t]}) \|_2 + R_1 \leq B \;\;\forall \mb{n}_{[1]} \in \mathcal{N}_{[1]}.
	\end{eqnarray*}
\end{theorem}
The proof of Theorem~\ref{prop:matrix_depen} is provided in the electronic companion. 
It proceeds by induction, similar to the proof of Theorem~\ref{prop:ellips_center}.
We start by proving the base case and then extend it to arbitrary cases.
Here, we use the uncertainty dependence described in Equation~\ref{eq:mu}. 
We capture the quadratic dependence by considering both the positive and negative cases in the resulting constraint. 
The reformulation in Theorem~\ref{prop:matrix_depen} contains a sum of 2-norms, making it computationally tractable and attractive for large-scale problems with few periods.
Similar to center dependence in Section~\ref{sec:ellips_center}, \(\mb{y}_{k}\) plays the key role in adapting the protection against connected matrices.
This advances the conventional RO paradigm because the presence of connectedness allows current worst-case realizations to affect future uncertainties, and Theorem~\ref{prop:matrix_depen} describes how \(\mb{y}_{k}\) modifies the 2-norm through the summation $\sum_{t=k+1}^{T} n_{k,t} \cdot \mb{y}_{t} (\mb{n}_{[t + 1]})\;r_{t}\sqrt{A_{k+1,t}f_{k}}$ in order to incorporate this affect.
Note that the reformulation in Theorem~\ref{prop:matrix_depen} has exponentially many constraints in $T^2$, limiting its applicability to problems with a few periods, as is often the case in many multiperiod problems.

\subsubsection*{Summary for Ellipsoidal Sets}
For center and matrix dependence, 
we showed how $\bs{\mu}_t(\mb{d}_{t-1})$ and $\bs{\Sigma}_t(\mb{d}_{t-1})$ can be affected by past uncertainties, and that the worst case is driven by future uncertainties through $\mb{y}_{k}$.
These results explicitly reveal the effect of connected uncertainties, namely that the past dependence establishes the \emph{location} (for $\bs{\mu}_t$ dependence via~\eqref{eq:center_update} or for $\bs{\Sigma}_t$ dependence  via~\eqref{eq:mu}), while the future connections determine the \emph{direction} of the worst case (via theorems~\ref{prop:ellips_center} for $\bs{\mu}_t$ or~\ref{prop:matrix_depen} for $\bs{\Sigma}_t$).
A setting with both $\bs{\mu}_t$ and $\bs{\Sigma}_t$ changing concurrently results in a highly nonconvex connection between the uncertainties, which is beyond our scope.

\subsection{Polyhedral Sets}
\label{sec:polyhedral}
Here, we extend these models to the \cu~setting, aiming to provide further intuition.
Specifically, we consider problems, where the uncertain constraint coefficients of each period reside in a polyhedral CU set and the RHS parameters depend on previous period uncertainty.
To reformulate the Constraint~\eqref{eq:prcs}, let the uncertain \(\mb{d}_t\)  reside in a polyhedral set 
\begin{equation}
\mU_t(\mb{d}_{t-1}) = \{\mb{d}_t \mid \mb{G}_t \mb{d}_t \geq \mb{g}_t + \bs{\Delta}_t \mb{d}_{t-1}\},
\label{eq:polyUR}
\tag{P}
\end{equation}
where the matrix \(\bs{\Delta}_t\)
is application based--for example, from time-series models.
Here, without loss of generality,  the parameters of \(\mU_1\) are considered as known
(\({\bs{\Delta}_1 = \mb{0}}\)).
This setting resembles the introductory Example~3, and the following theorem reformulates~\eqref{eq:prcs}.
\begin{theorem}
	\label{prop_poly_rhs}
	The robust counterpart of constraint~\eqref{eq:prcs} under the uncertainty set~\eqref{eq:polyUR} is
	\begin{align*}
	\sum_{t=1}^{T}\mb{q}_t^{\top}\mb{g}_t &\leq B &&\\
	\mb{q}_t^{\top}\mb{G}_t &= \mb{x}_t^{\top} + \mb{q}_{t+1}^{\top}\bs{\Delta}_{t+1} &&\forall t = 1,\dots,T,	\\
	\mb{q}_t &\leq 0 &&\forall t=1,\dots,T,
	\end{align*}
	where \(\bs{\Delta}_1 = 0\) and \(\;\bs{\Delta}_{T+1} = 0\).
\end{theorem}
The proof leverages duality and is relegated to the electronic companion.
In this theorem, the term \(\mb{q}_t^{\top}\mb{g}_t\) safeguards against all uncertainty realizations. 
The dual variable \(\mb{q}_t\) is appropriately adjusted by \(\mb{q}_t^{\top}\mb{G}_t = (\mb{x}_t + \bs{\Delta}_{t+1}^{\top}\mb{q}_{t+1})^{\top}\) to account for uncertainties.
The connectedness--that is, the nonzero $\bs{\Delta}_{t+1}$, contributes the second term $\bs{\Delta}_{t+1}^{\top} \mb{q}_{t+1}^{}$.
It is worthwhile pointing out that if the uncertainty dependence occurs on the left-hand side (LHS) (e.g., 
$\widetilde{\mU}_t(\mb{d}_{t-1}) = \{\mb{d}_t \mid \mb{G}_t(\mb{d}_{t-1}) \mb{d}_t \geq \mb{g}_t\}$), then the corresponding reformulation requires dual variables that are functions of $\mb{d}_{t-1}$, resulting in  an-infinite dimensional optimization problem.

Theorems~\ref{prop:ellips_center},~\ref{prop:matrix_depen}, and ~\ref{prop_poly_rhs} highlight that connected uncertainties can be modeled in a natural and tractable fashion via \cus.
These results show that modifying decisions in order to account for future worst cases is an instrumental component of reformulating RO problems with \cus.
An extension to affinely adjustable RO problems reveals that the resulting constraints are similar in nature to those in Theorem~\ref{prop_poly_rhs} (see electronic companion for a detailed discussion).

When uncertainties follow unknown distributions, the concept of connected sets can  also be applied to multi-period distributionally robust problems. 
Moment-based ambiguity sets lend themselves naturally, where, for example, the mean or covariance depends on uncertainty realizations from the previous period.
In a parallel fashion to the RO setting, we now discuss modeling with \cus~in a distributional environment for static solutions.

%
%%%%%%%%%%%%%%%%%%%%%%%%%%%%%%%%%%%%%%%%%%%%%%%%%%%%%%%%%%%%%%%%%%%%%%%%%%%%%%%%%%%%%%%%%%%%%%%%%%%%%%%%

\section{Connected Uncertainty with DRO}
\label{sec:connect_sets}

In distributionally robust settings, the uncertainty set models all distributions that satisfy the set specifying constraints.
For \cus, the parameters of these constraints may depend on previous realizations.
Consider the following example of \cus~that are specified by distributional moments,
\begin{align*}
\widetilde{\mU}_1 &= \bigg\{P_1 \in \mathcal{M} \Big| P_1(\mb{d}_1 \in \Xi_1)
= 1,\;  \big|\mathbb{E}_{P_{1}}[\mb{d}_1] - \bs{\mu}_1\big| \leq
\bs{\delta}_1,\;\\
&\hspace{6cm}\mathbb{E}_{P_1}[(\mb{d}_1 - \bs{\mu}_1)(\mb{d}_1
- \bs{\mu}_1)^{\top}] \preceq \bs{\Sigma}_1  \bigg\},\\
\widetilde{\mU}_2(\mb{d}_{1}) &= \bigg\{P_{2 \mid {1}} \in \mathcal{M} \Big| P_{2\mid {1}}(\mb{d}_{2} \in \Xi_2) = 1,\; \big|\mathbb{E}_{P_{2\mid {1}}}[\mb{d}_2] - \bs{\mu}_2(\mb{d}_{1})\big| \leq
\bs{\delta}_2,\;\\
&\hspace*{60mm} \mathbb{E}_{P_{2\mid
		{1}}}[(\mb{d}_2 - \bs{\mu}_2^0)(\mb{d}_2 - \bs{\mu}_2^0)^{\top}] \preceq
\bs{\Sigma}_2 \bigg\},
\end{align*}
where $P_{2 \mid 1}$ describes the conditional distribution of $\mb{d}_2$ given the realization of $\mb{d}_{1}$. 
In these uncertainty sets, the conditional mean $\mathbb{E}_{P_{2|1}}[\mb{d}_2]$ is bounded by $\bs{\mu}_2(\mb{d}_{1}) + \bs{\delta}_2$ and $\bs{\mu}_2(\mb{d}_{1}) - \bs{\delta}_2$ and the covariance matrix by $\bs{\Sigma}_2$, with the bounds updated based on the previous realization. 
These  sets naturally describe settings in which the uncertainty is modeled using time series. 
In what follows, we consider distributional uncertainty sets to be \emph{connected}, when the moments at any given period depend on the realizations from previous periods, as presented in the example of $\widetilde{\mU}_1$ and $\widetilde{\mU}_2(\mb{d}_1)$.

%
%%%%%%%%%%%%%%%%%%%%%%%%%%%%%%%%%%%%%%%%%%%%%%%%%%%%%%
%
The aim of this section is to reformulate the constraint
\begin{equation}
\tag{C-DRO}
\label{eq:cdro}
\mathbb{E}_P[\sum_{t=1}^{T}h_t(\mb{x}_t,\mb{d}_t)] \leq B.
\end{equation}
The expectation \(\mathbb{E}_P[\cdot]\) is taken over the joint distribution $P$ of all \(\mb{d}_t\),
and \(h_t(\mb{x}_t,\mb{d}_t)\) is a function of the decision variable \(\mb{x}_t \in \mathbb{R}^{n_t}\) and the uncertain parameter \(\mb{d}_t \in \mathbb{R}^{m}\).
Unless specified, we do not make any assumptions on the structure of \(h_t(\cdot,\!\cdot)\) beyond regularity conditions required for the existence of integrals.
The dimension of $\mb{d}_t$ shall be constant over time for the clarity of exposition.
Each \(\mb{d}_t\) has a distribution that lies in a different uncertainty set, which depends on the previous realization \(\mb{d}_{t-1}\).

\noindent\textbf{Remark:} While DRO problems with connected uncertainty sets appear similar to robust Markov decision processes (MDPs), they are different in nature. 
In robust MDPs, the timing of decision-making and uncertainty realization alternates, whereas in DRO with \cus~all decisions (over the entire time horizon) are made at the very first period, anticipating future uncertainties and before any realization of the uncertainty occurs.
This is necessary because in many applications, long term decisions have to be taken without full knowledge of future uncertainties, e.g., in risk management settings.

In what follows, we first provide a general reformulation for the constraint~\eqref{eq:cdro} before showing that the tractability can be improved with a conservative approximation.

\subsection{Mean dependence}
\label{sec:dro_reform}
To reformulate~\eqref{eq:cdro} for moment based \cus, consider the sets in which the bounds on the first moment parameter depend on the previous realization as 

\(\bs{\mu}_{t}(\mb{d}_{t-1}) = \mb{A}_{t} \mb{d}_{t-1} + \mb{b}_{t}\), given by
\begin{equation}
\tag{D}
\begin{aligned}
\widetilde{\mU}_t(\mb{d}_{t-1}) = \bigg\{P_{t \mid {t-1}} \in \mathcal{M} \Big| P_{t\mid {t-1}}(\mb{d}_{t} \in \Xi_t) = 1,\; \big|\;\mathbb{E}_{P_{t\mid {t-1}}}[\mb{d}_t] - \bs{\mu}_t(\mb{d}_{t-1})| \leq \bs{\delta}_t,\; \label{eq:d}\\
\quad\quad\quad\quad\quad\quad\quad\quad\quad\;\;\quad\quad\quad \mathbb{E}_{P_{t\mid
		{t-1}}}[(\mb{d}_t - \bs{\mu}_t^0)(\mb{d}_t - \bs{\mu}_t^0)^{\top}] \preceq
\bs{\Sigma}_t \bigg\}.
\end{aligned}
\end{equation}
This set contains the distributions of \(\mb{d}_{t}\) conditioned on \(\bs{\mb{d}}_{t-1}\).
The parameter \(\bs{\mu}_t^0\) denotes a constant estimate of the mean and is different than the true mean $\mathbb{E}_{P_{t\mid {t-1}}}[\mb{d}_t]$, which is unknown and bounded by the set. 
To prevent the dependence of the second moment terms on the previous realization, we also assume \(\bs{\mu}_t^0\) to be different from \(\bs{\mu}_t(\mb{d}_{t-1})\).
A possible value for \(\bs{\mu}_t^0\) can be the unconditional mean of a time series at time \(t\). 
The robust counterpart of constraint~\eqref{eq:cdro} can be expressed with the following proposition.
\begin{proposition}
	\label{prop:first_reform}
	Given the sets \(\widetilde{\mU}_1,\dots,\widetilde{\mU}_T(\mb{d}_{T-1})\) and their joint uncertainty set \(\widetilde{\mU}\), we have the following.
	Constraint~\eqref{eq:cdro}, given by 
	\[\sup_{P \in \widetilde{\mU}} \mathbb{E}_P\left[\sum_{t=1}^{T}h_t(\mb{x}_t,\mb{d}_t)\right] \leq B\]
	is equivalent to
	\begin{equation}
	\label{eq:C1R}
	\begin{aligned}
	\sup_{P_1 \in \widetilde{\mU}_1} \mathbb{E}_{P_1}\!\bigg[h_1(\mb{x}_1,\mb{d}_1) 
	+& \sup_{P_{2 \mid 1} \in \widetilde{\mU}_2(\mb{d}_1)}\!\!\bigg\{\mathbb{E}_{P_{2 \mid 1}}\bigg[h_2(\mb{x}_2,\mb{d}_2) + \dots \\
	+& \!\!\!\!\!\!\!\sup_{P_{T \mid T-1} \in   \widetilde{\mU}_T(\mb{d}_{T-1})}\!\!\left\{\mathbb{E}_{P_{T\mid{T-1}}}
	\bigg[h_T(\mb{x}_T,\mb{d}_T)\right]\bigg\}\bigg]\bigg\}\bigg] \leq B.
	\end{aligned}
	\end{equation}	
\end{proposition} 

The proof of Proposition~\ref{prop:first_reform}  is provided in the electronic companion. 
To ease the exposition, consider the function
\begin{equation}
\label{eq:stat_resfunc_t}
g_t(\mb{x}_{[t+1:T]},\mb{d}_t) \coloneqq \!\!\!\!\!\!\!\!\sup_{P_{t+1|t} \in \widetilde{\mU}_{t+1}(\mb{d}_t)}\mathbb{E}_{P_{t+1|t}}[h_{t+1}(\mb{x}_{t+1},\mb{d}_{t+1}) + g_{t+1}(\mb{x}_{[t+2:T]},\mb{d}_{t+1})].
\end{equation}
For brevity, we denote \(\sigp_t
\equiv \bs{\Sigma}_t - \bs{\mu}_t^0(\bs{\mu}_t^0)^{\top}\), and use the compact notations \(\tp_t,\tuq_t,\tlq_t\) and \(\tR_t\) to denote variables \(p_t(\mb{d}_{t-1}), \mb{q}_t^u(\mb{d}_{t-1}), \mb{q}^l_t(\mb{d}_{t-1})\) and \(\mb{R}_t(\mb{d}_{t-1})\), which are functions of the previous period realization to ensure compactness. 
The following theorem provides the reformulation. 
%%%%%%%%%%%%%%%%%%%%%%%%%%%%%%%
\begin{theorem}
	\label{thm:dro1}
	The constraint~\eqref{eq:cdro} can be reformulated as 
	\begin{align*}
	p_1 + (\uq_1 - \lq_1)^{\top}\bs{\mu}_1 + (\uq_1 + \lq_1)^{\top}\bs{\delta}_1 + \mb{R}_1 \cdot \sigp_1 &\leq B\\
	\alpha_t(\mb{d}_{t-1},\mb{d}_t) + \bs{\beta}_t(\mb{d}_{t-1},\mb{d}_t)^{\top}\mb{d}_t + \mb{d}_t^{\top}\tR_t \mb{d}_t - h_t(\mb{x}_t,\mb{d}_t) &\geq 0\\
	&\hspace{1mm} \forall (\mb{d}_{t-1},\mb{d}_t) \in \Xi_{t-1}\times \Xi_t\;\;\forall t\\
	\tp_T + (\tuq_T - \tlq_T  - 2\tR_{T} \bs{\mu}_T^0)^{\top}\mb{d}_T + \mb{d}_T^{\top}\tR_T \mb{d}_T -h_T(\mb{x}_T,\mb{d}_T)
	&\geq 0\\
	&\hspace{1mm}\forall (\mb{d}_{T-1},\mb{d}_T) \in \Xi_{T-1}\times \Xi_T\\
	\tuq_t,\tlq_t \geq 0, \; \tR_t &\succeq 0 \;\;\forall \mb{d}_{t-1} \in \Xi_{t-1}\\ 
	\tuq_T,\tlq_T \geq 0, \; \tR_T &\succeq 0  \;\;\forall \mb{d}_{T-1} \in \Xi_{T-1},
	\end{align*}
	where $t = 1,\dots,T-1$ and
	\begin{eqnarray*}
		\alpha_i(\mb{d}_{i-1},\mb{d}_i) &=& \tp_i - \tp_{i+1} - (\tuq_{i+1} - \tlq_{i+1})^{\top}\mb{b}_{i+1} - \tR_{i+1} \cdot \sigp_{i+1} - (\tuq_{i+1} + \tlq_{i+1})^{\top}\bs{\delta}_{i+1}\\
		\bs{\beta}_i(\mb{d}_{i-1},\mb{d}_i) &=& \tuq_i - \tlq_i - 2\tR_i \bs{\mu}_i^0 - \mb{A}_{i+1}^{\top}\tuq_{i+1} + \mb{A}_{i+1}^{\top}\tlq_{i+1}.
	\end{eqnarray*}
\end{theorem}
\proof{Proof:}
	The proof proceeds by induction.
	We first provide the reformulation for \(t=1\), then assume it to be true for \(t=k\), before proving it for \(t=k+1\). 
	
	\paragraph{Base case (\(t = 1\))}
	The original constraint~\eqref{eq:cdro} can be expressed as 
	\begin{align*}
	\sup_{P_1 \in \widetilde{\mU}_1} \mathbb{E}_{P_1}[h_1(\mb{x}_1,\mb{d}_1) + \sup_{P_{2|1} \in \widetilde{\mU}_2(\mb{d}_1)}&\mathbb{E}_{P_{2|1}}[h_2(\mb{x}_2,\mb{d}_2) + \dots \\
	&+  \sup_{P_{T|T-1} \in \widetilde{\mU}_T(\mb{d}_{T-1})} \mathbb{E}[h_T(\mb{x}_T,\mb{d}_T)]]] \leq B,
	\end{align*}
	which can be compressed as 
	$\sup_{P_1 \in \widetilde{\mU}_1} \mathbb{E}_{P_1}[h_1(\mb{x}_1,\mb{d}_1) + g_1(\mb{x}_{[2:T]},\mb{d}_1)] \leq B$.	
	This optimization problem can be expressed as the following moment problem 
	\begin{align*}
	\sup_{P_1 \in \mathcal{M}} &\; \int_{\Xi_1} (h_1(\mb{x}_1,\mb{d}_1) + g_1(\mb{x}_{[2:T]},\mb{d}_1)) dP_1\\
	&\; \int_{\Xi_1} dP_1 = 1\\
	&\; \bs{\mu}_1 - \bs{\delta}_1 \leq \int_{\Xi_1} \mb{d}_1 dP_1 \leq \bs{\mu}_1 + \bs{\delta}_1\\
	&\; \int_{\Xi_1} (\mb{d}_1 - \bs{\mu}_1^0) (\mb{d}_1 - \bs{\mu}_1^0)^{\top}dP_1 \preceq \bs{\Sigma_1}.
	\end{align*}
	The last constraint is equivalent to
	\(\mathbb{E}_{P_{t \mid t-1}}[\mb{d}_t \mb{d}_t^{\top}] - 2\bs{\mu}_t^0\mathbb{E}_{P_{t \mid t-1}}[\mb{d}_t^{\top}]  \preceq \bs{\Sigma}_t - \bs{\mu}_t^0(\bs{\mu}_t^0)^{\top}
	\). 
	Thus the dual problem is given by 
	\begin{align*}
	\inf_{p_1,\uq_1,\lq_1,\mb{R}_1} &\; p_1 + (\uq_1)^{\top}(\bs{\mu}_1 + \bs{\delta}_1) - (\lq_1)^{\top}(\bs{\mu}_1 - \bs{\delta}_1) + \mb{R}_1 \cdot \sigp_1\\
	\text{s.t.} \;\;\;\;\;&\; p_1 + (\uq_1 - \lq_1)^{\top}\mb{d}_1 - 2 (\bs{\mu}_1^0)^{\top}\mb{R}_1 \mb{d}_1 + \mb{d}_1^{\top}\mb{R}_1 \mb{d}_1 \geq h_1(\mb{x}_1,\mb{d}_1) \\
	&\hspace{7cm} + g_1(\mb{x}_{[2:T]},\mb{d}_1)\;\; \forall \mb{d}_1 \in \Xi_1\\
	&\; \uq_1, \lq_1 \geq 0\\
	&\; \mb{R}_1 \succeq 0,
	\end{align*}
	which proves the base case.
	
	\paragraph{Inductive case (\(t=k\))}
	For the \(t=k+1\) reformulation, consider the function 
	\begin{equation}
	\label{eq:resfunc}
	g_k(\mb{x}_{[k+1:T]},\mb{d}_k) = \sup_{P_{k+1|k} \in \widetilde{\mU}_{k+1}(\mb{d}_{k})}\mathbb{E}_{P_{k+1|k}}[h_{k+1}(\mb{x}_{k+1},\mb{d}_{k+1}) + g_{k+1}(\mb{x}_{[k+2:T]},\mb{d}_{k+1})].
	\end{equation}
	Now, the third constraint in the reformulation of Theorem~\ref{thm:dro1} for \(t=k\) can be expressed as 
	\begin{align*}
	\tp_k &+ (\tuq_k - \tlq_k - 2\tR_k \bs{\mu}_k^0)^{\top}\mb{d}_k + \mb{d}_k^{\top}\tR_k \mb{d}_k \geq h_k(\mb{x}_k,\mb{d}_k) \\
	& \qquad \qquad + \sup_{P_{k+1|k} \in \widetilde{\mU}_{k+1}(\mb{d}_k)}\mathbb{E}_{P_{k+1|k}}[h_{k+1}(\mb{x}_{k+1},\mb{d}_{k+1}) + g_{k+1}(\mb{x}_{[k+2:T]},\mb{d}_{k+1})] \\
	&\hspace{8.7cm}\forall \mb{d}_{k-1} \in \Xi_{k-1}\;\;\forall \mb{d}_k \in \Xi_k.
	\end{align*}
	Using the dual of~\eqref{eq:resfunc},
	we can write the above constraint as the following two constraints
	\begin{align*}
	& \tp_k+ (\tuq_k - \tlq_k - 2\tR_k \bs{\mu}_k^0)^{\top}\mb{d}_k + \mb{d}_k^{\top}\tR_k \mb{d}_k \\
	& \qquad \geq h_k(\mb{x}_k,\mb{d}_k) +  \tp_{k+1} + \tR_{k+1} \cdot \sigp_{k+1} + (\tuq_{k+1} - \tlq_{k+1})^{\top}\bs{\mu}_{k+1}(\mb{d}_k) \\
	& \qquad + (\tuq_{k+1} + \tlq_{k+1})^{\top}\bs{\delta}_{k+1}\;\; &&\hspace{-1.8cm}\forall \mb{d}_{k-1} \in \Xi_{k-1}\;\; \forall \mb{d}_k \in \Xi_k\\
	&\tp_{k+1}+ (\tuq_{k+1} - \tlq_{k+1}- 2\tR_{k+1}\bs{\mu}_{k+1}^0)^{\top}\mb{d}_{k+1} + \mb{d}_{k+1}^{\top}\tR_{k+1} \mb{d}_{k+1} \\ 
	&\qquad  \geq h_{k+1}(\mb{x}_{k+1},\mb{d}_{k+1}) + g_{k+1}(\mb{x}_{[k+2:T]},\mb{d}_{k+1}) &&\hspace{-1.8cm}\forall \mb{d}_{k} \in \Xi_{k} \;\; \forall \mb{d}_{k+1} \in \Xi_{k+1}\\
	& \tuq_{k+1},\tlq_{k+1} \geq 0,\;\tR_{k+1} \succeq 0 &&\hspace{-1.8cm}\forall \mb{d}_{k} \in \Xi_{k}.
	\end{align*}
	Substituting \(\bs{\mu}_{k+1}(\mb{d}_k) = \mb{A}_{k+1} \mb{d}_k + \mb{b}_{k+1}\), we rearrange the first constraint as 
	\begin{align*}
	&p_k+ (\tuq_k - \tlq_k - 2\tR_k\bs{\mu}_k^0 - \mb{A}_{k+1}^{\top}\tuq_{k+1} + \mb{A}_{k+1}^{\top}\tlq_{k+1})^{\top}\mb{d}_k + \mb{d}_k^{\top}\tR_k \mb{d}_k\\ 
	& \qquad \geq h_k(\mb{x}_k,\mb{d}_k) + \tp_{k+1} + (\tuq_{k+1} - \tlq_{k+1})^{\top} \mb{b}_{k+1} + (\tuq_{k+1} + \tlq_{k+1})^{\top}\bs{\delta}_{k+1} \\ 
	& \qquad + \tR_{k+1} \cdot \sigp_{k+1}\hspace{5.6cm} \;\;\forall \mb{d}_{k-1} \in \Xi_{k-1}\;\; \forall \mb{d}_k \in \Xi_k, 
	\end{align*}
	which can be written in a more compact form as
	\begin{align*}
	&\alpha_k(\mb{d}_{k-1},\mb{d}_k)+ \bs{\beta}_k(\mb{d}_{k-1},\mb{d}_k)^{\top}\mb{d}_k + \mb{d}_k^{\top}\tR_k \mb{d}_k - h_k(\mb{x}_k,\mb{d}_k) \geq 0 \\
	&\hspace{9cm}\forall \mb{d}_{k-1} \in \Xi_{k-1}\;\; \forall \mb{d}_k \in \Xi_k.
	\end{align*}
	We can now give the complete set of constraints for $t=k+1$ as 
	\begin{align*}
	p_1 + (\uq_1 - \lq_1)^{\top}\bs{\mu}_1 +  (\uq_1 + \lq_1)^{\top}\bs{\delta}_1 + \mb{R}_1 \cdot \sigp_1 &\leq B \\
	\alpha_t(\mb{d}_{t-1},\mb{d}_t) + \bs{\beta}_t(\mb{d}_{t-1},\mb{d}_t)^{\top}\mb{d}_t + \mb{d}_t^{\top}\tR_t \mb{d}_t - h_t(\mb{x}_t,\mb{d}_t) &\geq 0 \\
	&\hspace{-1cm}\forall \mb{d}_{t-1} \in \Xi_{t-1}\; \forall \mb{d}_t \in \Xi_t ,\;\;\forall t\\
	\tp_{k+1}+ (\tuq_{k+1} - \tlq_{k+1} - 2\tR_{k+1}\bs{\mu}_{k+1}^0)^{\top}\mb{d}_{k+1} + \mb{d}_{k+1}^{\top}\tR_{k+1} \mb{d}_{k+1} &\geq h_{k+1}(\mb{x}_{k+1},\mb{d}_{k+1})\\ 
	&\hspace{-4.5cm} + g_{k+1}(\mb{x}_{[k+2:T]},\mb{d}_{k+1}) \;\; \forall \mb{d}_{k} \in \Xi_{k} \quad\forall \mb{d}_{k+1} \in \Xi_{k+1}.\\
	\tuq_t,\tlq_t &\geq 0 \quad \forall \mb{d}_{t-1} \in \Xi_{t-1} \;\forall t\\ 
	\tR_t &\succeq 0 \quad \forall \mb{d}_{t-1} \in \Xi_{t-1} \;\forall t,
	\end{align*}
	The complete reformulation of the Constraint~\eqref{eq:C1R} can be obtained by applying the induction up to \(t=T\), where 
	\[g_{T-1}(\mb{x}_{T},\mb{d}_{T-1}) = \sup_{P_{T \mid T-1} \in \widetilde{\mU}_T(\mb{d}_{T-1})}\mathbb{E}_{P_{T \mid T-1}} h_T(\mb{x}_T,\mb{d}_{T}).\] \Halmos
\endproof
%%%%%%%%%%%%%%%%%%%%%%%%%%%%%%%
%
Theorem~\ref{thm:dro1} provides a prescription for how to modify the protection term when uncertainties are connected.
Notice that the variables in Theorem~\ref{thm:dro1} depend on the uncertainty realization in the previous period.
Furthermore, the variables $\alpha_i$ and $\bs{\beta}_i$  conjoin the variables of the current period with that of the future period.
Therefore, Theorem~\ref{thm:dro1} modifies the protection in order to account for the dependence amongst the uncertainties.
That means when, for example, autocorrelated uncertainties are modeled in the traditional (nonconnected) way, the protection terms are not appropriately modified, which can lead to constraint violations.

\subsection{Conservative Reformulation}
\label{sec:dro_approx}  

Since the reformulation in Theorem~\ref{thm:dro1} is an infinite dimensional optimization problem, which can be computationally challenging, the following theorem provides a conservative reformulation that is tractable.

%%%%%%%%%%%%%%%%%%%%%%%%%%%%%%%%%%%%%%%%%%%%%%%%%%%%%%%%%%%%%%%%%%%%%%%%%%%%%%
\begin{theorem}
	\label{thm:stat_dro}
	The Constraint~\eqref{eq:cdro} can be conservatively reformulated as 
	\begin{align*}
	p_1 + (\uq_1 - \lq_1)^{\top}\bs{\mu}_1 + (\uq_1 + \lq_1)^{\top}\bs{\delta}_1 +  \mb{R}_1 \cdot \sigp_1 &\leq B&\\
	\alpha_t + \bs{\beta}_t^{\top}\mb{d}_t + \mb{d}_t^{\top}\mb{R}_t \mb{d}_t - h_t(\mb{x}_t,\mb{d}_t) &\geq 0 &&\hspace{-8mm}\forall \mb{d}_t \in \Xi_t \;\;\forall t = 1,\dots,T-1\\
	p_T + (\uq_T - \lq_T - 2\mb{R}_{T} \bs{\mu}_T^0)^{\top}\mb{d}_T + \mb{d}_T^{\top}\mb{R}_T \mb{d}_T &\geq h_T(\mb{x}_T,\mb{d}_T) &&\forall \mb{d}_T \in \Xi_T\\
	\uq_t,\lq_t &\geq 0, \; \mb{R}_t \succeq 0  &&\forall t=1,\dots,T,
	%	\mb{R}_t &\succeq 0 &&\forall t=1,\dots,T,
	\end{align*}
	where 
	\begin{eqnarray*}
		\alpha_i &=& p_i - p_{i+1} - (\uq_{i+1} - \lq_{i+1})^{\top}\mb{b}_{i+1} - \mb{R}_{i+1} \cdot \sigp_{i+1} - (\uq_{i+1} + \lq_{i+1})^{\top}\bs{\delta}_{i+1}\\
		\bs{\beta}_i &=& \uq_i - \lq_i - 2\mb{R}_{i} \bs{\mu}_i^0 - \mb{A}_{i+1}^{\top}\uq_{i+1} + \mb{A}_{i+1}^{\top}\lq_{i+1}.
	\end{eqnarray*}
\end{theorem}
\proof{Proof:}
	The proof proceeds by induction.
	
	\paragraph{Base case (\(t = 1\))}
	The original Constraint~\eqref{eq:cdro} can be expressed as 
	\begin{align*}
	\sup_{P_1 \in \widetilde{\mU}_1} \mathbb{E}_{P_1}[h_1(\mb{x}_1,\mb{d}_1) + \!\!\!\!\!\sup_{P_{2|1} \in \widetilde{\mU}_2(\mb{d}_1)}\mathbb{E}_{P_{2|1}}[h_2(\mb{x}_2,\mb{d}_2) &+ \dots\\
	&+  \!\!\!\!\!\sup_{P_{T|T-1} \in \widetilde{\mU}_T(\mb{d}_{T-1})} \mathbb{E}[h_T(\mb{x}_T,\mb{d}_T)]]] \leq B,
	\end{align*}
	which can be shortened to 
	\(\sup_{P_1 \in \widetilde{\mU}_1} \mathbb{E}_{P_1}[h_1(\mb{x}_1,\mb{d}_1) + g_1(\mb{x}_{[2:T]},\mb{d}_1)] \leq B\), and rewritten as
	\begin{align*}
	\sup_{P_1 \in \mathcal{M}} &\; \int_{\Xi_1} \left(h_1(\mb{x}_1,\mb{d}_1) + g_1(\mb{x}_{[2:T]},\mb{d}_1)\right) dP_1\\
	&\; \int_{\Xi_1} dP_1 = 1\\
	&\; \bs{\mu}_1 - \bs{\delta}_1 \leq \int_{\Xi_1} \mb{d}_1 dP_1 \leq \bs{\mu}_1 + \bs{\delta}_1\\
	&\; \int_{\Xi_1} (\mb{d}_1 - \bs{\mu}_1^0) (\mb{d}_1 - \bs{\mu}_1^0)^{\top}dP_1 \preceq \bs{\Sigma_1}.
	\end{align*}
	The dual of this moment problem proves the base case via
	\begin{align*}
	\inf_{p_1,\uq_1,\lq_1,\mb{R}_1} &\; p_1 + (\uq_1)^{\top}(\bs{\mu}_1 + \bs{\delta}_1) - (\lq_1)^{\top}(\bs{\mu}_1 - \bs{\delta}_1) + \mb{R}_1 \cdot \sigp_1\\
	\text{s.t.} \;\;\;\;\;&\; p_1 + (\uq_1 - \lq_1)^{\top}\mb{d}_1 - 2 (\bs{\mu}_1^0)^{\top}\mb{R}_1 \mb{d}_1 + \mb{d}_1^{\top}\mb{R}_1 \mb{d}_1 \geq h_1(\mb{x}_1,\mb{d}_1) + g_1(\mb{x}_{[2:T]},\mb{d}_1) \\
	&\hspace{10cm}\forall \mb{d}_1 \in \Xi_1\\
	&\; \uq_1, \lq_1 \geq 0, \; \mb{R}_1 \succeq 0.
	\end{align*}
	
	\paragraph{Inductive case (\(t=k\)).} 
	We assume the constraints in Theorem~\ref{thm:stat_dro} hold for \(t=k\) and prove the reformulation for \(t=k+1\). 
	These constraints can be expressed as 
	\begin{align*}
	p_1 + (\uq_1 - \lq_1)^{\top}\bs{\mu}_1 +  (\uq_1 + \lq_1)^{\top}\bs{\delta}_1 + \mb{R}_1 \cdot \sigp_1 &\leq B \\
	\alpha_t + \bs{\beta}_t^{\top}\mb{d}_t + \mb{d}_t^{\top}\mb{R}_t \mb{d}_t - h_t(\mb{x}_t,\mb{d}_t) &\geq 0 \hspace{3.3cm} \forall \mb{d}_t \in \Xi_t ,\;\;\forall t\\
	p_k+ (\uq_k - \lq_k - 2\mb{R}_{k} \bs{\mu}_k^0)^{\top}\mb{d}_k + \mb{d}_k^{\top}\mb{R}_k \mb{d}_k &\geq h_k(\mb{x}_k,\mb{d}_k) 
	+ g_k(\mb{x}_{[k+1:T]},\mb{d}_k)\\
	&\hspace{4cm} \forall \mb{d}_k \in \Xi_k,
	\end{align*}
	where $t= 1,\dots,k-1$ and \(\alpha_t\) and  \(\bs{\beta}_t\) are given in the theorem.
	For \(t=k+1\), 
	consider the function 
	\begin{equation*}
	\label{eq:stat_resfunc}
	g_k(\mb{x}_{[k+1:T]},\mb{d}_k) = \sup_{P_{k+1|k} \in \widetilde{\mU}_{k+1}(\mb{d}_k)}\mathbb{E}_{P_{k+1|k}}[h_{k+1}(\mb{x}_{k+1},\mb{d}_{k+1}) + g_{k+1}(\mb{x}_{[k+2:T]},\mb{d}_{k+1})],
	\end{equation*}
	which is bounded above by its dual for any feasible \((p_{k+1},\uq_{k+1},\lq_{k+1},\mb{R}_{k+1}) \in \mathcal{P}_{k+1}\) as 
	\begin{align*}
	g_k(\mb{x}_{[k+1:T]},\mb{d}_k) \leq  p_{k+1} + \mb{R}_{k+1} \cdot \sigp_{k+1} \\&\hspace{-2.8cm}+ (\uq_{k+1} - \lq_{k+1})^{\top}\bs{\mu}_{k+1}(\mb{d}_k) + (\uq_{k+1} + \lq_{k+1})^{\top}\bs{\delta}_{k+1}
	&\forall \mb{d}_k \in \Xi_k.
	\end{align*}
	Then, the last constraint in the reformulation for \(t=k\) can be conservatively expressed as 
	\begin{equation}
	\label{eq:conservative_approx_step}
	\begin{aligned}
	&p_k + (\uq_k - \lq_k - 2\mb{R}_{k} \bs{\mu}_k^0)^{\top}\mb{d}_k + \mb{d}_k^{\top}\mb{R}_k \mb{d}_k\\ 
	&\qquad \geq h_k(\mb{x}_k,\mb{d}_k) +  p_{k+1} + \mb{R}_{k+1} \cdot \sigp_{k+1} +   
	(\uq_{k+1} - \lq_{k+1})^{\top}\bs{\mu}_{k+1}(\mb{d}_k) \\
	& \qquad + (\uq_{k+1} + \lq_{k+1})^{\top}\bs{\delta}_{k+1} \;\; &&\hspace{-26mm} \forall \mb{d}_k \in \Xi_k  \\
	& p_{k+1}+ (\uq_{k+1} - \lq_{k+1} - 2\mb{R}_{k+1} \bs{\mu}_{k+1}^0)^{\top}\mb{d}_{k+1} + \mb{d}_{k+1}^{\top}\mb{R}_{k+1} \mb{d}_{k+1} \\ 
	& \qquad \geq h_{k+1}(\mb{x}_{k+1},\mb{d}_{k+1}) +  g_{k+1}(\mb{x}_{[k+2:T]},\mb{d}_{k+1}) \;\; &&\hspace{-26mm}  \forall \mb{d}_{k+1} \in \Xi_{k+1}.
	\end{aligned}
	\end{equation}
	Using \(\bs{\mu}_{k+1}(\mb{d}_k) = \mb{A}_{k+1} \mb{d}_k + \mb{b}_{k+1}\), we can rearrange the first constraint as 
	\begin{align*}
	& p_k+ (\uq_k - \lq_k - 2\mb{R}_{k} \bs{\mu}_k^0 - \mb{A}_{k+1}^{\top}\uq_{k+1} + \mb{A}_{k+1}^{\top}\lq_{k+1})^{\top}\mb{d}_k + \mb{d}_k^{\top}\mb{R}_k \mb{d}_k \\
	& \hspace{2mm} \geq h_k(\mb{x}_k,\mb{d}_k)  + p_{k+1} + (\uq_{k+1} - \lq_{k+1})^{\top} \mb{b}_{k+1} + (\uq_{k+1} + \lq_{k+1})^{\top}\bs{\delta}_{k+1}  
	+ \mb{R}_{k+1} \cdot \sigp_{k+1} \\
	&\hspace{10.5cm} \forall \mb{d}_k \in \Xi_k,
	\end{align*}
	which can be written compactly as
	$\alpha_k+ \bs{\beta}_k^{\top}\mb{d}_k + \mb{d}_k^{\top}\mb{R}_k \mb{d}_k - h_k(\mb{x}_k,\mb{d}_k) \geq 0 \;\;\forall \mb{d}_k \in \Xi_k$.
	With this, the complete set of constraints for $t=k+1$ are 
	\begin{align*}
	p_1 + (\uq_1 - \lq_1)^{\top}\bs{\mu}_1 +  (\uq_1 + \lq_1)^{\top}\bs{\delta}_1 + \mb{R}_1 \cdot \sigp_1 &\leq B \\
	\alpha_t + \bs{\beta}_t^{\top}\mb{d}_t + \mb{d}_t^{\top}\mb{R}_t \mb{d}_t - h_t(\mb{x}_t,\mb{d}_t) &\geq 0\\
	&\hspace{-1cm}\forall \mb{d}_{t-1} \in \Xi_{t-1}\;\; \forall \mb{d}_t \in \Xi_t ,\;\;\forall t\\
	p_{k+1}+ (\uq_{k+1} - \lq_{k+1} - 2\mb{R}_{k+1} \bs{\mu}_{k+1}^0)^{\top}\mb{d}_{k+1} + \mb{d}_{k+1}^{\top}\mb{R}_{k+1} \mb{d}_{k+1} &\geq h_{k+1}(\mb{x}_{k+1},\mb{d}_{k+1}) \\
	&\hspace{-2.5cm}+ g_{k+1}(\mb{x}_{[k+2:T]},\mb{d}_{k+1}) \;\; \forall \mb{d}_{k+1} \in \Xi_{k+1}.
	\end{align*}
	The complete reformulation of the robust counterpart of~\eqref{eq:cdro} can be obtained by repeating the induction up to \(t=T\), using 
	\[g_{T-1}(\mb{x}_{T},\mb{d}_{T-1}) = \sup_{P_{T \mid T-1} \in \widetilde{\mU}_T(\mb{d}_{T-1})}\mathbb{E}_{P_{T \mid T-1}} h_T(\mb{x}_T,\mb{d}_{T}).\] \Halmos
\endproof
%%%%%%%%%%%%%%%%%%%%

In the reformulation of Theorem~\ref{thm:stat_dro}, $\alpha_i$ and $\bs{\beta}_i$ conjoin the variables of the current and the future periods, which is analogous to the results of Theorems~\ref{prop:ellips_center} and~\ref{thm:dro1}.
This modification of the protection term accounts for the connection between the periods. 
When the support set is an ellipsoid and the objective function is piecewise linear concave in \(\mb{d}_t\), it is possible to use the S-Lemma to reformulate the constraints of the problems as semi-definite constraints, which can then be solved by commercially available solvers. 

{\bf Remark:}
Although Theorem~\ref{thm:stat_dro} provides a conservative reformulation of the Constraint~\eqref{eq:C1R}, under certain conditions, it is also possible for this reformulation to be exact. 
This tightness occurs, 
when $g_{k-1}(\mb{x}_{[k:T]},\mb{d}_{k-1})$ is a linear function of the uncertain component $\mb{d}_{k-1}$ over the set $\Xi_{k-1}$.

The linearity arises when there exists a single optimal solution for all $\mb{d}_{k-1}$.
Note that the verification of this condition requires solving a sequence of bilinear optimization problems, which is challenging.

\subsection{Parameter Estimation}
\label{sec:uncertainty_estimation}
In the CU models for both RO and DRO problems, the connection between the uncertainties at different time periods is modeled as a linear function of the parameters. 
This enables us to leverage existing parameter estimation approaches to establish the model parameters.
For example, when data exhibits a dependence on the center across different time periods, an autoregressive analysis can reveal the parameters of the model.
On the other hand, when the variance of the data exhibits a changing nature over time, all centers can be preestimated by regression models, and the residuals can be fitted to the general form in~\eqref{eq:mu} with linear regressions.

In Section~\ref{sec:dro_example}, we numerically demonstrate the advantages of \cus~through a portfolio-optimization problem, where the realized wealth exhibits a lower standard deviation, supporting the notion of robust decisions.
In the next two sections, we probe the performance of the proposed \cus~modeling frameworks through two stylized RO and DRO problems.
Each demonstration serves to model a broad range of decision-making settings.
For a clean experiment, randomly generated data are used  to directly relate the findings to the properties of our proposed models.

%%%%%%%%%%%%%%%%%%%%%%%%%%%%%%%%%%%%%

\section{RO Application: Knapsack Problem}
\label{sec:appl-3:-knaps}
Knapsack problems offer an insightful benchmark for optimization problems because they arise as sub-problems in many common applications~\citep{freville2004multidimensional,monaci2013robust}.
Consider a knapsack problem with known objective but uncertain constraint coefficients (a.k.a. weights).
The problem spans two periods, each with an uncertain weight, where the second period weight value depends on the first realized weight. 
These uncertain parameters arise from an auto-regressive model. 
Such a problem can be described by
\begin{equation}
\label{eq:ks}
\tag{KS}
\begin{aligned}
\max_{\mb{x}_1,\mb{x}_2}\;&\; \mb{c}_1^{\top}\mb{x}_1 + \mb{c}_2^{\top}\mb{x}_2\\
\text{s.t.}\;&\; \mb{d}_1^{\top}\mb{x}_1 + \mb{d}_2^{\top}\mb{x}_2 \leq B 
\quad\quad\quad \forall\mb{d}_2\in\mU_2(\mb{d}_1),\; \forall\mb{d}_1\in\mU_1 \\
&\; \mb{x}_1 \in \{0,1\}^{m_1}, \; \mb{x}_2 \in \{0,1\}^{m_2}.
\end{aligned}
\end{equation}
The two binary decisions \(\mb{x}_1\) and \(\mb{x}_2\)  with known objective \(\mb{c}_1,\mb{c}_2\) and uncertain weight coefficients \(\mb{d}_1,\mb{d}_2\) correspond to periods one and two.
Both decisions are taken before either of the weights are realized, i.e., \(\mb{x}_1\) and \(\mb{x}_2\) are here-and-now decisions.
The uncertain parameters $\mb{d}_t$ for $t=1,2$ are modeled to reside in ellipsoidal sets with given covariance matrices and known first period center $\bs{\mu}_1$. 
The uncertainty dependence can be modeled by allowing the center of the second period ellipsoid \(\bs{\mu}_2(\mb{d}_1)\) to depend on the realization of the first period weights $\mb{d}_1$ as 
\begin{equation*}
\label{eq:cond_mean}
\bs{\mu}_2(\mb{d}_1) = \bs{\Phi}\bs{\mu}_1 + \bs{\Psi}\mb{d}_1.
\end{equation*}
Here, the parameters \(\bs{\Phi}\) and \(\bs{\Psi}\) are matrices capturing the dependence  and which may be determined through time series modeling.
We emphasize that the solutions to Problem~\eqref{eq:ks} are static, while the uncertainties are connected and, as such stagewise.
This uncertainty model is parallel to the discussion in Section~\ref{sec:ellips_center}.

%%%%%%%%%%%%%%%%%%%%%%%%%%%%%%%%%%%%%%%%%%%%%%%%%%%%%%
Given the mean of either period \(\bs{\mu}_1\) and \(\bs{\mu}_2(\mb{d}_1)\), the residual uncertainties \(\bs{\epsilon}_1\) and \(\bs{\epsilon}_2\) in
\vspace{-2mm}
\begin{equation}
\label{eq:uncertain_realn}
\begin{aligned}
\mb{d}_1 = \bs{\mu}_1 + \bs{\epsilon}_1\;\text{and } \;
\mb{d}_2 = \bs{\Phi}\bs{\mu}_1 + \bs{\Psi}\mb{d}_1 + \bs{\epsilon}_2
\end{aligned}
\vspace{-2mm}
\end{equation}
are characterized by a normal distribution with mean \(\mb{0}\) and covariance \(\bs{\Sigma}\) such that \(\mb{L}\mb{L}^{\top} = \bs{\Sigma}\). 
Then, the corresponding \emph{first period} uncertainty set is
\vspace{-2mm}
\begin{equation*}
\mU_1 = \left\{\mb{d}_1 \mid \mb{d}_1 = \bs{\mu}_1 + \mb{L}\mb{u} \;:\;\|\mb{u}\|_2 \leq r_1\right\}.
\vspace{-2mm}
\end{equation*}
The second period mean given the first-period realization \(\bs{\mu}_2(\mb{d}_1)\) 
has the same covariance matrix.
Consequently, the \emph{second period} \cu~set is
\vspace{-2mm}
\begin{equation*}
\mU_2({\mb{d}_1}) = \left\{\mb{d}_2 \mid \mb{d}_2 =  \bs{\Phi}\bs{\mu}_1 + \bs{\Psi}\mb{d}_1 + \mb{L} \mb{w}, \;\;\|\mb{w}\|_2 \leq r_2 \right\}.
\vspace{-2mm}
\end{equation*}

In the experiment, we compare the performance of the connected model to the standard RO model, namely, when uncertainty dependence is not taken into account. 
The first period has the same uncertainty set \(\mU_1\).
For the second-period set, the parameters are \(\bs{\mu}_2 = \bs{\mu}_1 \) and \(\mb{L}\mb{L}^{\top} = \bs{\Sigma_2} = \bs{\Sigma}\).
Therefore, the \emph{nonconnected} (NC) second-period set is
\vspace{-2mm}
\begin{equation*}
\mU_{2,\text{NC}} = \left\{\mb{d}_2 \mid\mb{d}_2 = \bs{\mu}_2 + \mb{L} \mb{w}, \;\;\|\mb{w}\|_2 \leq r_2 \right\}.
\vspace{-2mm}
\end{equation*}
Note that the key difference between $\mathcal{U}_2({\mb{d}_1})$ and $\mathcal{U}_{2,\text{NC}}$
is that the center of \(\mU_2(\mb{d}_1)\) is \(\bs{\Phi}\bs{\mu}_1 + \bs{\Psi}\mb{d}_1\), whereas that of \(\mU_{2,\text{NC}}\) is \(\bs{\mu}_2\)--that is, 
not updated according to the realization of $\mb{d}_1$.
The covariance matrices for both remain the same. 
We now describe the experimental setting.

\subsubsection*{Numerical Experiments:}
\label{sec:numer-exper-kpbu}
Four experimental modules are conducted for uncertain $\mb{d}_1$, and $\mb{d}_2$, which are generated by~\eqref{eq:uncertain_realn} for fixed values of \(\bs{\Phi} = \mb{I}\) and \(\bs{\Psi} = \lambda \mb{I}\) with the identity matrix \(\mb{I}\) and a scalar \(\lambda\). 
The parameter \(\bs{\mu}_1\) is set to \(\mb{e}\), and \(\bs{\Sigma}\) is generated randomly. 
We consider a problem with \(20\) items in each period \(m_1 = m_2 = 20\), following these steps. 
\begin{enumerate}[label={(\roman*)}]\setlength{\itemindent}{5mm}
	\item Generate $k$ samples of \(\mb{c}_1, \mb{c}_2\) and $k$ estimates of \(\bs{\mu}_1\) and \(\bs{\Sigma}\), each using $l$ samples of \(\mb{d}_1\).
	\item For each estimate, solve~\eqref{eq:ks} for $s$ different uncertainty set sizes, yielding $k\times s$ solutions with $\lambda = 0.5$ and for a fixed set size $r_1 = r_2 = 2$ with changing $\lambda$.
	\item For each solution, generate $n$ samples of $\mb{d}_1, \mb{d}_2$ from~\eqref{eq:uncertain_realn} to probe constraint satisfaction.
	\item Average over $k\times s$ objective values and $k\times s$ fractions of constraint satisfaction.
\end{enumerate}
In module (i), let $l=500$ and $k=30$. 
Sample averages are used to estimate \(\bs{\mu}_1\) and \(\bs{\Sigma}\).
The coefficients $\mb{c}_1$ and $\mb{c}_2$ are sampled from normal distributions centered at  \(\mb{e}\) and \(1.25 \times \mb{e}\), respectively with covariance matrix \(\frac{1}{100}\bs{\Sigma}\). 
For the constraint RHS, the budget is \(B=40\).

\noindent
In (ii), the experiment is conducted for $s=20$ different uncertainty set sizes
\(r_1 = r_2 \in [0, 4]\) and for another 20 different correlation values with $\lambda \in [-1, 1]$. 

\noindent
In (iii), the constraint satisfaction is measured using $n=500$ samples. 

To probe each setting, the knapsack problem is solved for both CU and NC sets.
This is parallel to the difference between models in Figure~\ref{fig:set_pic}(b) and (d) of the introductory Example 2.
\begin{figure}[h!]
	\centering
	\begin{minipage}[t]{65mm}
		\includegraphics[width=65mm]{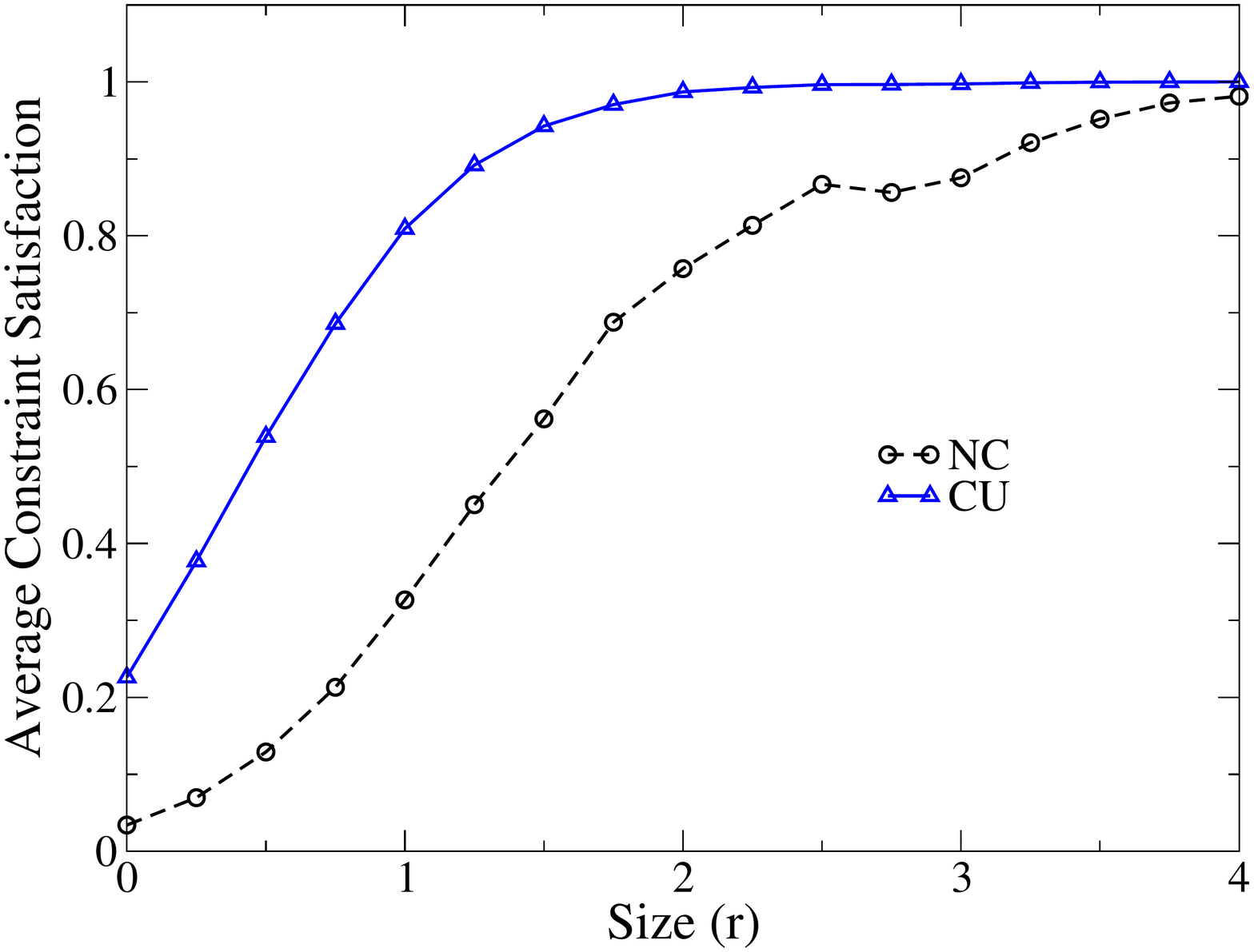}\\
	\end{minipage}
	\hspace{-5mm}
	\begin{minipage}[t]{65mm}	
		\includegraphics[width=65mm]{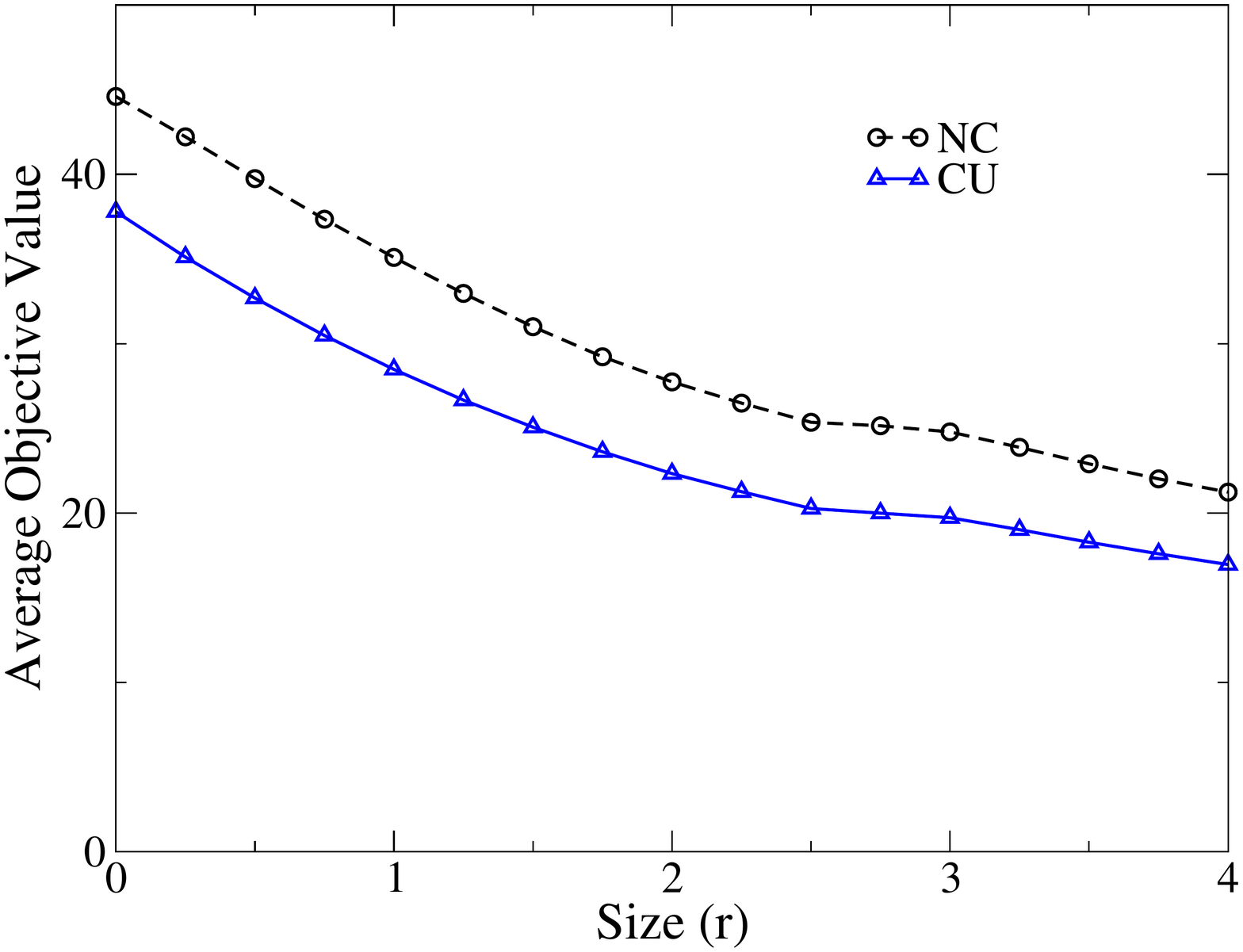}
	\end{minipage}
	\caption{Comparison of Connected and Nonconnected sets for the Robust Knapsack Problem at different set sizes:
		Fraction of Constraint Satisfaction (Left) and 
		Objective Value (Right)}	
	\label{fig:ts_05_cs_ao}
\end{figure}

\begin{figure}[h!]
	\centering
	\begin{minipage}[t]{65mm}
		\includegraphics[width=65mm]{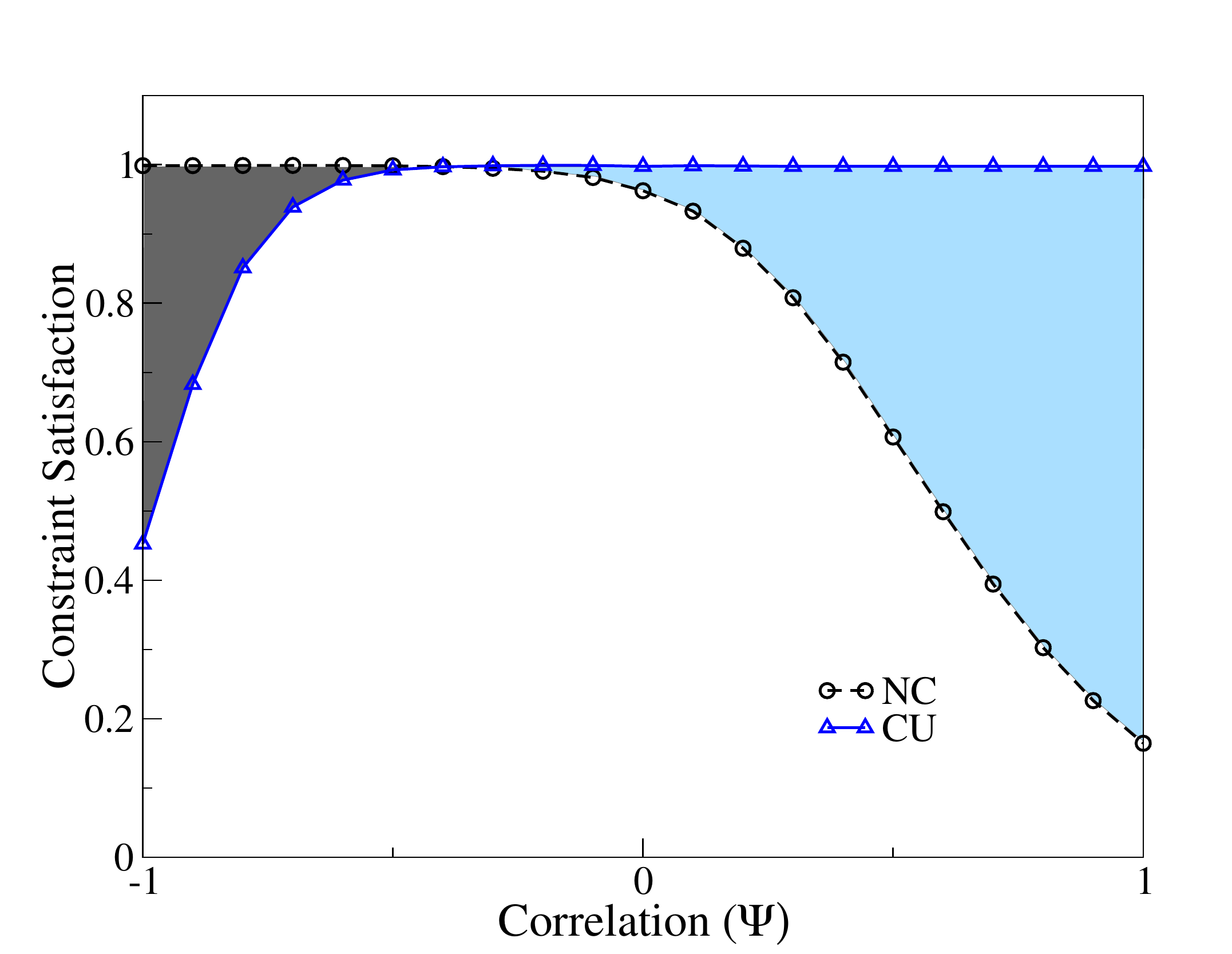}\\
	\end{minipage}
	\hspace{-5mm}
	\begin{minipage}[t]{65mm}	
		\includegraphics[width=65mm]{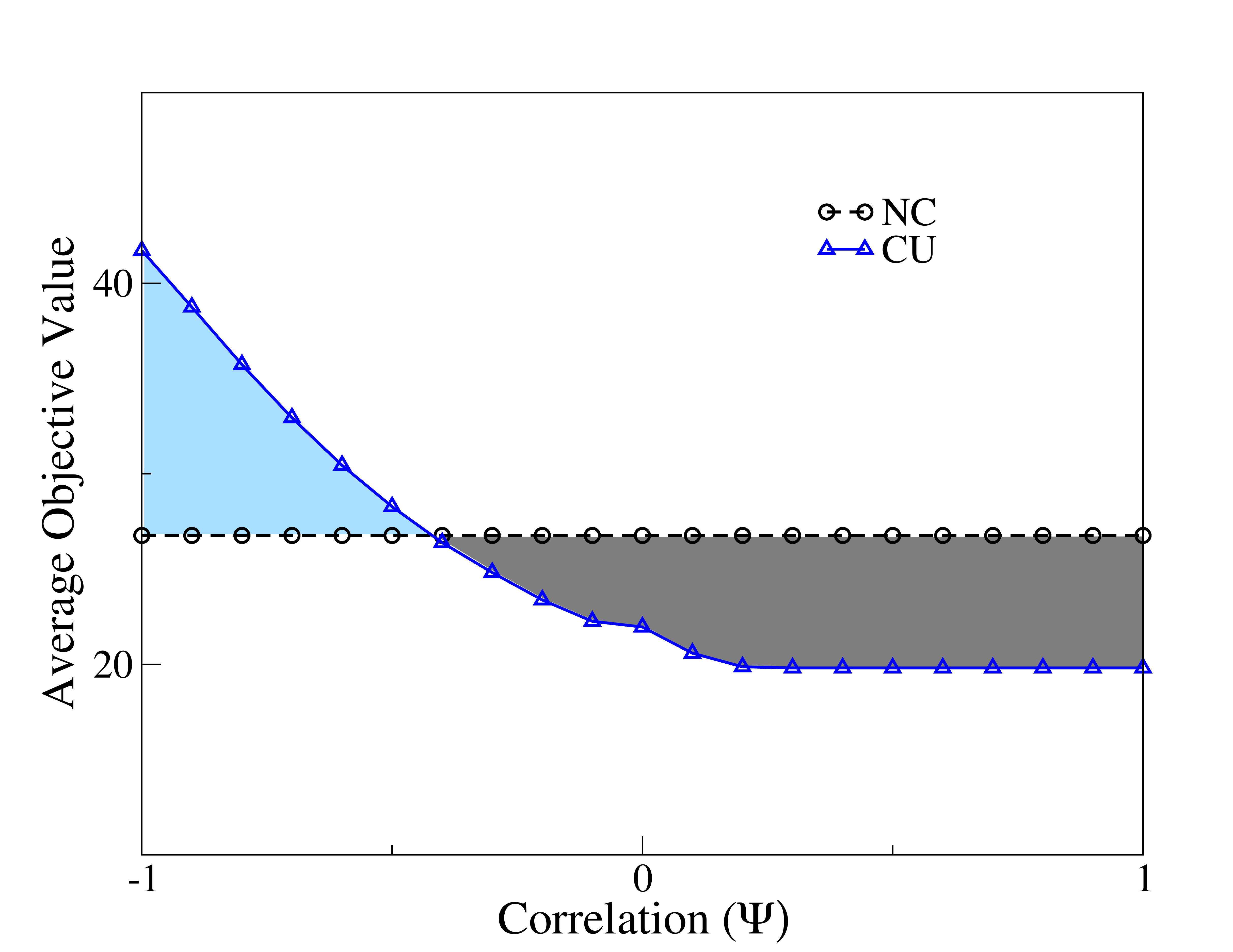}
	\end{minipage}
	\caption{Comparison of Connected and Nonconnected sets for the Robust Knapsack Problem at different temporal correlations (dark - NC Preferable and light - CU Preferable):
		Fraction of Constraint Satisfaction (Left) and 
		Objective Value (Right)}	
	\label{fig:ijoo_cs_ao}
\end{figure}
\begin{figure}[h!]
	\centering
	\vspace{-5mm}
	\includegraphics[width=65mm]{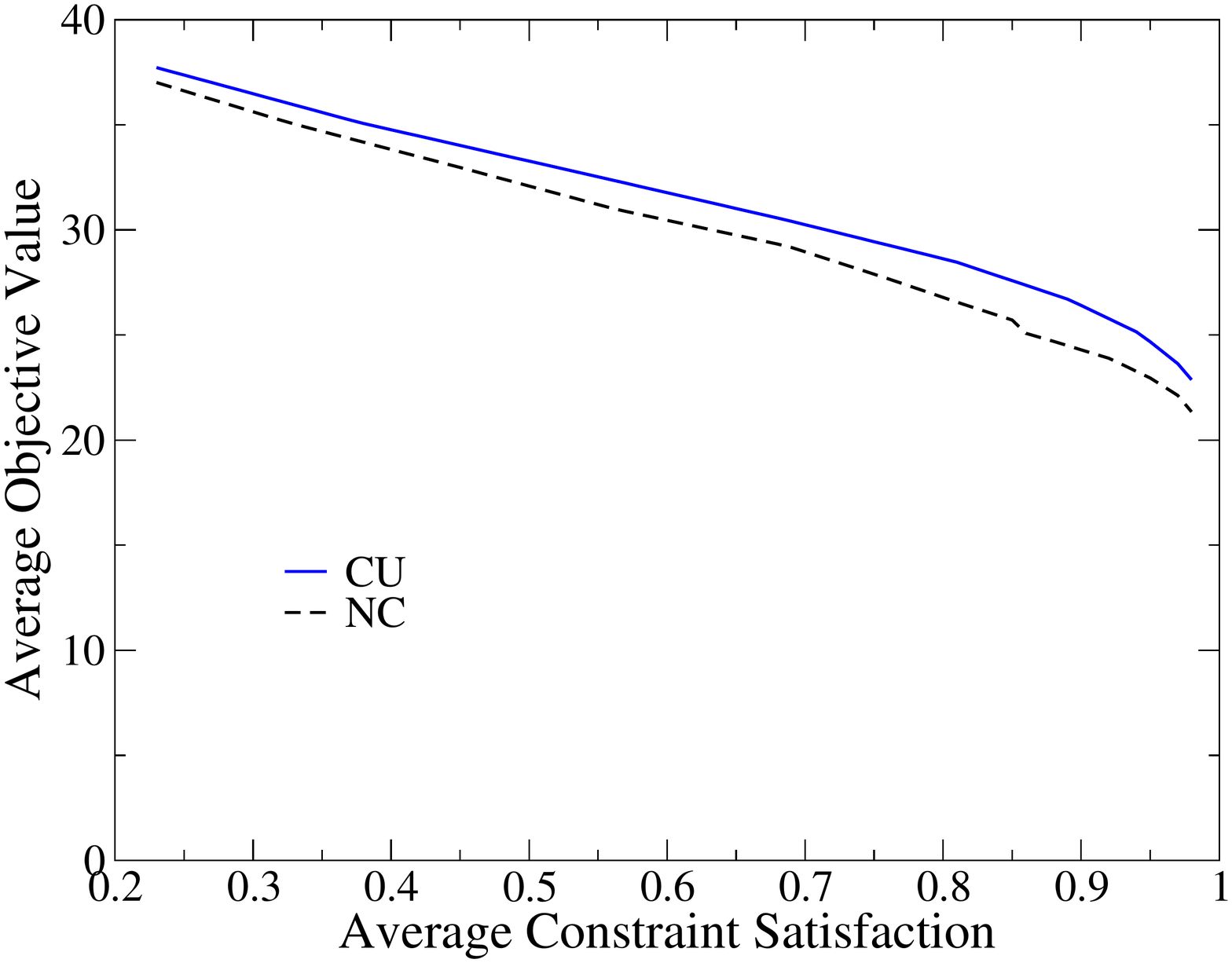}\hspace{-2mm}
	\includegraphics[width=65mm]{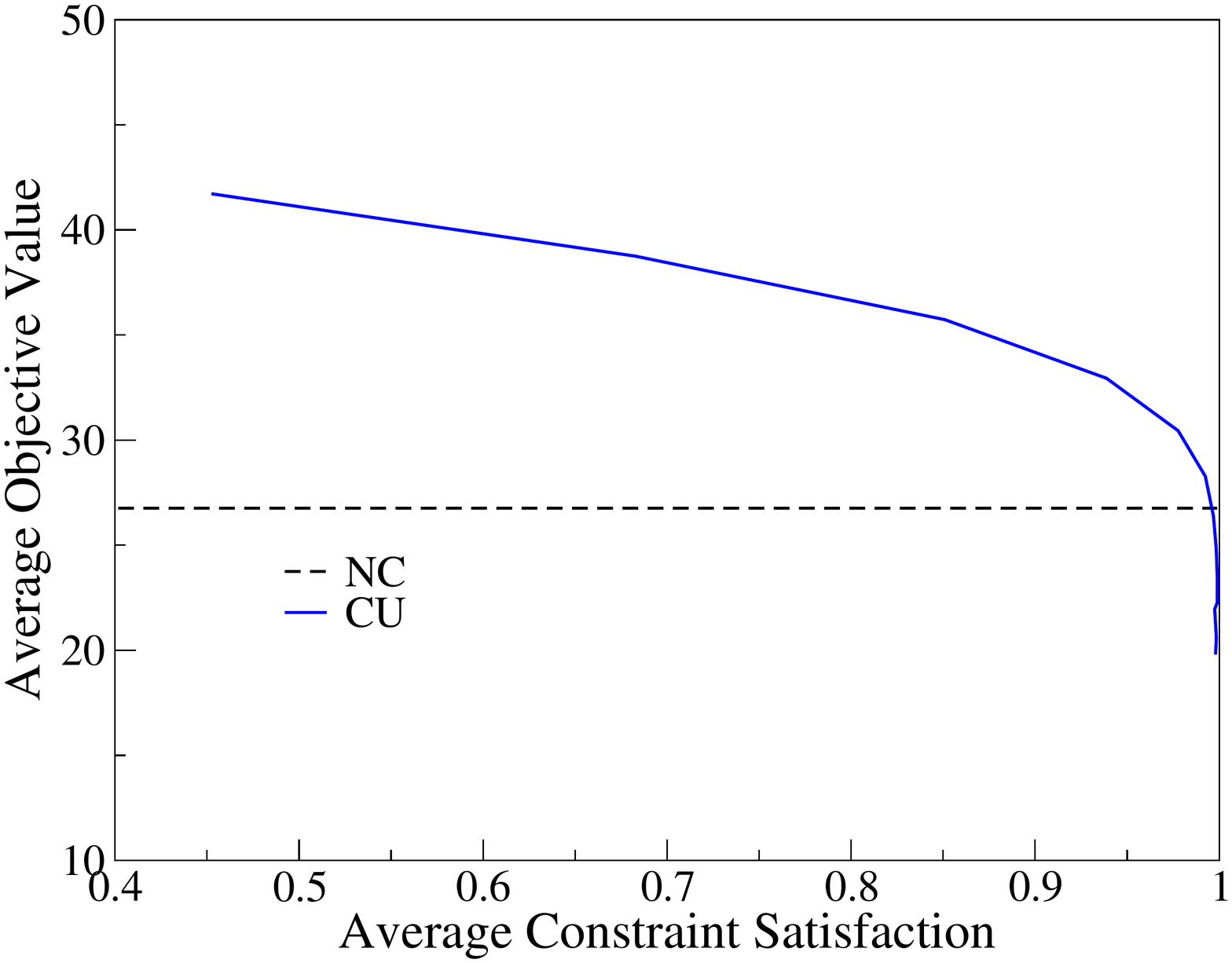}
	\caption{\label{fig:ts_05_aocs_sn}Comparison of Objective vs.  Constraint Satisfaction for Connected and Nonconnected sets:
		For Changing Uncertainty Set Size $r$ (Left) and For Changing Temporal Correlation $\Psi$ (Right)}
\end{figure}

Figure~\ref{fig:ts_05_cs_ao} shows the average constraint satisfaction (left) and the average objective value (right) for a varying size \(r\) of the uncertainty set.
The data points at \(r_1 = r_2 = 0\) correspond to the nominal problem with no uncertainty. 
Figure~\ref{fig:ijoo_cs_ao} presents the average constraint satisfaction (left) and the average objective value (right) as a function of the temporal correlation $\Psi$ (specifically $\lambda$). 
The shaded regions indicate where the different models CU or NC dominate the other on the basis of the constraint satisfaction or objective function value. 	
Figure~\ref{fig:ts_05_aocs_sn} shows the average objective value for any level of constraint satisfaction for different uncertainty set sizes (left) and for different temporal correlation values (right). 

We now summarize our observations.
\begin{itemize}
	\item\textbf{\emph{Effect of Uncertainty Set Size}}: 
	For both models, constraint satisfaction increases and average objective value decreases with $r$.
	The average objective is lower for CU, because connectedness magnifies the worst-case.
	Note that the objective value is only measured if constraints are satisfied.
	\item\textbf{\emph{CU vs. NC}}: For any $r$, CU solutions have higher constraint satisfaction than NC solutions, if correlation $\bs{\Psi}$ is positive, and higher objective function value, when the correlation is negative (see Figure~\ref{fig:ijoo_cs_ao}.  
	This is because connected sets account for dependency on the first period and provide additional protection beyond that of NC sets. 
	For any level of constraint satisfaction, the average objective of CU is better than NC (see Figure~\ref{fig:ts_05_aocs_sn}). 
	This can be observed either by changing the uncertainty set sizes or temporal correlations.
	\item\textbf{\emph{First vs Second period solutions}}:
	For a single estimate in (ii), Figure~\ref{fig:ts_05_aocs_sn} (right) shows that the optimal solution gradually concentrates only on \(\mb{x}_1\) for CU and only on \(\mb{x}_2\) for NC as $r$ increases.
	For NC, this is because $\mb{c}_2$ tends to be higher. 
	For CU, the second period weights are magnified due to connectedness, without corresponding benefits.
	Overall more components of \(\mb{x}_1\) and \(\mb{x}_2\) are non-zero for NC, resulting in higher objective but worse constraint satisfaction and vice versa.
	\item\textbf{\emph{Negative Correlation}}:
	When consecutive uncertainties are negatively correlated ($\bs{\Psi} < 0$), CU sets achieve lower constraint satisfaction but higher objective value than NC sets at any $r$ (see the electronic companion).
	This is because the worst-case uncertainty is dampened for CU sets by the negative correlation.
	For any level of constraint satisfaction, both models perform similarly.
	The solutions concentrate on \(\mb{x}_2\) for increasing $r$ because of higher $\mb{c}_2$ (see the figures in the electronic companion).
\end{itemize}

\subsection*{Summary}
CU sets exhibit clear advantages over nonconnected sets for autocorrelated uncertainties. 
Although CU has lower objective function values when autocorrelation is positive and lower constraint satisfaction when the auto correlation is negative, it outperforms NC for any given level of constraint satisfaction. 
The numerical analysis of modeling with \cus~highlights the improvements in constraint satisfaction over popular uncertainty sets.
Note that an alternative NC model with growing uncertainty sets over time, as displayed in Figure~\ref{fig:set_pic}(c), also does not improve over CU.
Therefore, we recommended to model with CU sets, particularly when uncertainties are positively correlated.

%%%%%%%%%%%%%%%%%%%%%%%%%%%%%%%%%%%%%%%%%%%%%%%%%%%%%%%%%%%%%%%%%%%%%%
\section{DRO Application: Portfolio Optimization}
\label{sec:dro_example}
The area of portfolio optimization has seen sizable advancements by systematically managing uncertainty, and DRO lends itself as a natural modeling framework.
In particular, different types of ambiguity sets have been proposed, such as moment-constrained sets~\citep{delage2010distributionally,natarajan2010tractable}, sets bounded by the Wasserstein metric~\citep{esfahani2018data}, etc.
In this section, we evaluate the performance of our \cu~framework on portfolio-optimization problems, resembling those by~\citet{delage2010distributionally}.
We investigate two instances of this problem.

\subsubsection*{Multi-asset Portfolio Problem}
Using historical data on five stocks, we probe our approach in a real-world setting. 
We model the asset returns with a time series fitted to the historic data for two periods of one week each.
We solve this problem for randomly selected days to evaluate average performance on actually realized returns.
We compare the CU model to two competing DRO models: one with same ambiguity sets for both periods, and one with the second-period ambiguity set centered at the expected second-period return (unconditioned).
We observe that the standard deviation of wealth for the CU model is lower than both DRO models, while achieving comparable returns.
This is because the CU model captures the compounding worst-case effects of autocorrelation, while managing conservatism.
We relegate the numerical details to the electronic companion.

\subsubsection*{Two-Asset Portfolio Problem} To gain further intuition in a controlled environment, we also simulate a two-asset portfolio optimization problem with synthetic data.
The goal of this experiment is to probe the long-term benefits of anticipating the behavior of uncertainties when making decisions, as opposed to short-term decision making based on a fixed uncertainty model.
Specifically, we consider a two-period problem, in which the choice of the portfolio has to be made at the beginning.
Such a problem can be expressed as
\begin{equation*}
\begin{aligned}
\max_{\mb{x}_1,\mb{x}_2} &\inf_{P_1 \in \widetilde{\mU}_1} \mathbb{E}_{P_1}\left[u_1(\mb{x}_1,\mb{d}_1) + \inf_{P_2 \in \widetilde{\mU}_2(d_1)}\mathbb{E}_{P_2}[u_2(\mb{x}_2,\mb{d}_2)]\right]\\
\text{s.t.}&\;\;\;\; \mb{e}^{\top}\mb{x}_t = 1 &&\forall t = 1,2\\
&\;\;\;\; \mb{x}_t \geq \mb{0} &&\forall t = 1,2,
\end{aligned}
\end{equation*}
where $u_t$ are utility functions, $\mb{x}_t$  decision variables, $\mb{d}_t$ return realizations, and $\widetilde{\mU}_t$ are the respective distributional uncertainty sets. 
The uncertainty sets $\widetilde{\mU}_1$ and $\widetilde{\mU}_2(\mb{d}_1)$ are as defined in~\eqref{eq:d}.
Each period corresponds to one week, and at the end of the first week, the assets can be reassigned.
However, the reallocation decision has to be specified at the beginning when solving the problem.

\begin{figure}[h!]
	\vspace{-5mm}		
	\centering
	\includegraphics[width=5.6cm,page=16]{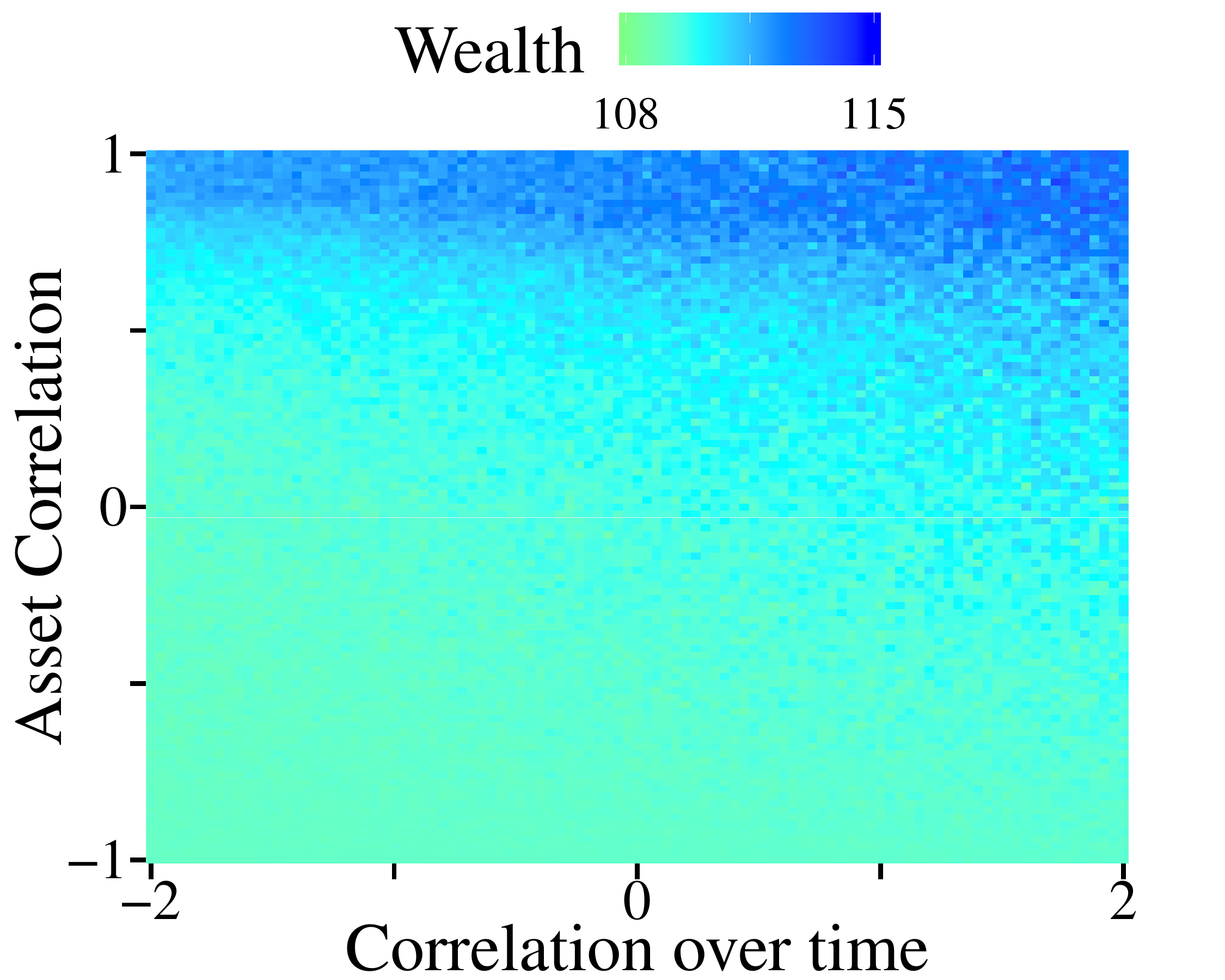}
	\includegraphics[width=5.6cm,page=18]{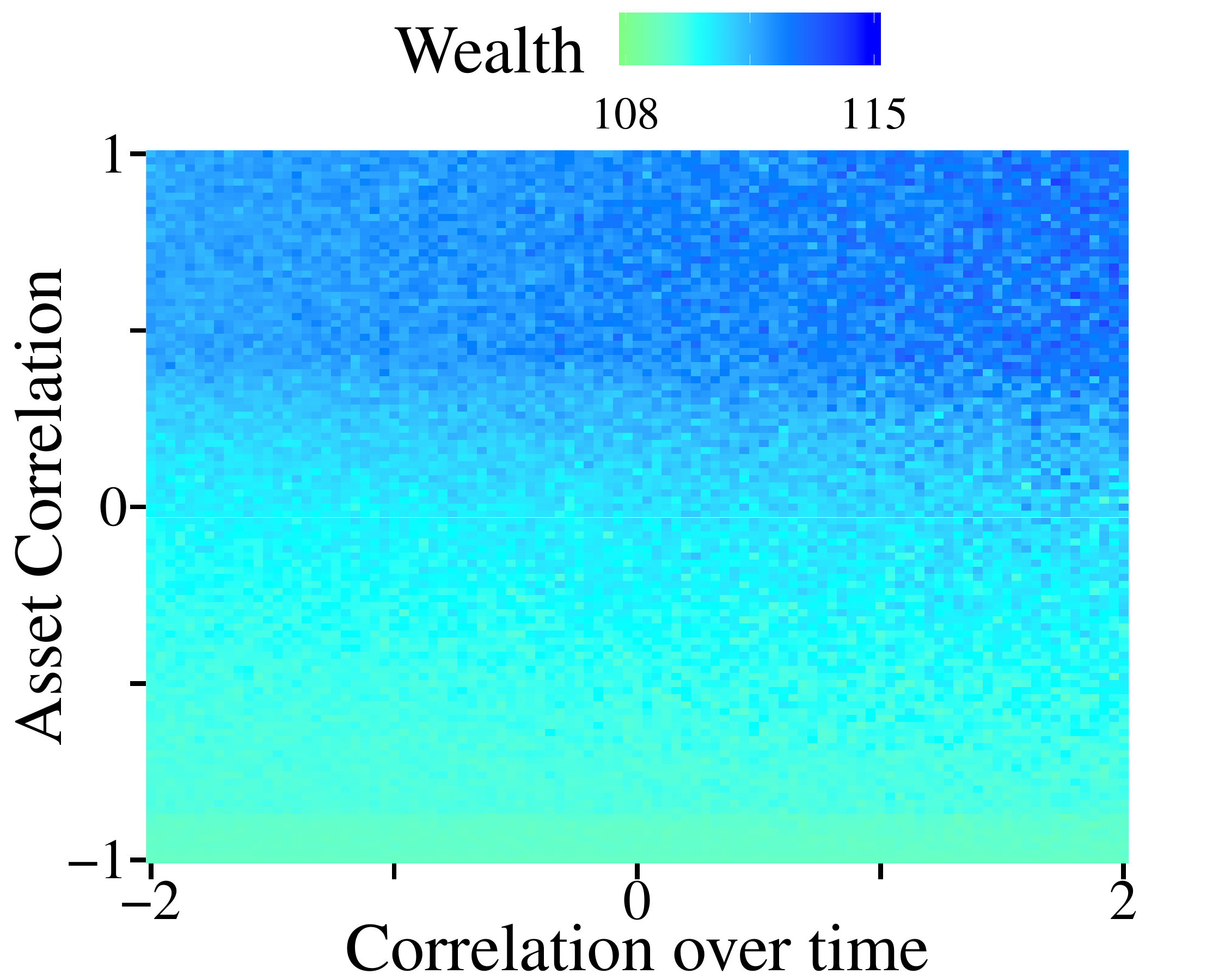}
	\caption{\label{fig:avg_wealth_3}Average Realized Wealth for CU and DRO at the end of Period 2.}
\end{figure}

To illustrate the impact of correlation among assets and correlation over time,
consider the above problem with a concave utility function for each period given by \(\min(1.5 r,0.015 + r, 0.06 + 0.2 r)\). 
Here, \(r = x_1d_{11} + x_2d_{12}\) represents the portfolio return in each period.
For this problem, the expected return in the first period resides within a box centered at \(\bs{\mu}_1 = [0.03,0.06]\) with size $\bs{\delta} = [0.02,0.02]$. 
The expected return in the second period lies inside a similar set, whose center depends on the realized return of the first period as \(\bs{\mu}_1 + \omega \cdot (\mb{d}_1 - \bs{\mu}_1)\), where \(\omega\) is a parameter controlling the correlation over time. 
Furthermore, the covariance of the return in both periods is bounded by \(\bs{\Sigma}_1\) 
\begin{align*}
\bs{\Sigma}_1= 
\left[\begin{array}{l l l}
0.005 & & 0.005\cdot\rho \\ 
0.005\cdot \rho & & 0.005
\end{array}\right], 
\end{align*}
where \(\rho\) measures the correlation among the asset returns. 
We vary the parameters \(\omega\) over the range \([-2,2]\) and \(\rho \in[-1,1]\) to study the impact of correlations on the following metrics of the portfolio model:
\begin{enumerate}
	\item \emph{Asset allocation:} to observe the behavior of decisions and to develop intuition (Figures~\ref{fig:asset2_cu_allocation}). 
	\item \emph{Average realized wealth:} to test the performance under realistic (random) setting (Figure~\ref{fig:avg_wealth_3}).	
	\item \emph{Difference in standard deviation of the CU and DRO models:} to probe sensitivity (Figure~\ref{fig:diff_cu_dro_3}).
\end{enumerate}

We now discuss the outcomes of the numerical experiments with regard to these metrics.
Additional metrics are discussed in the electronic companion.
In each of the upcoming figures, the horizontal axis represents the correlation over time \(\omega\), and  the vertical axis the correlation among the assets.
We will divide the discussion into these two aspects. 
The color in the plots represents the metric under consideration and its values are displayed in the legend.

\subsection{Asset Allocation}
\begin{figure}[h!]
	\vspace{-2mm}
	\centering		
	\includegraphics[width=5.6cm,page=2]{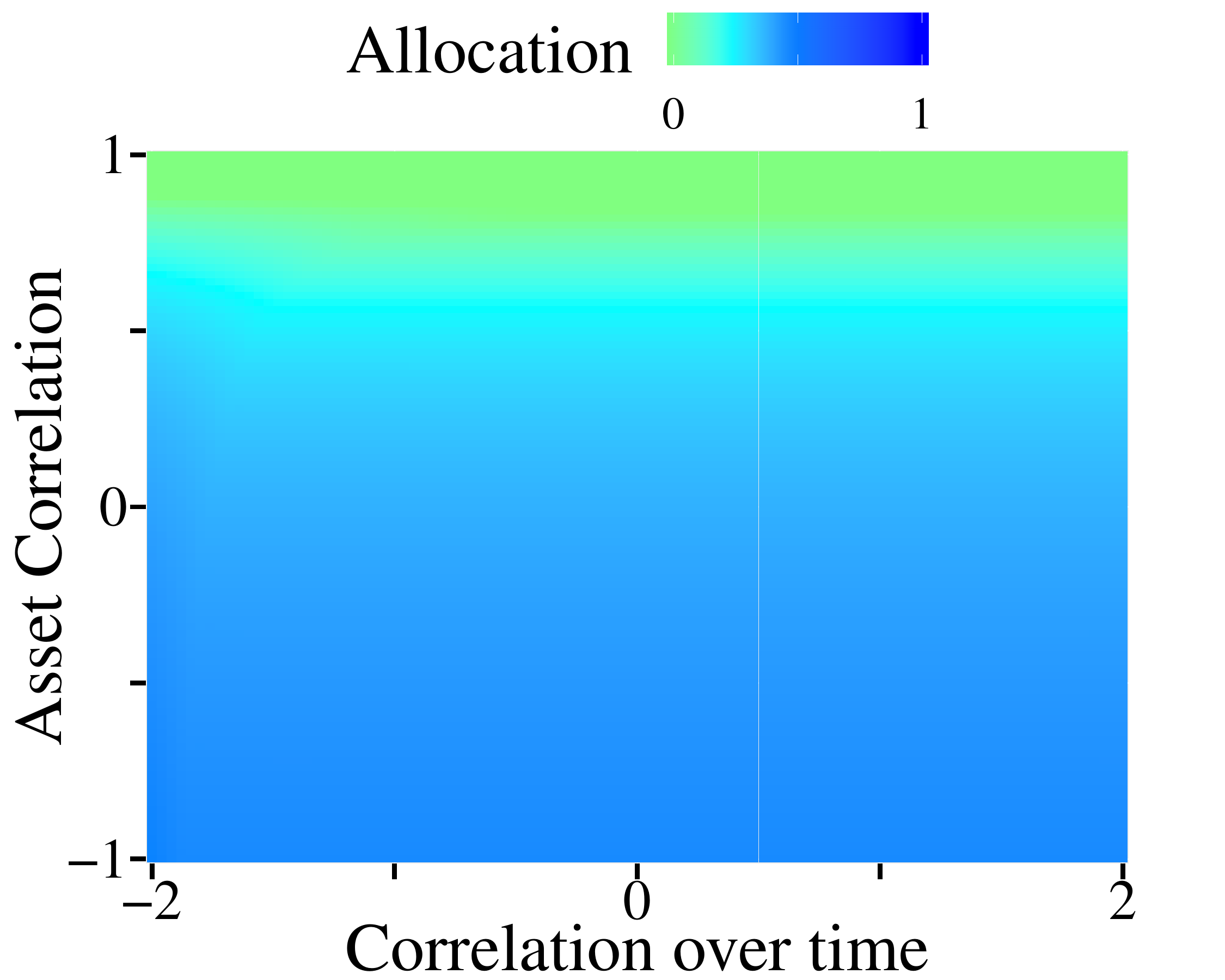}
	\hspace{-5mm}
	\includegraphics[width=5.6cm,page=4]{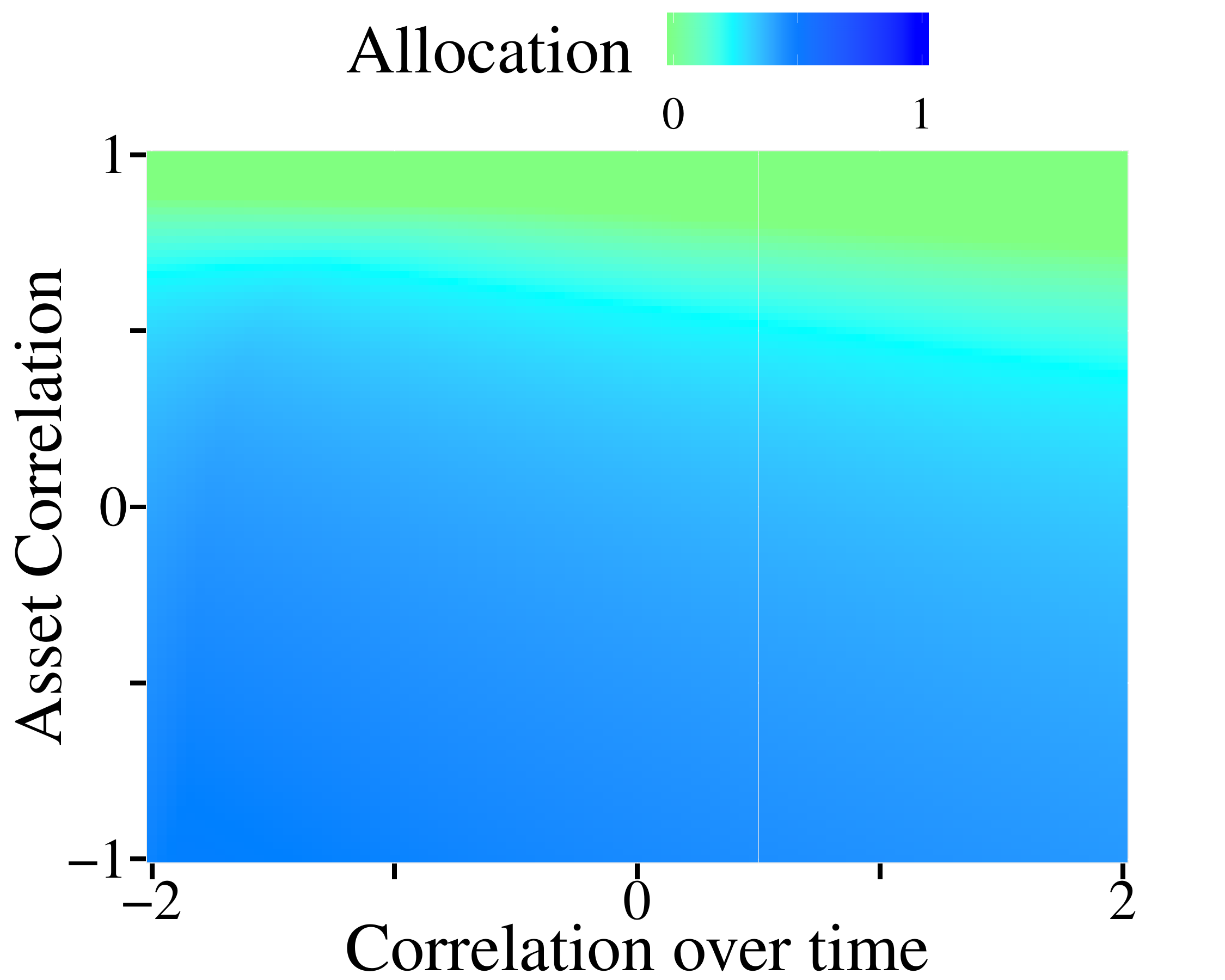}
	\hspace{-5mm}
	\includegraphics[width=5.6cm,page=6]{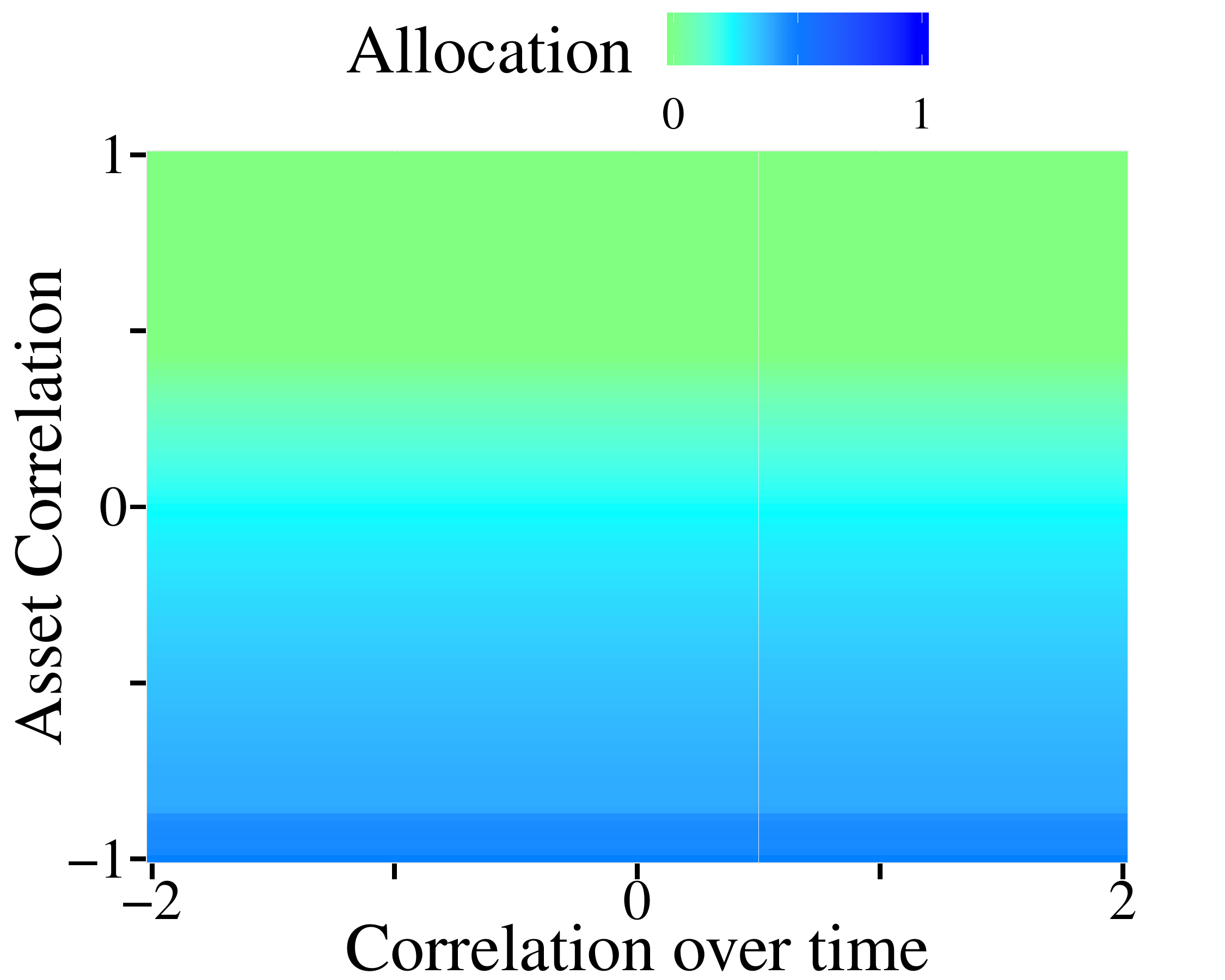}
	\caption{\label{fig:asset2_cu_allocation}CU Allocation to Asset 1 for CU Period 1 (Left), CU Period 2 (Center), and DRO Period 1 (Right).}
\end{figure}

Figure~\ref{fig:asset2_cu_allocation} displays the asset allocation for asset 1 for the CU model for period 1 (left), period 2 (center), and the DRO model in period 1 (right), which is the same as period 2. 
The sum of allocations is one. 
Figure~\ref{fig:asset2_cu_allocation} displays that both the CU and DRO models allocate more wealth to asset 2.
This is due to the fact that asset 2 has higher expected returns.

\subsubsection*{Asset Correlation}
For both CU and DRO models, the allocation to asset 2 is reduced if the assets are positively correlated. 
This is because a positive correlation reduces the protection against the worst case that would arise by spreading the wealth among the two assets. 
A negative correlation among the assets increases this protection, resulting in a broader spread. 

\subsubsection*{Correlation over time}
We also observe that as the correlation over time increases, the CU model concentrates on asset 2. 
This is because it has higher expected return and the benefit of hedging against bad returns by spreading the wealth among both assets is reduced by the positive correlation. 
The correlation over time mostly affects the assets in the second period. 
A higher correlation leads to more allocations to asset 2 since it has a higher expected value and, thus, leads to better worst case performance than asset 1. 	

\subsection{Average wealth}
\subsubsection*{Asset Correlation}
Figure~\ref{fig:avg_wealth_3} shows the average wealth at the end of the second period. 
It can be observed that high correlation between the assets leads to higher wealth. 
This contrasts with the worst-case performance and is because positive correlations also promote high realized returns. 
The wealth increases faster for the DRO model due to lower conservatism of the model. 
\subsubsection*{Correlation over time}
The realized wealth increases with correlation over time as it leads to higher realized returns. 
This behavior is more visible for the CU as compared to the DRO model. 

\subsection{Difference in CU and DRO models}
\begin{figure}[h!]
	\vspace{-2mm}	
	\centering
	\includegraphics[width=5.6cm,page=36]{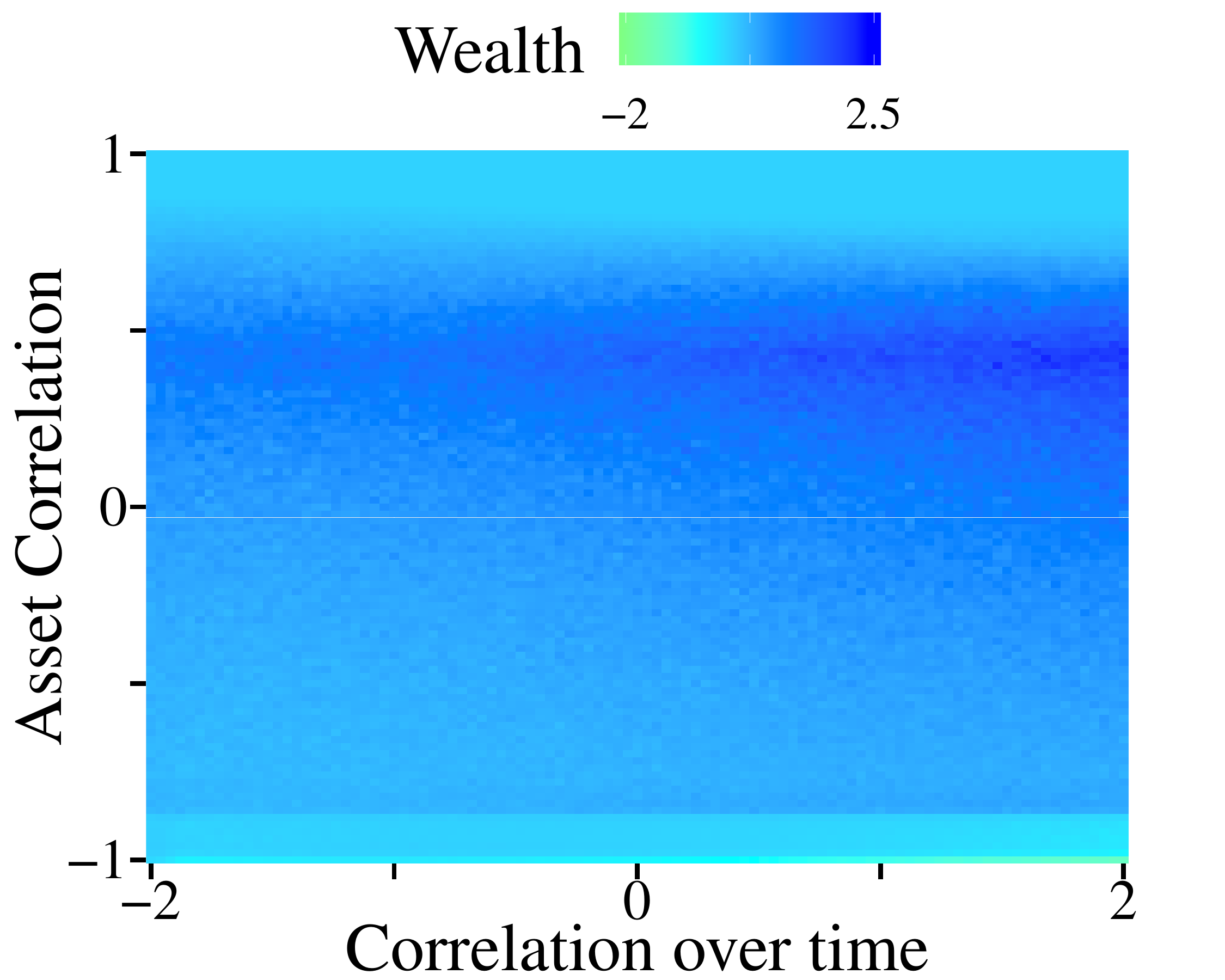}
	\includegraphics[width=5.6cm,page=38]{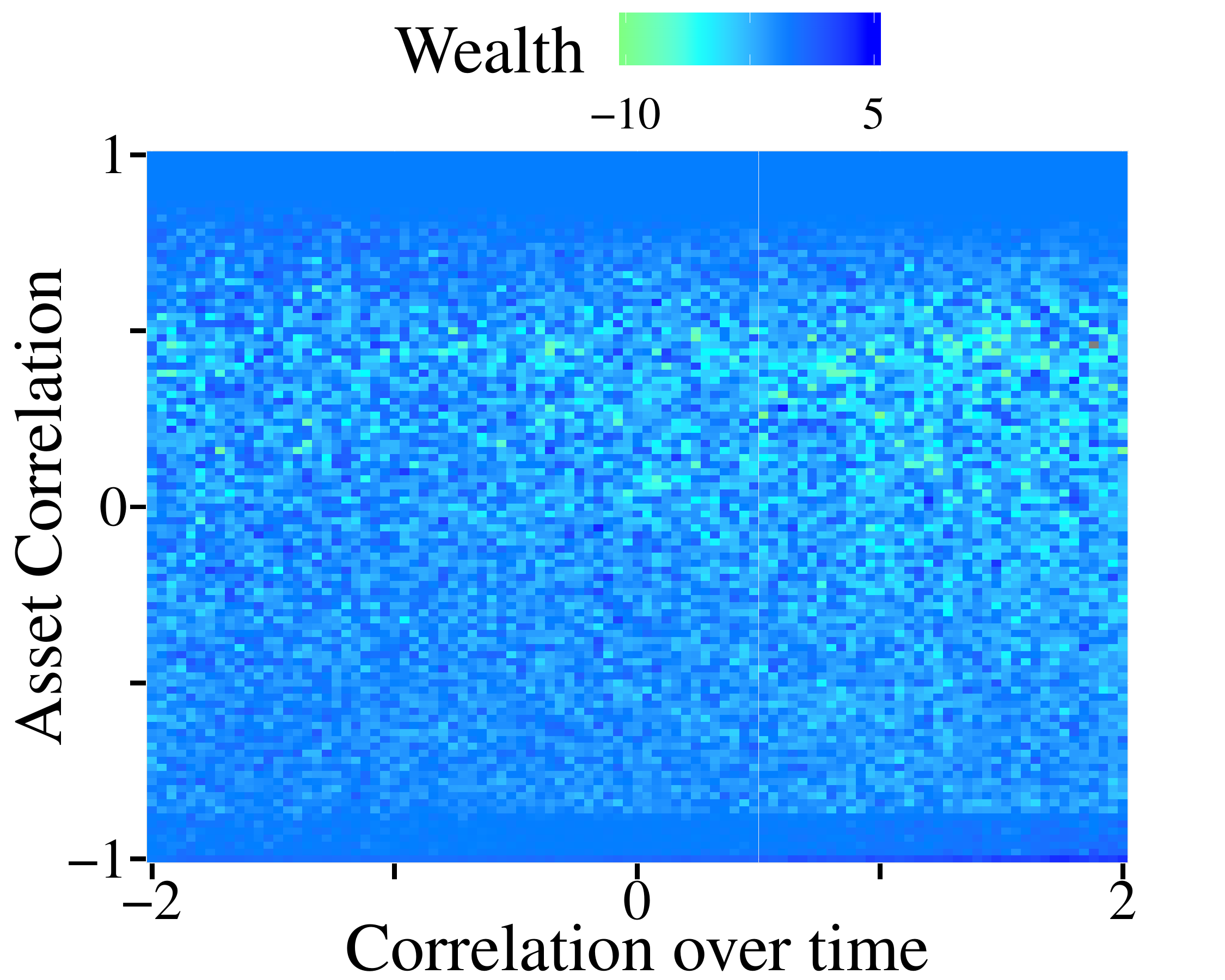}	\caption{\label{fig:diff_cu_dro_3}Difference in Wealth Standard Deviation (Left) and Difference in Worst-Case Wealth (Right) for Period 2.}
\end{figure}

Figure~\ref{fig:diff_cu_dro_3} (left) shows the difference in wealth standard deviation between the CU and DRO models, and Figure~\ref{fig:diff_cu_dro_3} (right) the difference in worst-case wealth between these two models.
\subsubsection*{Asset Correlation}
The difference in standard deviations is always positive, increases with the asset correlation, and achieves a peak at 0.4. 
That means that the CU model always leads to lower standard deviation, highlighting robustness of solutions.
Also, the worst-case wealth of the CU model tends to perform better as the correlation increases, and achieving best performance at 0.4.

\subsubsection*{Correlation over time}
The benefits of the CU model in terms of lower standard deviation and better worst-case wealth increases with correlation over time. 

%%%%%%%%%%%%%%%%%%%%%%%%%%%%%%%%%%%%%%%%%%%%%%%%%%%%%%%%%%%%%%%%%%%%%%%%%%%%%%%%%%%%%%%%%%%%%%%%%%%%%%%%

\section{Conclusion}
\label{sec:5_CO}
In multistage optimization under uncertainty, the model of the underlying uncertainty plays a decisive role for both the solution quality and the computational effort.
Current techniques usually assume a collection of predetermined sets, which are independent of each other.
Given that in many applications, the uncertainty is connected over multiple periods, these methods are either over conservative, suffering from sizable suboptimality, or do not adequately capture all realizations of the uncertainty, risking the feasibility of solutions.
We introduce a new modeling framework in which the uncertainty model depends on past realizations.
This work extends the efficacy of robust optimization, as well as distributionally robust optimization, to the connected uncertainty paradigm.
We study commonly used constraints and uncertainty set dependencies and provide their respective reformulations.
Specifically for robust settings, tractable reformulations are developed for polyhedral and ellipsoidal uncertainty sets with linear and quadratic dependence, as they occur in a wide range of decision-making processes.
When the uncertainty is modeled in a distributionally robust fashion, 
we provide reformulations for moment-constrained ambiguity sets with the mean depending on the past, as in time-series settings.
Numerical experiments on a knapsack problem exhibit sizable improvements in constraint satisfaction and better objective performance for any fixed level of constraint satisfaction, when uncertainties are modeled with the proposed connected sets.
Similarly, a distributionally robust portfolio-optimization problem achieves approximately the same expected returns at narrower standard deviations, highlighting the robustness of solutions, when modeling with \cus.
Thus, we recommend using connected uncertainty sets for autocorrelated or time-varying uncertainties because they improve the performance over nonconnected sets. 
%.
Since in multi-period problems, uncertainties are naturally connected across periods,  
the proposed approach offers a general modeling framework that can be applied to numerous operational applications.

\bibliographystyle{informs2014} 
\bibliography{DyUS_short} 

%%%%%%%%%%%%%%%%%
\end{document}

\RUNAUTHOR{Nohadani and Sharma}
\RUNTITLE{Optimization under Connected Uncertainty}

	%%%%%%%%%%%%%%%%
	\ECHead{Electronic Companion to
	``Optimization under Connected Uncertainty''}
	
	\section{Proof of Theorem~3}
	\label{proof_th_prop_poly_rhs}
	The following is the proof for the robust counterpart for the polyhedral uncertainty set.
	\proof{Proof:}
	For (P), we define the joint uncertainty set \(\mU\) as
	\begin{eqnarray*}
		\mU &=& \left\{(\mb{d}_1^{\top},\mb{d}_2^{\top},\dots,\mb{d}_T^{\top})^{\top} \mid \mb{d}_1 \in \mU_1,\;\mb{d}_t \in \mU_t(\mb{d}_{t-1}) \;\forall t=2,\dots,T \right\}\\
		&=& \left\{(\mb{d}_1^{\top},\mb{d}_2^{\top},\dots,\mb{d}_T^{\top})^{\top} \mid \mb{G}_1 \mb{d}_1 \geq \mb{g}_1,\; \mb{G}_t \mb{d}_t \geq \mb{g}_t + \bs{\Delta}_t \mb{d}_{t-1}\; \forall t=2,\dots,T\right\},
	\end{eqnarray*}
	combining the polyhedral sets for each period.
	The robust counterpart of~(C-RO) becomes
	\[\sum_{t=1}^{T}\mb{d}_t^{\top}\mb{x}_t \leq B \;\;\forall (\mb{d}_1^{\top},\mb{d}_2^{\top},\dots,\mb{d}_T^{\top})^{\top} \in \mU, \]
	which can be rewritten as 
	\[\max_{(\mb{d}_1^{\top},\mb{d}_2^{\top},\dots,\mb{d}_T^{\top})^{\top} \in \mU} \sum_{t=1}^{T}\mb{d}_t^{\top}\mb{x}_t \leq B.\]
	The LHS is computed by
	\begin{align*}
		\max_{\mb{d}_t}&\;\; \sum_{t=1}^{T}\mb{d}_t^{\top}\mb{x}_t\\
		\text{s.t. } &\;\; \mb{G}_1 \mb{d}_1 \geq \mb{g}_1 \\
		&\;\; \mb{G}_t \mb{d}_t \geq \mb{g}_t + \bs{\Delta}_t \mb{d}_{t-1} \quad\;\;\forall t=2,\dots,T.
	\end{align*}
	Using duality, the robust counterpart of~(C-RO) is  
	\begin{align*}
		\sum_{t=1}^{T} \mb{q}_t^{\top}\mb{b}_t &\leq B\\
		\mb{q}_t^{\top}\mb{A}_t  - \mb{q}_{t+1}^{\top}\bs{\Delta}_{t+1} &= \mb{x}_{t}^{\top}, \quad\forall t=1,\dots,T\\
		\mb{q}_t &\leq 0 \quad\;\; \forall t=1,\dots,T,
	\end{align*}
	which is the desired result.  \Halmos
	\endproof

	%%%%%%%%%%%%%%%%%%%%%%%%%%%%%%%%%%%%%%%%%%%%
	
	\section{ARO with Ellipsoidal CU Sets for Connected Uncertainty with RO}
	\label{sec:aro_ellipsoidal}
	In this section, we reformulate the connected constraints from problem~(CU-ARO) for ellipsoidal uncertainty sets, where the center depends on the previous period realization as
	\begin{equation}
		\mathcal{U}_t(\mb{d}_{t-1}) = \{\mb{d}_t \mid \mb{d}_t = \bs{\mu}_t(\mb{d}_{t-1}) + \mb{L}_t \mb{u}_t : \|\mb{u}_t\|_2 \leq r_t\},
	\end{equation}
	where \(\mb{L}_t \mb{L}_t^{\top} = \bs{\Sigma}_t\).
	The dependence of $\bs{\mu}_t$  on the previous period realization is given by
	\begin{equation}
		\bs{\mu}_t(\mb{d}_{t-1}) = \mb{A}_t \bs{\mu}_{t-1}(\mb{d}_{t-2}) + \mb{F}_t \mb{d}_{t-1} + \mb{c}_t.
	\end{equation}
	For this setting, the following theorem provides the reformulation for \cus~with affine decision policies.
	
	\begin{theorem}
		The two-period adjustable optimization problem~(CU-ARO) has a tractable reformulation, when the uncertainty affects the center and the fully adaptive decisions are replaced by affine decision rules. 
	\end{theorem}
	\proof{Proof:}
	We replace \(\mb{x}_2(\mb{d}_1)\) with the affine decision rule \(\mb{x}_2(\mb{d}_1) = \mb{X}_2\mb{d}_1\) and expand \(\mb{b}_1(\mb{d}_1) = \mb{B}_1\mb{d}_1\) and \(\mb{b}_2(\mb{d}_2) = \mb{B}_2\mb{d}_2\).
	We focus on reformulating the second constraint, which is affected by the connected uncertainty  and whose robust problem is
	\begin{align*}
		\max_{\mb{d}_1,\mb{d}_2} &\; \mb{B}_{2,i}^{\top}\mb{d}_2 - [\mb{A}_{22}\mb{X}_{2}]_{i}^{\top}\mb{d}_1\\
		\text{s.t.} \;\;&\; \mb{d}_1 = \bs{\mu}_1 + \mb{L}_1 \mb{u}_1\\
		&\; \mb{d}_2 = \bs{\mu}_2(\mb{d}_1) + \mb{L}_2 \mb{u}_2\\
		&\; \|\mb{u}_t\|_2 \leq r_t \quad \forall t=1,2.
	\end{align*}
	Substituting \(\mb{d}_1\) and \(\mb{d}_2\), we rewrite the problem as 
	\begin{align*}
		\max_{\mb{d}_1,\mb{d}_2} &\;\; \mb{B}_{2,i}^{\top}(\mb{A}_2 \bs{\mu}_1 + \mb{F}_2 \mb{d}_1 + \mb{c}_2)+ \mb{B}_{2,i}^{\top}\mb{L}_2\mb{u}_2 - [\mb{A}_{22}\mb{X}_{2}]_{i}^{\top}\bs{\mu}_1 - [\mb{A}_{22}\mb{X}_{2}]_{i}^{\top}\mb{L}_1 \mb{u}_1\\
		\text{s.t.} \;\;& \|\mb{u}_t\|_2 \leq r_t \quad\forall t=1,2.
	\end{align*}
	Using the dual, the complete second constraint of~(CU-ARO) is given by
	\begin{align*}
		\mb{A}_{21,i}^{\top}\mb{x}_1 &\geq \mb{B}_{2,i}^{\top}(\mb{A}_2 \bs{\mu}_1 + \mb{F}_2\bs{\mu}_1 +  \mb{c}_2) - [\mb{A}_{22}\mb{X}_{2}]_{i}^{\top}\bs{\mu}_1 + r_2 \|\mb{B}_{2,i}^{\top}\mb{L}_2\|_2 + r_1\|\mb{B}_{2,i}^{\top}\mb{F}_2 \mb{L}_1 \\ 
		&- [\mb{A}_{22}\mb{X}_{2}]_{i}^{\top}\mb{L}_1\|_2.
	\end{align*}
	The remaining constraints in~(CU-ARO) is then reformulated in a similar manner, leading to a tractable reformulation.  \Halmos
	\endproof	
	
	%%%%%%%%%%%%%%%%%%%%%%%%%%%%%%%%%%%%%%%%%%%%%%%%%%%%%%%%
	\section{Extensions to Affinely adaptive Decisions for Connected Uncertainty with RO}
	In previous sections, we focused on the effect of connected uncertainties on \emph{static} and here-and-now decisions.
	In many problems, decisions can also adapt to previously revealed uncertainties.
	To this end, ARO has applications in inventory management~\citep{solyali2015impact} and unit commitment~\citep{lorca2017multistage}, amongst others. 
	ARO problems, however, are known to be NP-complete, which can be circumvented by
	leveraging affine or piece-wise static decision rules~\citep{hanasusanto2015k}. 
	
	In this section, we extend the notion of modeling with \cus~to \emph{adaptive} decisions and provide reformulations for a two-period ARO problem under affine decision rules. 
	Such a problem can be expressed as
	\begin{equation}
		\tag{CU-ARO}
		\label{eq:cu_aro}
		\hspace*{-5mm}
		\begin{aligned}
			\min_{\substack{\mb{x}_1\\\mb{x}_2(\mb{d}_1)}} \!
			\max_{\substack{\mb{d}_1 \in \mU_1\\\mb{d}_2 \in \mU_2(\mb{d}_1)}} &\; \mb{c}_1^{\top}\mb{x}_1 +  \mb{c}_2^{\top}\mb{x}_2(\mb{d}_1)\\
			\text{s.t.} \;\;\;&\; \mb{A}_{11} \mb{x}_1 \geq \mb{b}_1(\mb{d}_1) &&\forall \mb{d}_1 \in \mU_1 \\
			&\; \mb{A}_{21}\mb{x}_1 + \mb{A}_{22} \mb{x}_2(\mb{d}_1) \geq \mb{b}_2(\mb{d}_2) &&\forall \mb{d}_2 \in \mU_2(\mb{d}_1) \, \forall \mb{d}_1 \in \mU_1\\
			&\; \mb{x}_2(\mb{d}_1) \geq \mb{0} &&\forall \mb{d}_1 \in \mU_1\\
			&\; \mb{x}_1 \geq \mb{0}.
		\end{aligned}
	\end{equation}
	We focus on polyhedral \cus, where the RHS depends on the previous period realization as
	\begin{equation*}
		\mU_1 = \{\mb{d}_1 \;|\; \mb{G}_1 \mb{d}_1 \geq \mb{g}_1\}, \quad\quad
		\mU_2(\mb{d}_1) = \{ \mb{d}_2 \;|\; \mb{G}_2 \mb{d}_2 \geq \mb{g}_2 + \bs{\Delta}_1 \mb{d}_1\}.
	\end{equation*}
	For this setting, the following theorem shows the impact of \cus~with affine decision policies on the fully adaptive decisions \(\mb{x}_2(\mb{d}_1) = \mb{x}_2^0 + \mb{X}_2\mb{d}_1\).
	For brevity, we incorporated the constant term \(\mb{x}_2^0\) into the matrix \(\mb{X}_2\) by modifying the uncertainty set and setting the first component of \(\mb{d}_1\) to be equal to \(1\).
	Then we can write  
	\(\mb{x}_2(\mb{d}_1) = \mb{X}_2\mb{d}_1\).

	\begin{theorem}
		\label{thm:cu_adaptive}
		The two-period ARO problem~(CU-ARO) has a tractable reformulation, when the uncertainty affects the RHS linearly and the variables are replaced by affine decision rules. 	
	\end{theorem}
	\proof{Proof:}
	We replace \(\mb{x}_2(\mb{d}_1)\) with the affine decision rule \(\mb{x}_2(\mb{d}_1) = \mb{X}_2 \mb{d}_1\) and expand \(\mb{b}_1(\mb{d}_1) = \mb{B}_1\mb{d}_1\) and \(\mb{b}_2(\mb{d}_2) = \mb{B}_2\mb{d}_2\).
	We focus on reformulating the second constraint which is affected by the connected uncertainty  and whose robust problem is
	\begin{equation*}
		\begin{aligned}
			\max_{\mb{d}_1,\mb{d}_2} &\; \mb{B}_{2,i}^{\top}\mb{d}_2 - [\mb{A}_{22}\mb{X}_{2}]_{i}^{\top}\mb{d}_1\\
			\text{s.t.} \;\;&\; \mb{G}_1 \mb{d}_1 \geq \mb{g}_1\\
			&\; \mb{G}_2 \mb{d}_2 \geq \mb{g}_2 + \Delta_1 \mb{d}_1.
		\end{aligned}
		%\vspace{-2mm}
	\end{equation*}
	Using the dual, the complete second constraint of~(CU-ARO) is given by
	\vspace{-1mm}
	\begin{equation*}
		\begin{aligned}
			\mb{A}_{21}\mb{x}_1 &\geq \mb{P}^{\top}\mb{g}_1 + \mb{Q}^{\top}\mb{g}_2\\
			\mb{P}^{\top}\mb{G}_1 + \mb{Q}^{\top}\bs{\Delta}_1  &= \mb{B}_{2}\\
			\mb{Q}^{\top}\mb{G}_2 &= - \mb{A}_{22}\mb{X}_{2}\\
			\mb{P},\mb{Q} &\leq \mb{0}.
		\end{aligned}
		\vspace{-1mm}
	\end{equation*}
	The columns of \(\mb{P}\) and \(\mb{Q}\) correspond to dual variables of the original problem. 
	The remaining constraints in~(CU-ARO) can be reformulated similarly, leading to a tractable reformulation. \Halmos
	\endproof
	
	Theorem~2 extends the result of modeling with standard sets for an adaptive setting, where such reformulations are possible for RHS uncertainty, to \cus.
	When uncertainties affect the LHS, even problems with standard and non-connected sets are difficult to reformulate~\citep{ben2009robust}.
	Therefore, we forego their extension to the CU setting.

	\section{Proof of Theorem~2}
	\proof{Proof:}
	The reformulation proceeds through induction. 
	The robust counterpart of~(C-RO) with respect to \(\mb{d}_{k+1},\dots,\mb{d}_{T}\) is
	\begin{align*}
		\label{eq:prop_MD_lem}
		&s_{k} + \Theta_{k+1} + \sum_{t = k+1}^{T} r_{t}\sqrt{A_{k+1,t} \mb{y}_{t}^{\top}\bs{\Sigma}_{k+1}\mb{y}_{t}} + R_{k+1} \leq B \quad \forall \mb{n}_{[k+1]} \in \mathcal{N}_{[k+1]},\nonumber\\
		\text{with } & s_{k} = \sum\limits_{t=1}^{k}\mb{d}_t^{\top}\mb{x}_t ,\;\;
		\Theta_{k+1} = \!\!\!\!\!\!\sum\limits_{t={k+1}}^{T} \!\! \bs{\mu}_{t}^{\top}\mb{x}_{t} \text { and } s_0 = 0.
	\end{align*}
	The proof proceeds parallel to that of Theorem~1.
	\paragraph{Base Case:}
	The Constraint~(C-RO) can be written as 
	$s_{T-1} + \mb{d}_T^{\top}\mb{x}_T \leq B$.
	Substituting \({\mb{d}_T = \bs{\mu}_T + \mb{L}_T\mb{u}_T}\).
	(C-RO) must hold for all \(\mb{u}_T\) with \(\|\mb{u}_T\|_2 \leq r_T\), and thus it also holds for 
	$s_{T-1} + \bs{\mu}_T^{\top}\mb{x}_T + r_T\|\mb{L}_T^{\top}\mb{x}_T\|_2 \leq B$.
	Using the equality \(\|\mb{L}_T^{\top}\mb{x}_T\|_2 =  
	\sqrt{\mb{x}_T^{\top}\bs{\Sigma}_T\mb{x}_T} \) and 
	\(\mb{y}_T = \mb{x}_T\), we obtain
	\[s_{T-1} + \Theta_T + r_T\|\mb{L}_T^{\top}\mb{y}_T\|_2 \leq B \Leftrightarrow s_{T-1} + \Theta_T + r_T\sqrt{\mb{y}_T^{\top}\bs{\Sigma}_T\mb{y}_T} \leq B.\]
	Since \(R_{T+1}\) is assumed to be zero, we have achieved the desired result. 
	\paragraph{Inductive case:}
	Assume that the result is true for \(k+1\)--that is, the reformulation with respect to \(\mb{d}_T,\mb{d}_{T-1},\dots \mb{d}_{k+1}\)
	is given by
	\begin{align*}
		s_{k} &+ \Theta_{k+1} + \sum_{t=k+1}^{T} r_{t}\sqrt{A_{k+1,t} \mb{y}_{t}^{\top}\bs{\Sigma}_{k+1}\mb{y}_{t}} + R_{k+1} \leq B \quad \forall \mb{n}_{[k+1]} \in \mathcal{N}_{[k+1]}.
	\end{align*}
	Now, we have to prove that the reformulation with respect to \(\mb{d}_{k}\) is given by 
	\begin{align*}
		s_{k-1} &+ \Theta_{k} + \sum_{t=k}^{T} r_{t}\sqrt{A_{k,t} \mb{y}_{t}^{\top} \bs{\Sigma}_{k}\mb{y}_{t}} + R_{k} \leq B \quad \forall \mb{n}_{[k]} \in \mathcal{N}_{[k]}.
	\end{align*}
	Substituting \(s_{k} = s_{k-1} + \mb{d}_{k}^{\top}\mb{x}_{k}\) and \(\mb{d}_{k} = \bs{\mu} _{k} + \mb{L}_{k}\mb{u}_{k}\) in the inductive assumption
	\begin{align*}
		&s_{k-1} + \Theta_{k+1} + \bs{\mu}_{k}^{\top}\mb{x}_{k} + \mb{u}_{k}^{\top}\mb{L}_{k}^{\top}\mb{x}_{k}\\
		&\hspace{0.61cm}+\sum_{t=k+1}^{T} r_{t}\sqrt{A_{k+1,t} \mb{y}_{t}^{\top}\bs{\Sigma}_{k+1}\mb{y}_{t}} + R_{k+1} \leq B \;\;\forall \mb{n}_{[k+1]} \in \mathcal{N}_{[k+1]}.
	\end{align*}
	Utilizing~(5) for \(t=k\) and \(\Theta_{k} = \bs{\mu}_{k}^{\top}\mb{x}_{k} + \Theta_{k+1}\), the constraint can be reformulated as 
	\begin{align*}
		& s_{k-1} + \Theta_{k}  + \mb{u}_{k}^{\top}\mb{L}_{k}^{\top}\mb{x}_{k} \\ 
		&+\sum_{t=k+1}^{T} r_{t}\bigg(A_{k+1,t}a_k \mb{y}_{t}^{\top}\bs{\Sigma}_{k}\mb{y}_{t} + A_{k+1,t}f_k{({(\mb{d}_{k} - \bs{\mu}_{k})}^{\top}\mb{y}_{t})}^2 \\
		&+ A_{k+1,t}\mb{y}_{t}^{\top}\mb{C}_k \mb{y}_{t}\bigg)^{\frac{1}{2}} 
		+  R_{k+1} \leq B \;\;\forall \mb{n}_{[k+1]} \in \mathcal{N}_{[k+1]}.
	\end{align*}
	Using the inequality,  \(\sqrt{a^2+b^2} \leq |a| + |b| \), 
	a more conservative constraint is
	\begin{align*}
		&s_{k-1} + \Theta_{k}  + \mb{u}_{k}^{\top}\mb{L}_{k}^{\top}\mb{x}_{k} 
		+ \sum_{t=k+1}^{T} r_{t}\sqrt{A_{k+1,t}a_k \mb{y}_{t}^{\top}\bs{\Sigma}_{k}\mb{y}_{t}} \\
		&+ \sum_{t=k+1}^{T}r_{t}\sqrt{A_{k+1,t}f_k}\left|{(\mb{d}_{k} - \bs{\mu}_{k})}^{\top}\mb{y}_{t}\right| \\
		&+ \sum_{t=k+1}^{T}r_{t}\sqrt{A_{k+1,t}\mb{y}_{t}^{\top}\mb{C}_k \mb{y}_{t}}+ R_{k+1} \leq B \;\;\forall \mb{n}_{[k+1]} \in \mathcal{N}_{[k+1]}.
	\end{align*}
	The above constraint holds for \(\left|{(\mb{d}_{k} - \bs{\mu}_{k})}^{\top}\mb{y}_{t}\right|\), if it holds for both the positive and the negative possible values.
	We use \(n_{k,t}\), which can take values \(1\) or \(-1\) to express the above as
	\begin{align*}
		&s_{k-1} + \Theta_{k}  + \mb{u}_{k}^{\top}\mb{L}_{k}^{\top}\mb{x}_{k} + \sum_{t=k+1}^{T}r_{t}\sqrt{A_{k+1,t}a_k \mb{y}_{t}^{\top}\bs{\Sigma}_{k}\mb{y}_{t}} \\
		& + \sum_{t=k+1}^{T}n_{k,t} \cdot r_{t}\sqrt{A_{k+1,t}f_k}{(\mb{d}_{k} - \bs{\mu}_{k})}^{\top}\mb{y}_{t} \\
		&+ \sum_{t=k+1}^{T}r_{t}\sqrt{A_{k+1,t}\mb{y}_{t}^{\top}\mb{C}_k \mb{y}_{t}} + R_{k+1} \leq B \\
		&\hspace{4cm}\forall n_{k,t} \in \{1,-1\},\;\forall t = k+1,\dots,T,\; \forall \mb{n}_{[k+1]} \in \mathcal{N}_{[k+1]}.
	\end{align*}
	Note that the use of $n_{k,t}$ allows to write the $2^{T-k}$ constraints in this concise form.
	Combining the set \(\mathcal{N}_{[k+1]}\) with \(\{1,-1\}^{T-k}\),  
	using \(\mb{d}_{k} - \bs{\mu}_{k} = \mb{L}_{k}\mb{u}_{k}\), and collecting the terms involving \(\mb{u}_{k}^{\top}\mb{L}_{k}^{\top}\), we obtain
	\begin{align*}
		&s_{k-1} + \Theta_{k}  + \mb{u}_{k}^{\top}\mb{L}_{k}^{\top}(\mb{x}_{k} + \sum_{t=k+1}^{T}n_{k,t}\cdot r_{t}\sqrt{A_{k+1,t}f_k}\mb{y}_{t})\\ 
		&+ \sum_{t=k+1}^{T} r_{t}\sqrt{A_{k+1,t}a_k \mb{y}_{t}^{\top}\bs{\Sigma}_{k}\mb{y}_{t}} + \sum_{t=k+1}^{T}r_{t}\sqrt{A_{k+1,t}\mb{y}_{t}^{\top}\mb{C}_k \mb{y}_{t}}\\ &+ R_{k+1} \leq B \;\;\forall \mb{n}_{[k]} \in \mathcal{N}_{[k]}.
	\end{align*}
	With the definition of \(\mb{y}_{k}\), we can write 
	\begin{eqnarray*}
		s_{k-1} + \Theta_{k}  &+& \mb{u}_{k}^{\top}\mb{L}_{k}^{\top}\mb{y}_{k} + \sum_{t=k+1}^{T} r_{t}\sqrt{A_{k+1,t}a_k \mb{y}_{t}^{\top}\bs{\Sigma}_{k}\mb{y}_{t}}\\
		&+& \sum_{t=k+1}^{T}r_{t}\sqrt{A_{k+1,t}\mb{y}_{t}^{\top}\mb{C}_k \mb{y}_{t}} + R_{k+1} \leq B \;\;\forall \mb{n}_{[k]} \in \mathcal{N}_{[k]}.
	\end{eqnarray*}
	Since this constraint must hold for all \(\mb{u}_{k}\) with \(\|\mb{u}_{k}\|_2 \leq r_{k}\), taking the maximum over \(\mb{u}_{k}\) and substituting \(\|\mb{L}_{k}^{\top}\mb{y}_{k}\|_2 = \sqrt{\mb{y}_{k}^{\top}\bs{\Sigma}_{k}\mb{y}_{k}}\), we obtain
	\begin{equation}
		\begin{aligned}
			&s_{k-1} + \Theta_{k}  + r_{k} \sqrt{\mb{y}_{k}^{\top}\bs{\Sigma}_{k}\mb{y}_{k}}\\
			&+ \sum_{t=k+1}^{T} r_{t}\sqrt{A_{k+1,t}a_k \mb{y}_{t}^{\top}\bs{\Sigma}_{k}\mb{y}_{t}}\\
			&+ \sum_{t=k+1}^{T}r_{t}\sqrt{A_{k+1,t}\mb{y}_{t}^{\top}\mb{C}_k \mb{y}_{t}} + R_{k+1} \leq B \;\;\forall \mb{n}_{[k]} \in \mathcal{N}_{[k]}. 
		\end{aligned}
		\label{eq:elips7}
	\end{equation}
	Using \(A_{k+1,t} a_k = \prod_{j=k+1}^{t}a_j a_k = A_{k,t}\), and the fact that $A_{k,k} = 1$
	we can rewrite~(3) as 
	\begin{align*}
		s_{k-1} + \Theta_{k}  &+ \sum_{t=k}^{T} r_{t}\sqrt{A_{k,t}\; \mb{y}_{ t}^{\top}\bs{\Sigma}_{k}\mb{y}_{t}}\\
		&\hspace{-0.4cm}+ \sum_{t=k+1}^{T}r_{t}\sqrt{A_{k+1,t}\mb{y}_{t}^{\top}\mb{C}_k \mb{y}_{t}} + R_{k+1} \leq B \;\;\forall \mb{n}_{[k]} \in \mathcal{N}_{[k]},
	\end{align*}
	by extending the first summation.
	Using the definition of \(R_{k}\), 
	we obtain the desired result,
	\begin{align*}
		&s_{k-1} + \Theta_{k}  + \sum_{t=k}^{T} r_{t}\sqrt{A_{k,t}\; \mb{y}_{t}^{\top}\bs{\Sigma}_{k}\mb{y}_{t}} + R_{k} \leq B \;\;\forall \mb{n}_{[k]} \in \mathcal{N}_{[k]}.
	\end{align*}
	The complete reformulation follows immediately by substituting \(k=1\). 
	\Halmos
	\endproof

	%%%%%%%%%%%%%%%%%%%%%%%%%%%%%%%%%%%%%%%%%%%
	
	\section{Robust Counterpart of~(C-DRO)} 
	\label{sec:rc_c_dro}
	The following proposition better illustrates the connection between Constraint~(C-DRO) and its robust counterpart~(6).
	Given the uncertainty sets \(\widetilde{\mU}_1,\dots,\widetilde{\mU}_T^{T-1}\), their joint uncertainty set is defined by
	\begin{equation*}
		\widetilde{\mU} = \left\{P \mid P = P_1 \times \dots \times P_{T \mid T-1},\; P_1 \in \widetilde{\mU}_1,\;P_{t \mid t-1} \in \widetilde{\mU}_{t}(\mb{d}_{t-1}) \;\forall \mb{d}_{t-1} \in \Xi_{t-1} \;\forall t\right\}.
	\end{equation*}
	The joint set, \(\widetilde{\mU}\) is the set of all distributions \(P\) with the marginals lying in the specified uncertainty sets $\widetilde{\mU}_{t}(\mb{d}_{t-1})$.
	\begin{proposition}
		Given the sets \(\widetilde{\mU}_1,\dots,\widetilde{\mU}_T(\mb{d}_{T-1})\) and their joint uncertainty set \(\widetilde{\mU}\), the robust counterpart of Constraint~(C-DRO), given by 
		\begin{equation}
			\label{eq:c1rc}
			\sup_{P \in \widetilde{\mU}} \mathbb{E}_P\left[\sum_{t=1}^{T}h_t(\mb{x}_t,\mb{d}_t)\right] \leq B
		\end{equation}
		is equivalent to 
		\begin{align*}
			&\sup_{P_1 \in \widetilde{\mU}_1} \mathbb{E}_{P_1}\!\bigg[h_1(\mb{x}_1,\mb{d}_1) +
			\!\!\!\!\sup_{P_{2 \mid 1} \in \widetilde{\mU}_2(\mb{d}_1)}\!\!\bigg\{\mathbb{E}_{P_{2 \mid 1}}\bigg[h_2(\mb{x}_2,\mb{d}_2) + \dots \\
			&\hspace{4cm}+\sup_{P_{T \mid T-1} \in   \widetilde{\mU}_T(\mb{d}_{T-1})}\!\!\!\!\!\!\!\!\!\left\{\mathbb{E}_{P_{T\mid{T-1}}}
			\left[h_T(\mb{x}_T,\mb{d}_T)\right]\right\}\bigg]\bigg\}\bigg] \leq B.
		\end{align*}	
	\end{proposition} 
	%
	%%%%%%%%%%%%%%%%%%%%%%%%%%%%%%%%%%%%%%%%%%%%%%%%%%%%%%%%%%%%%%%%%%%%
	\proof{Proof:}
	%%%%%%%%%%%%%%%%%%%%%%%%%%%%%%%%%
	We first show the forward direction of the equivalence.
	Observe that for any distribution \(P\in\widetilde{\mU}\) with the corresponding marginal distribution \(P_1\), and for each \(\mb{d}_{t-1}\) with the conditional distribution \(P_{t \mid t-1}\), it holds that
	\begin{equation}
		\begin{aligned}
			\label{eq:prop1exp}
			&\mathbb{E}_{P_1}\left[h_1(\mb{x}_1,\mb{d}_1) +
			\left\{\mathbb{E}_{P_{2 \mid 1}}\left[h_2(\mb{x}_2,\mb{d}_2) + \dots + 
			\left\{\mathbb{E}_{P_{T\mid{T-1}}}
			\left[h_T(\mb{x}_T,\mb{d}_T)\right]\right\}\right]\right\}\right] \\
			=\;&\mathbb{E}_P\left[\sum_{t=1}^{T}h_t(\mb{x}_t,\mb{d}_t)\right]. 
		\end{aligned}
	\end{equation}
	For a small \(\epsilon > 0\), let \(P^* \in \widetilde{\mU}\) be such that
	\begin{equation}
		\label{eq:c1rc2}
		\mathbb{E}_{P^*}\left[\sum_{t=1}^{T}h_t(\mb{x}_t,\mb{d}_t)\right] \geq \sup_{P \in \widetilde{\mU}} \mathbb{E}_P\left[\sum_{t=1}^{T}h_t(\mb{x}_t,\mb{d}_t)\right] - \epsilon.
	\end{equation}
	That means the LHS of~(6) is $\epsilon$-optimal to the LHS of~(4).
	Since \(P^*\) is in \(\widetilde{\mU}\), there exist marginal distributions \(P_1^*,\dots,P_{T|T-1}^*\) such that \(P_1^*\) lies in \(\widetilde{\mU}_1\) and the conditional distribution of \(\mb{d}_{t}\), \(P_{t \mid t-1}^*\)  lies in \(\widetilde{\mU}(\mb{d}_{t-1})\), where \(\mb{d}_{t-1} \in \Xi_{t-1}\). This holds true for all \(t=2,\dots,T\).
	Using~(5), this means
	\begin{align*}
		&\sup_{P_1 \in \widetilde{\mU}_1} \mathbb{E}_{P_1}\!\bigg[h_1(\mb{x}_1,\mb{d}_1) +
		\!\!\!\!\sup_{P_{2 \mid 1} \in \widetilde{\mU}_2(\mb{d}_1)}\!\!\bigg\{\mathbb{E}_{P_{2 \mid 1}}\bigg[h_2(\mb{x}_2,\mb{d}_2) + \dots \\
		& \hspace{5cm}+ 
		\!\!\!\!\!\!\!\sup_{P_{T \mid T-1} \in   \widetilde{\mU}_T(\mb{d}_{T-1})}\!\!\!\!\!\!\left\{\mathbb{E}_{P_{T\mid{T-1}}}
		\left[h_T(\mb{x}_T,\mb{d}_T)\right]\right\}\bigg]\bigg\}\bigg]  \\
		& \geq \mathbb{E}_{P_1^*}\left[h_1(\mb{x}_1,\mb{d}_1) +
		\left\{\mathbb{E}_{P_{2 \mid 1}^*}\left[h_2(\mb{x}_2,\mb{d}_2) + \dots + 
		\left\{\mathbb{E}_{P_{T\mid{T-1}}^*}
		\left[h_T(\mb{x}_T,\mb{d}_T)\right]\right\}\right]\right\}\right] \\
		&= \mathbb{E}_{P^*}[\sum_{t=1}^{T}h_t(\mb{x}_t,\mb{d}_t)]\\
		& \geq  \sup_{P \in \widetilde{\mU}} \mathbb{E}_P[\sum_{t=1}^{T}h_t(\mb{x}_t,\mb{d}_t)] - \epsilon.
	\end{align*}
	
	Now for the opposite side of the inequality, let \(P_1^*,P_{2 \mid 1}^*\) up to \(P_{T \mid T-1}^*\) be \(\epsilon\)-optimal to~(6).
	The constraints of~(6) ensure that \(P_1^* \in \widetilde{\mU}_1\) and \(P_{t \mid t-1}^* \in \widetilde{\mU}_t(\mb{d}_{t-1})\).
	Since \(\widetilde{\mU}\) is the set of all joint distributions with these specified marginals, for the joint distribution \(P^*\) with these marginals, the equation~(5) will hold.
	However since \(P_1^*\) and \(P_{t \mid t-1}^*\) are \(\epsilon\)-optimal, this means that
	\begin{align*}
		&\sup_{P_1 \in \widetilde{\mU}_1} \mathbb{E}_{P_1}\!\bigg[h_1(\mb{x}_1,\mb{d}_1) +
		\!\!\!\!\sup_{P_{2 \mid 1} \in \widetilde{\mU}_2(\mb{d}_1)}\!\!\bigg\{\mathbb{E}_{P_{2 \mid 1}}\bigg[h_2(\mb{x}_2,\mb{d}_2) + \dots \\
		&\hspace{5cm}+\sup_{P_{T \mid T-1} \in   \widetilde{\mU}_T(\mb{d}_{T-1})}\bigg\{\mathbb{E}_{P_{T\mid{T-1}}}
		\left[h_T(\mb{x}_T,\mb{d}_T)\right]\bigg\}\bigg]\bigg\}\bigg]  \\
		& \leq \mathbb{E}_{P_1^*}\!\!\left[h_1(\mb{x}_1,\mb{d}_1) +
		\left\{\mathbb{E}_{P_{2 \mid 1}^*}\!\!\left[h_2(\mb{x}_2,\mb{d}_2) + \dots + 
		\left\{\mathbb{E}_{P_{T\mid{T-1}}^*}
		\left[h_T(\mb{x}_T,\mb{d}_T)\right]\right\}\right]\right\}\right] + \epsilon \\
		&= \mathbb{E}_{P^*}\left[\sum_{t=1}^{T}h_t(\mb{x}_t,\mb{d}_t)\right] + \epsilon\\
		& \leq  \sup_{P \in \widetilde{\mU}} \mathbb{E}_P[\sum_{t=1}^{T}h_t(\mb{x}_t,\mb{d}_t)] + \epsilon.
	\end{align*}
	This gives the opposite inequality, and the result follows by letting $\epsilon$ shrink to zero. \Halmos
	%%%%%%%%%%%%
	\endproof  
	%%%%%%%%%%%%%%%%%%%%%%%%%%%%%%$
	%	
	Note that this result does not depend on the structure of the moment based uncertainty set and can be extended to sets \(\widetilde{\mU}_1,\dots,\widetilde{\mU}_T(\mb{d}_{T-1})\) of any structure, in which the parameters depend on previous realizations.
	An important part of this result is the additive nature of the constraint.

	%%%%%%%%%%%%%%%%%%%%%
	\section{RO Application: Knapsack Problem with Negative Correlations}
	\label{sec:neg_cor}
\begin{figure}[t]
		\centering
		\includegraphics[width=75mm]{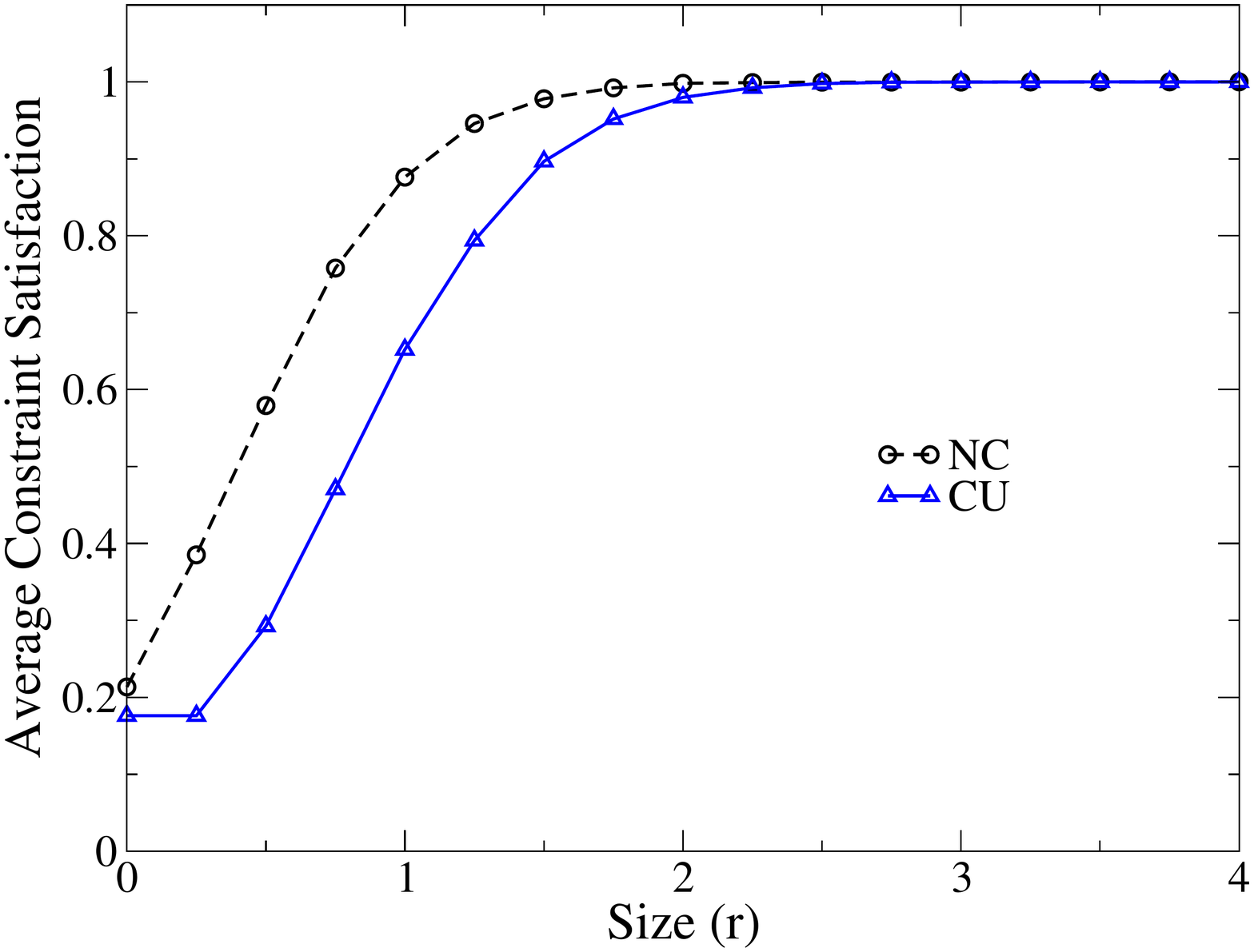}\hspace{1mm}
		\includegraphics[width=75mm]{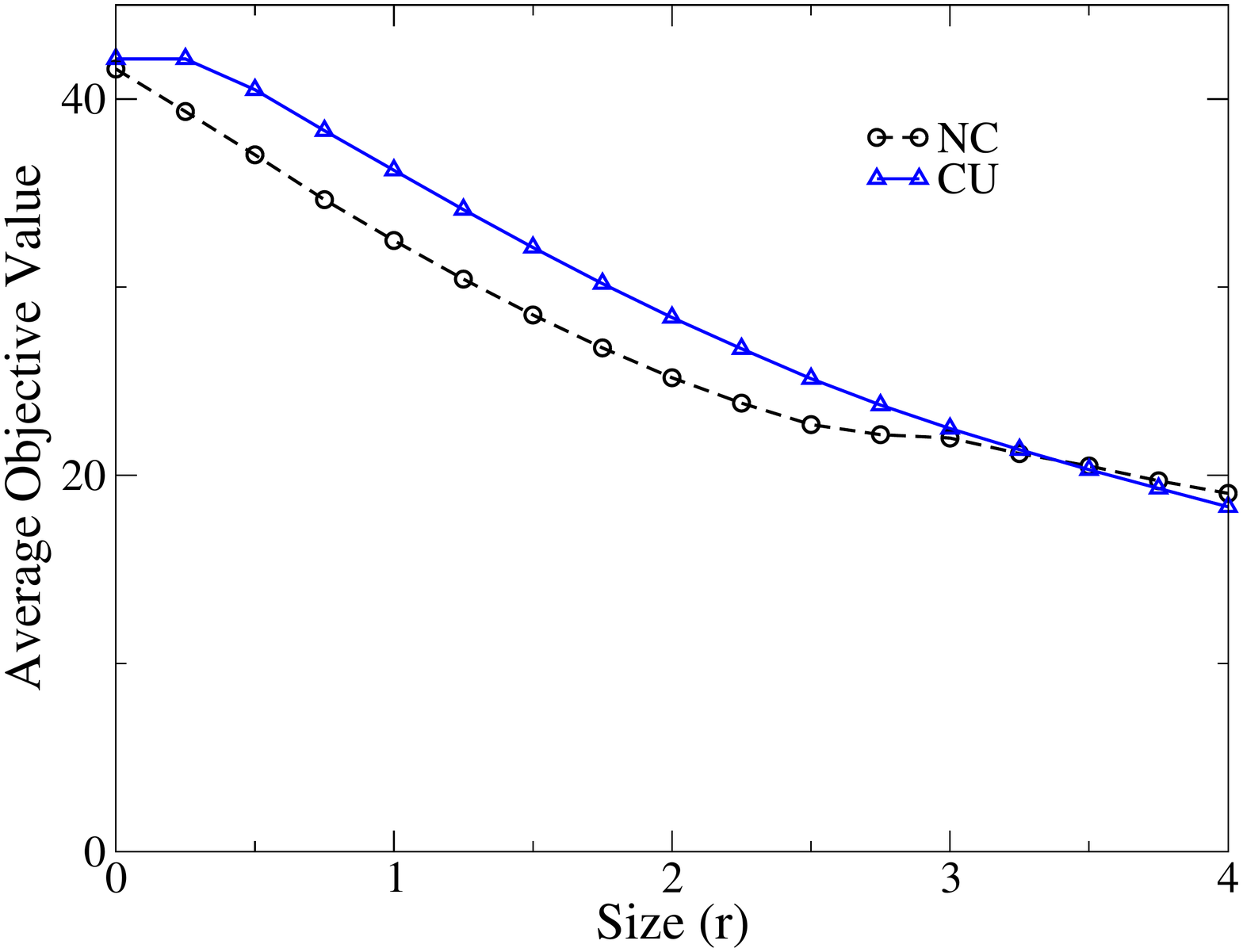}
		\caption{\label{fig:ts_n02_cs_ao}Comparison of Connected and Nonconnected sets for the Robust Knapsack Problem at different set sizes:
			Average Fraction of Constraint Satisfaction (Left) Average Objective Value (Right)}
	\end{figure}

\begin{figure}[t]
		\centering
		\includegraphics[width=75mm]{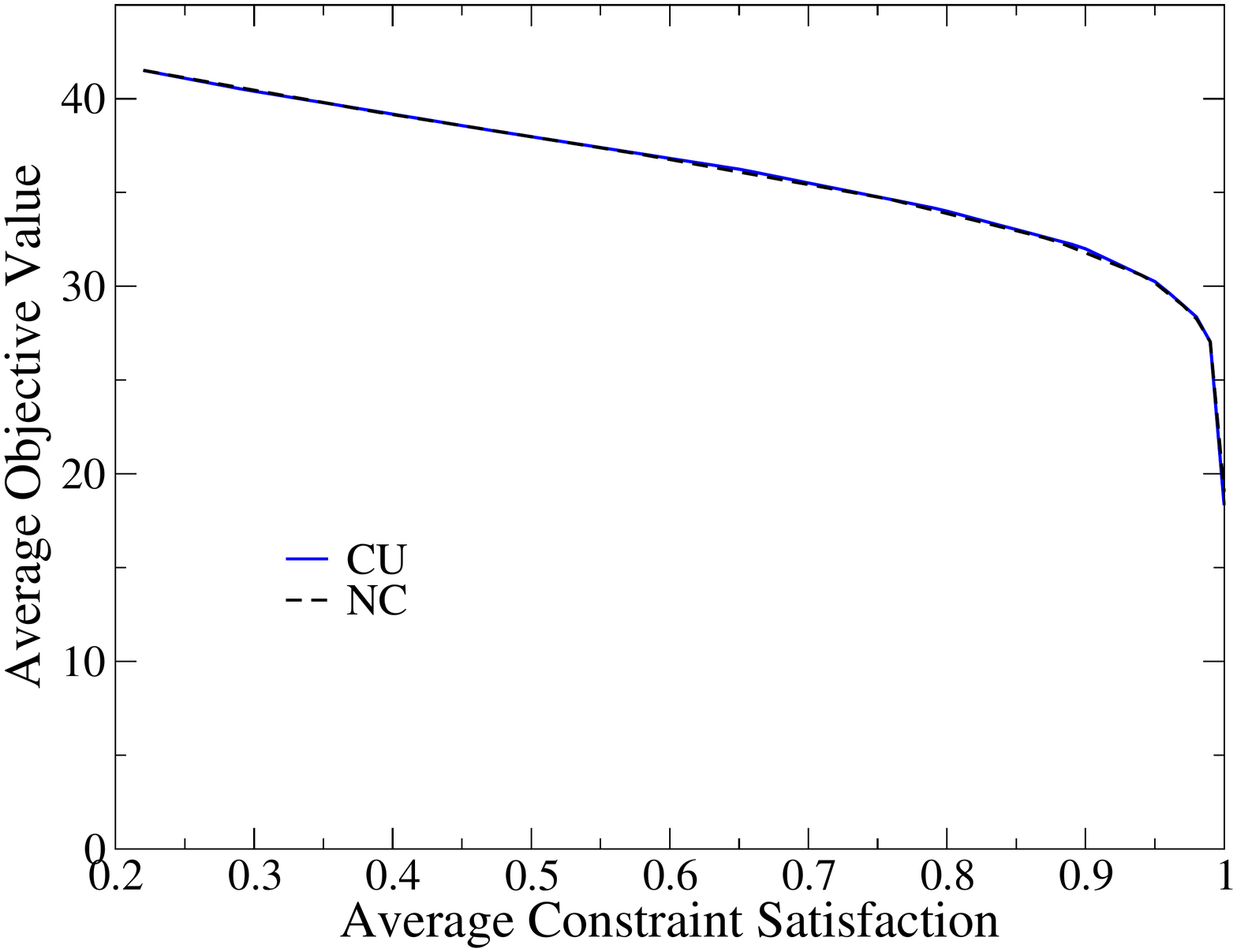}\hspace{1mm}
		\includegraphics[width=75mm]{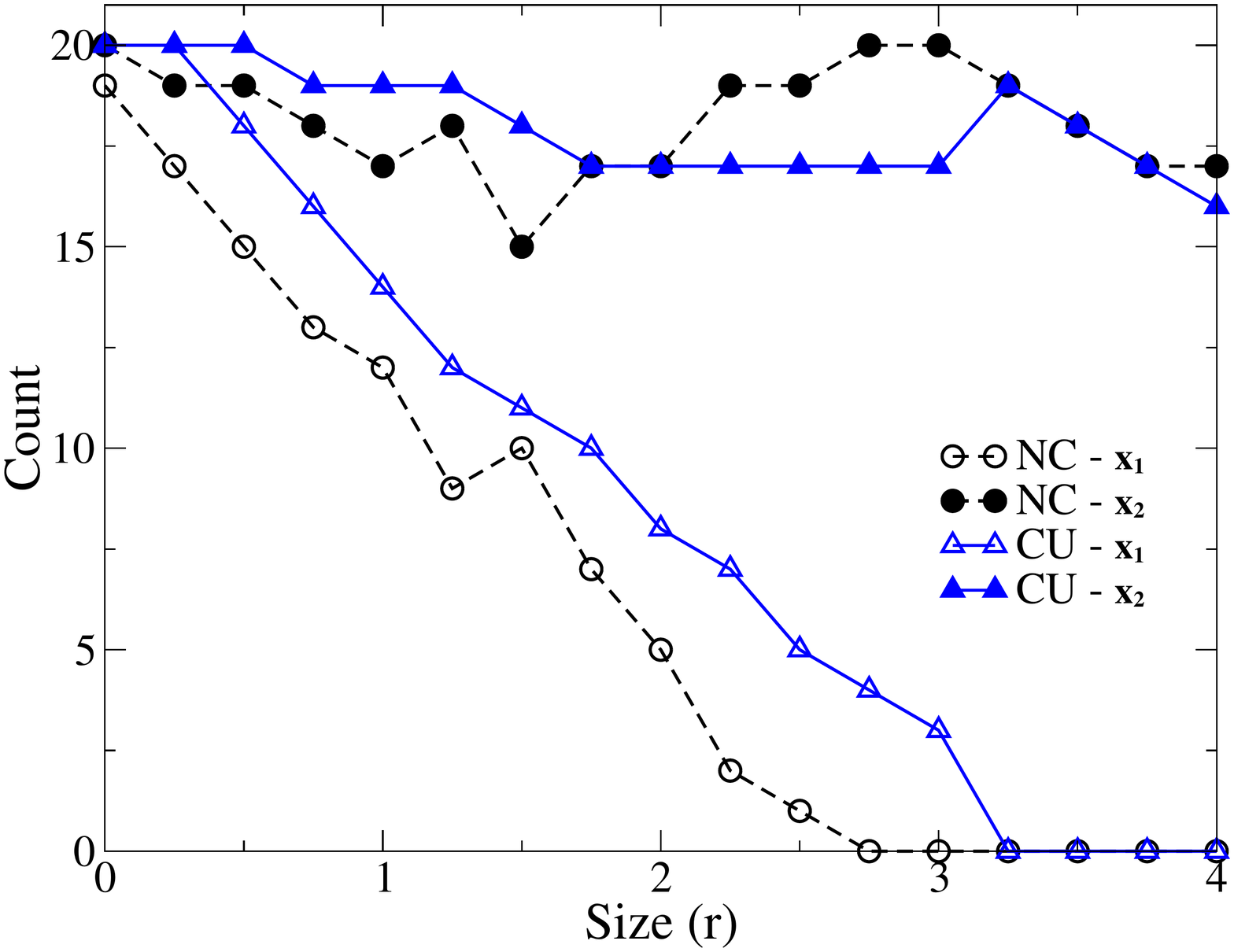}
		\caption{\label{fig:ts_n02_aocs_sn}Comparison of Connected and Nonconnected sets:
			Average Objective Value vs. Constraint Satisfaction (Left) and For a single iteration, the number of non-zero variables of a period (Right)}
	\end{figure}

	We repeat the experiments conducted for the Knapsack problem with a negative correlation among the weights over time--that is, \(\bs{\Psi} = -0.2 \times \mb{I}\).
	Figure~\ref{fig:ts_n02_cs_ao} presents the value of average objective for any level of constraint satisfaction (left) and the number of non-zero allocations vs. set size.
	Figure~\ref{fig:ts_n02_aocs_sn} displays how the average constraint satisfaction (left) and the average objective value (right) change with size of the uncertainty set \(r\).  
	These results show that when uncertainties are negatively correlated, CU solutions achieve a better objective value but a lower constraint satisfaction than NC solutions. 
	However, for any level of constraint satisfaction both have similar performance. 
	
	\begin{itemize}
		\item\textbf{\emph{Effect of Uncertainty Set Size}}: 
		For both models, constraint satisfaction increases and average objective value decreases with $r$.
		Note that the objective value is only measured, if constraints are satisfied.
		\item\textbf{\emph{CU vs. NC}}: For any $r$, CU solutions have lower constraint satisfaction than NC solutions.  
		CU solutions also have higher average objective value. 
		This is because connected sets depend on the first period, and the negative correlation instead of magnifying the worst case causes one of the realizations to take a lower value. 
		This reduces the protection but increases the objective value. 
		
		Note that for any level of satisfaction, the average objective of CU is almost the same as NC (see Figure~\ref{fig:ts_n02_aocs_sn} left).
		\item\textbf{\emph{First vs Second period solutions}}:
		For a single estimate in (ii) shown in Figure~\ref{fig:ts_n02_aocs_sn} (right), the optimal solution gradually concentrates only on \(\mb{x}_2\) for both CU and NC as $r$ increases.
		For NC, this is because $\mb{c}_2$ tends to be higher. 
		For CU, the negative correlation prevents the magnification of the weights for both periods as such second period allocations which have higher objective coefficients can be selected. 
	\end{itemize}
	
	
	\section{DRO Application: Portfolio Optimization}
	\label{sec:dro_app}
	We study a practical portfolio optimization problem on historical stock data. 
	For our portfolio, we choose among 5 stocks.
	The experiment is repeated 150 times for randomly selected dates.
	At each date, we compute the weekly returns for the previous 100 weeks based on stock price data and fit them to a time series model. 
	The two weeks following the selected dates serve to evaluate the performance of the model. 
	
	To capture risk aversion, in each experiment, we maximize a concave piecewise linear utility function \(u= \min(1.5 r, 0.015 + r, 0.06 + 0.1 r)\), where \(r\) denotes the return on the portfolio.
	It is assumed that \(\overline{\bs{\mu}}_{2}  = \bs{\mu}_{2}(\mb{d}_1) + \bs{\delta}\) and \(\underline{\bs{\mu}}_{2}  = \bs{\mu}_{2}(\mb{d}_1) - \bs{\delta} \) with \(\bs{\mu}_{2}(\mb{d}_1) = \bs{\mu}_0 + \mb{A}\bs{\mu}_1 + \mb{B}\mb{d}_1\).
	The vector \(\bs{\mu}_0\) and the matrices \(\mb{A}\) and \(\mb{B}\) are estimated from data via a vector autoregressive moving average with lag 1.
	The parameters \(\bs{\Sigma}_1\) and \(\bs{\Sigma}_2\) are set to the residual covariance matrices.
	The value of \(\bs{\mu}_1\) is the return point-estimate at the end of the first week and \(\bs{\delta}\) is three times the standard error of this estimate (in order to cover almost all realized means under normality assumption).
	
	To probe the performance of the proposed CU model, we compare it against the standard DRO and a modified DRO, described as following:
	\begin{description}\setlength{\itemindent}{5mm}
		\item CU: Model with connected uncertainty set,
		\item DRO-1: Model with \(\bs{\mu}_2 = \bs{\mu}_1\), and
		\item DRO-2:  Model with \(\bs{\mu}_2 = \mb{A}\bs{\mu}_1 + \mb{B}\bs{\mu}_1\).
	\end{description}
	DRO-1 represents the standard model as used widely in the literature. 
	It is computationally attractive because of a simpler reformulation. 
	However, it does not take into account potential connections to previous periods.
	We propose a modified version (DRO-2), which captures some of the potential connections to previous periods using the unconditional mean for the second period.
	
	In these experiments, the parameters \(\bs{\mu}_i\) and \(\bs{\Sigma}_i\) are defined based on past weeks, and the returns are computed over the future two weeks.
	We evaluate the solution quality from these three models by comparing the returns.
	Each experiment starts with an initial wealth of \(W_0=\mathdollar 100\), which is recomputed using the total return on the portfolio
	\begin{equation*}
		W_{t+1} = W_t \cdot (\widehat{\mb{d}}_t^\top\mb{x}_t), \; t=1,2,
	\end{equation*}
	The realized demand $\widehat{\mb{d}}_t$ is taken from the actual stock data.
	This wealth (W) is then averaged over the different experiments and reported in Table~\ref{table:portfolio_wealth}, along with the standard deviation (Std) as a measure for robustness of the solution.
	
	\begin{table}[h]
		\centering
		\begin{tabular}{|c|c|c|c|c|c|c|}
			\hline
			Model & \multicolumn{2}{c|}{Period 1} & \multicolumn{2}{c|}{Period 2} & \multicolumn{2}{c|}{Period 3} \\
			& W & Std & W & Std & W & Std \\ \hline
			CU & 100 & 0 & 100.1 & 1.97 & 101.2 & 5.38 \\ 
			\hline
			DRO1 & 100 & 0 & 99.9 & 2.36 & 100.9 & 5.92 \\ 
			\hline 
			DRO2 & 100 & 0 & 99.9 & 2.36 & 100.9 & 6.06 \\ 
			\hline 
		\end{tabular}
		\caption{\label{table:portfolio_wealth}Average Wealth (W) and its Standard Deviation (Std) over time for various models}
	\end{table} 
	
	Table~\ref{table:portfolio_wealth} displays the performances of the three different models.
	For all three, we observe a positive average return at the end of two weeks, which is attributed to the mean positive return on the random samples for the five stocks.
	Furthermore, the standard deviation (spread) of the sample paths grows as time elapses, because the inherent uncertainty in returns is compounding over time.
	
	When comparing the three models, we observe marginal differences in the average wealth, reflecting the relatively short time horizon.
	However, we observe the wealth standard deviation for the CU model to be lower than both of the DRO models.
	This is attributed to the fact that the CU model captures the compounding worst-case effect of connected periods and yields a more conservative solution.
	the CU model reduces the wealth standard deviation for wealth by 0.417 (17.6\% of the original wealth standard deviation) for period 1 and 0.563 (9.5\%) for period 2 compared to DRO-1.
	Other measures of deviation, such as interquartile range, reveal a similar decrease in wealth.
	Therefore, the CU model is able to select assets that are less volatile in order to provide more robust allocations.

	Note that solving the CU model is computationally more demanding than either of the DRO models, because of latter's convex subproblems.
	Therefore, depending on the application, the advantages of the CU model have to outweigh the additional computational burden.
	Furthermore, the DRO-2 model represents another way to leverage the connection between uncertainties of different periods.
	However, it only accounts for the average effect between periods and not the worst case. 
	
	\section{DRO Application: Portfolio Optimization II}
	\label{sec:dro_app_2}
\begin{figure}
		\centering		
		\includegraphics[width=5.6cm,page=9]{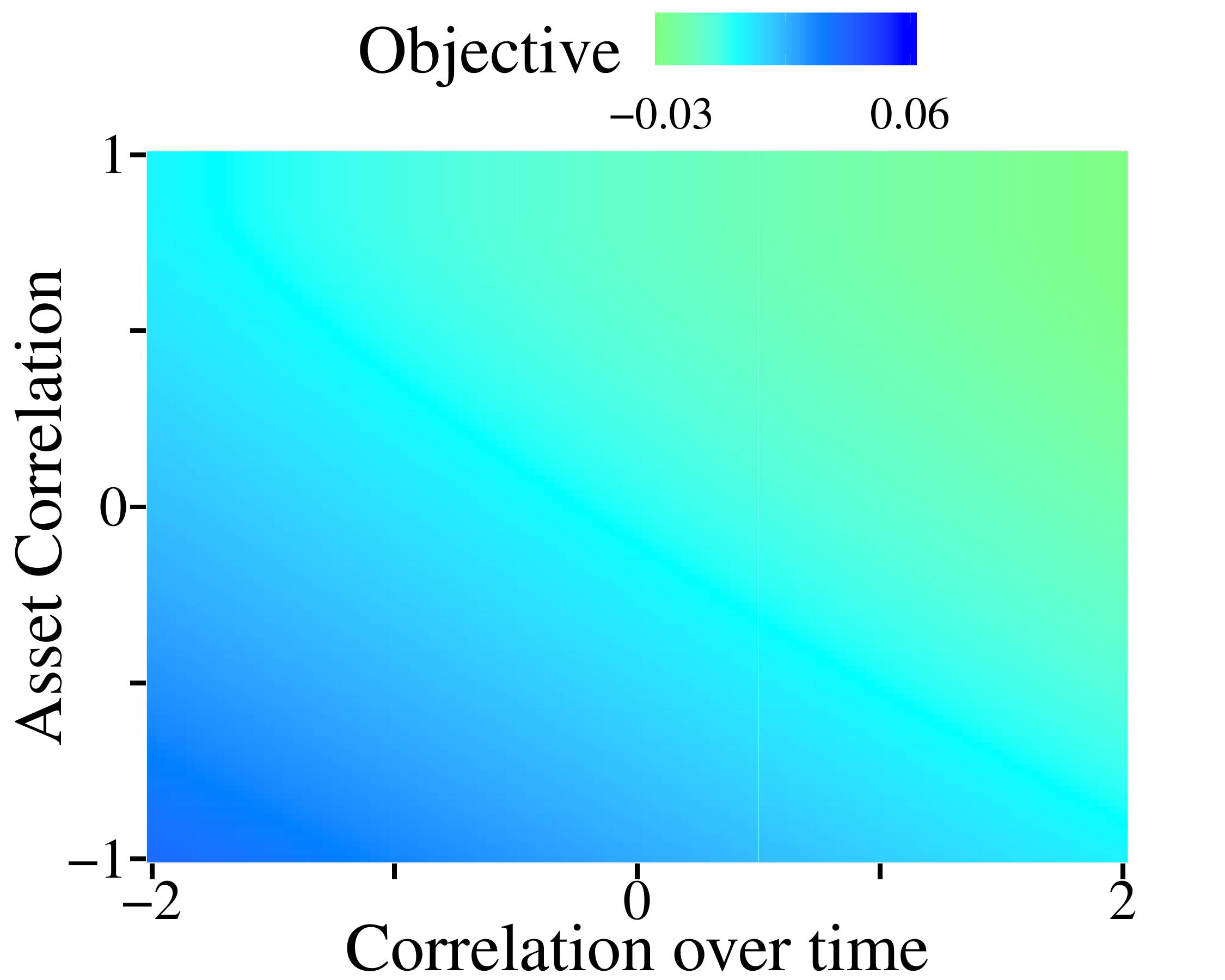}
		\includegraphics[width=5.6cm,page=11]{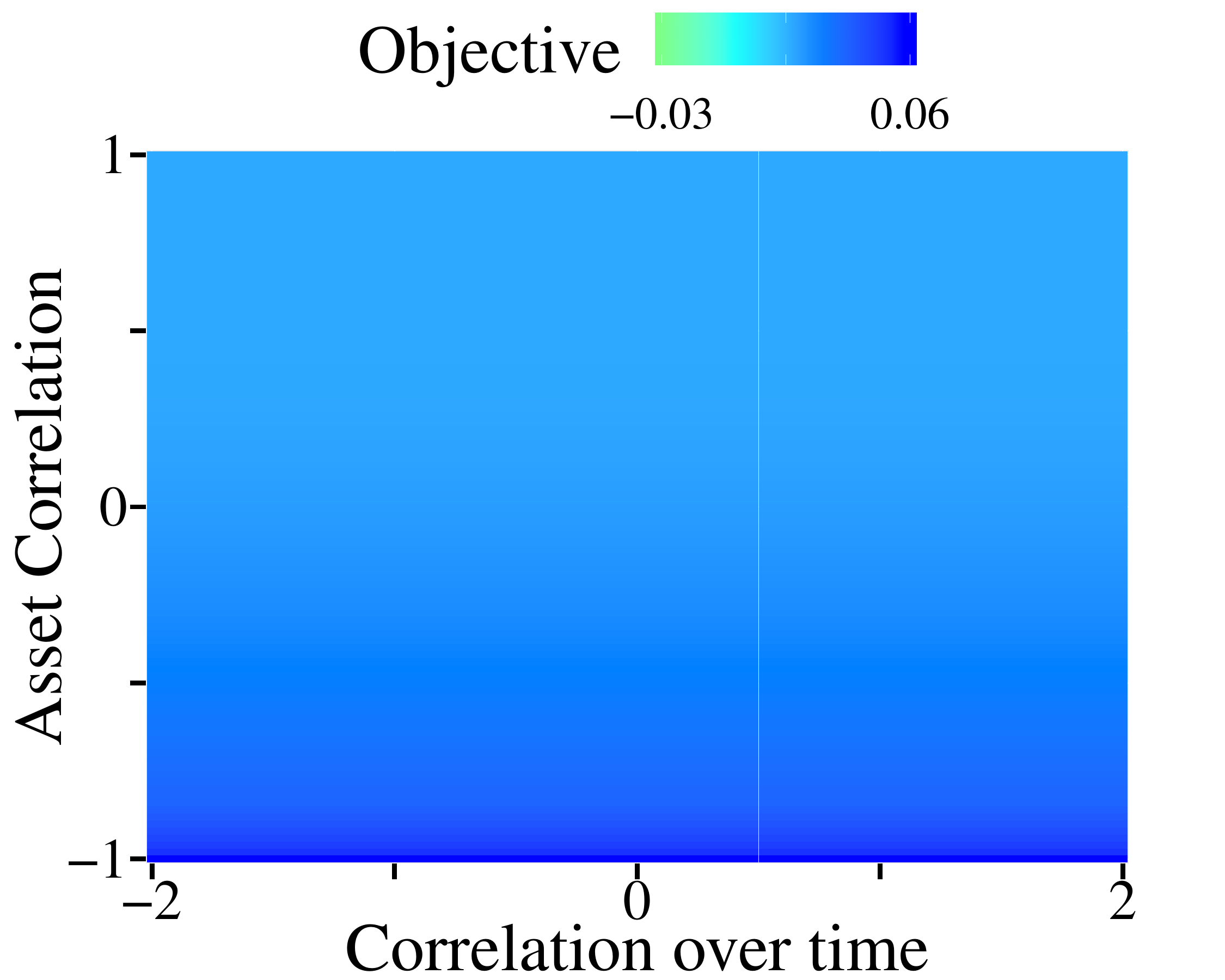}
		\caption{\label{fig:cu_objective}CU and DRO Problem Objective}
	\end{figure}

	\begin{enumerate}
		\item \emph{Objective function:} to measure performance of the model in the nominal setting (Figure~\ref{fig:cu_objective}).
		\item \emph{Worst case realized wealth:} to study the effect of adversarial settings (Figure~\ref{fig:worst_wealth_3}).
		\item \emph{Standard deviation of wealth:} to evaluate the robustness of solutions (Figure~\ref{fig:wealth_sd_3}).
	\end{enumerate}
	
	\subsection*{Objective}
	\subsubsection*{Asset Correlation}
	Figure~\ref{fig:cu_objective} indicates that positive correlations among the assets reduce the objective performance for the CU model. 
	This is due to the fact that positive correlation worsens the worst-case returns everywhere simultaneously, whereas a negative correlation prevents this from occurring as a result the objective (which is worst-case by definition) reduces.
	
	\subsubsection*{Correlation over time}
	Figure~\ref{fig:cu_objective} also presents similar behavior for correlation over time leading to worse objective performance for the same reasons as asset correlation. 
	The DRO model does not capture the correlation over time but it also performs better when the assets are negatively correlated. 
	
	\subsection*{Standard Deviation and Worst Case Wealth}
	\begin{figure}
		\centering		
		\includegraphics[width=5.6cm,page=22]{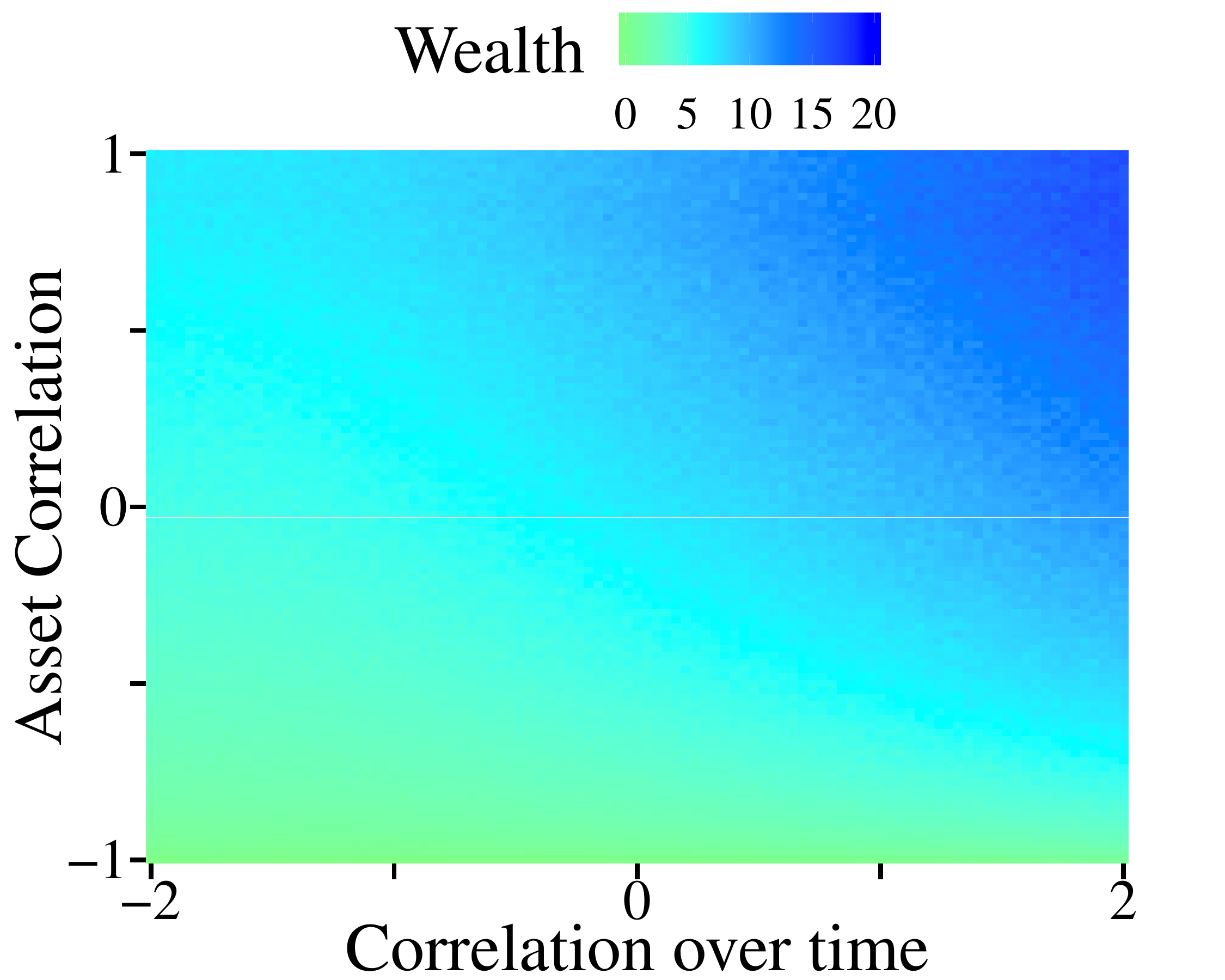}
		\includegraphics[width=5.6cm,page=24]{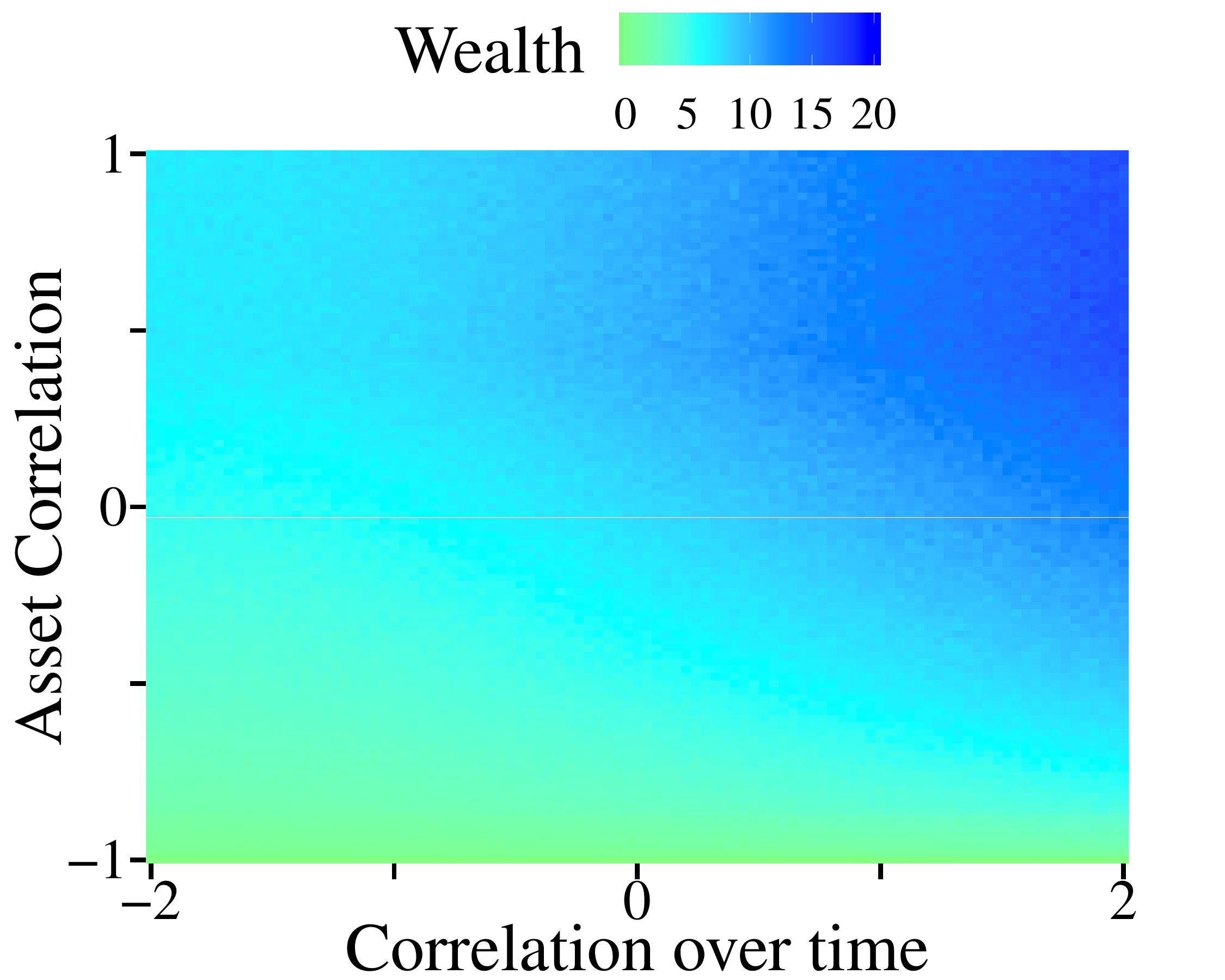}
		\caption{\label{fig:wealth_sd_3}Realized Wealth Standard Deviation for CU and DRO at end of Period 2}
	\end{figure}
	Figure~\ref{fig:wealth_sd_3} shows the standard deviation for the realized wealth at the end of the second period. 
	\subsubsection*{Asset Correlation}
	The standard deviation is higher, when the assets are positively correlated among themselves and over time. 
	
	\begin{figure}		
		\centering
		\includegraphics[width=5.6cm,page=28]{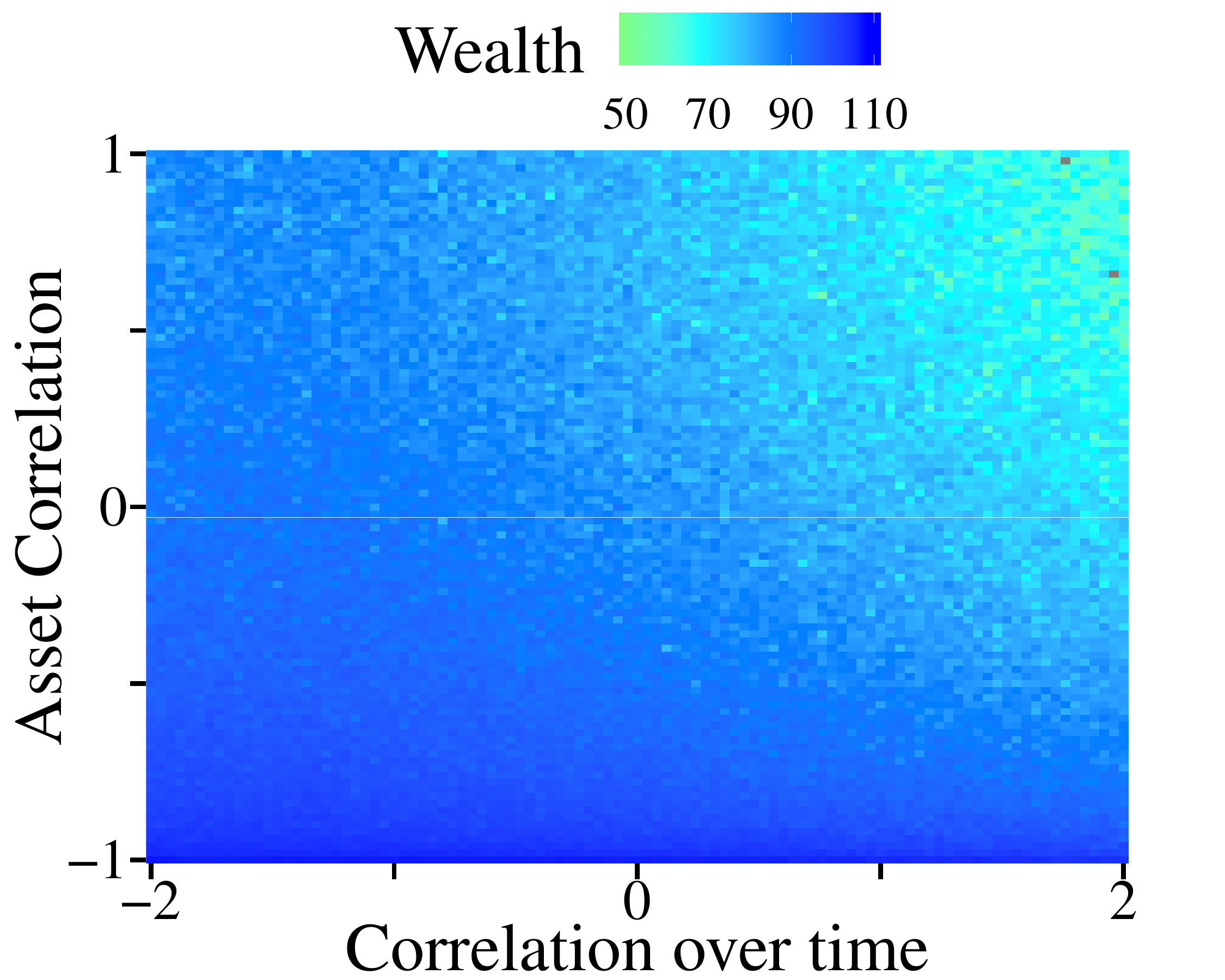}
		\includegraphics[width=5.6cm,page=30]{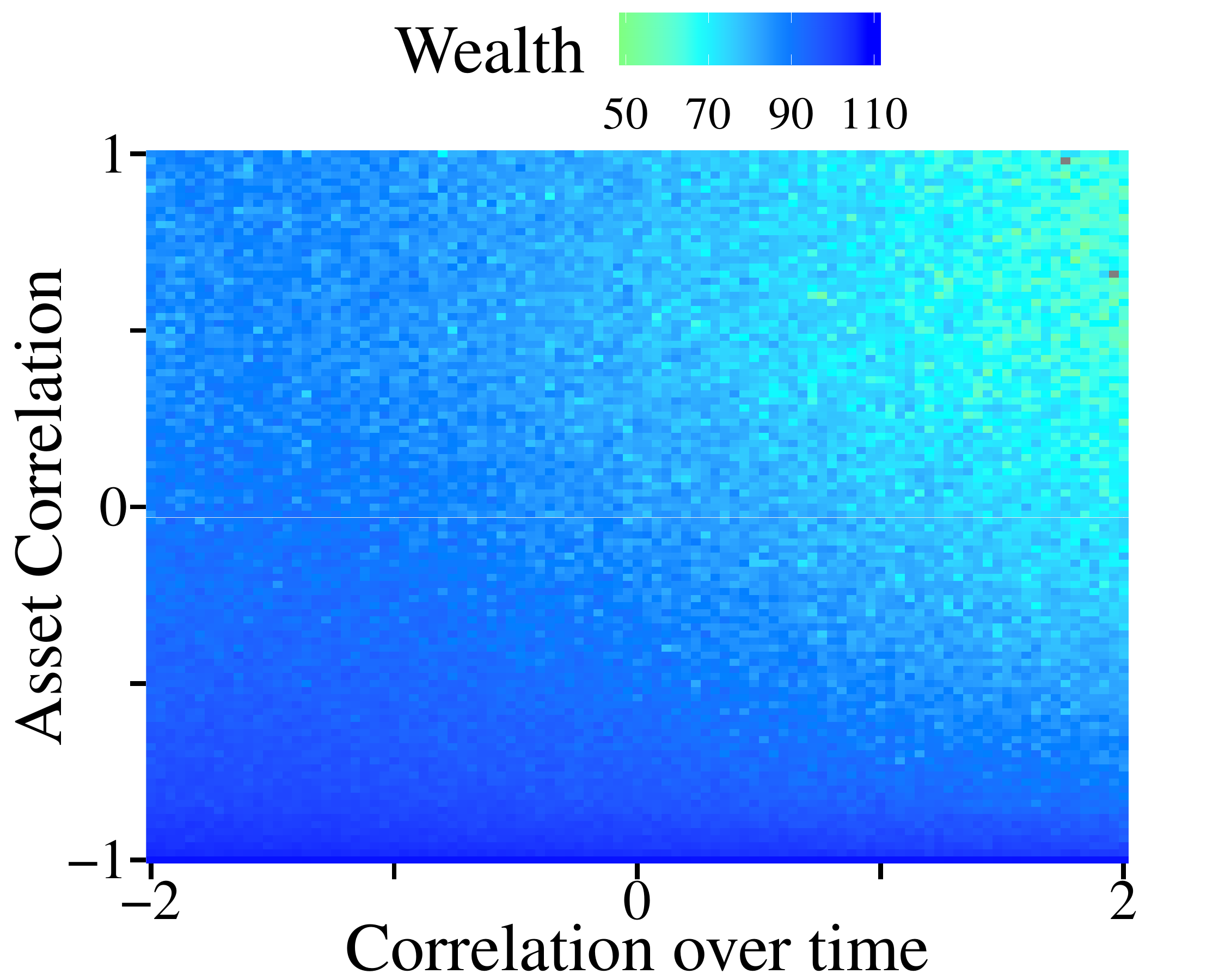}
		\caption{\label{fig:worst_wealth_3}Worst Case Realized Wealth for Period 2}
	\end{figure}
	Figure~\ref{fig:worst_wealth_3} shows the worst case realized wealth at the end of the second period.
	The worst case wealth is higher, when assets are negatively correlated among themselves and over time. 
	This is because the negative correlation makes it unlikely for both returns to be very low. 
	The pattern observed here is similar to that in Figure~\ref{fig:cu_objective}, with the worst case wealth being worse when assets are positively correlated. 
	
\bibliographystyle{informs2014} 
\bibliography{DyUS_short}

	%%%%%%%%%%%%%%%%%